%&amstex
% This is an AMS-TeX file and should be compiled
% using AMS-TeX.
% the command should be:  amstex file  or something like that
\magnification=1200

\loadmsam
\loadmsbm
\loadeufm
\loadeusm
\UseAMSsymbols
%\usepackage{amsby}
%\input amssym.def
%\hsize=6.00 true in
%\hoffset=.50 true in
%\voffset=-0.1 true in
%\vsize=8.75 true in

\font\BIGtitle=cmr10 scaled\magstep3
\font\bigtitle=cmr10 scaled\magstep1
\font\boldsectionfont=cmb10 scaled\magstep1
\font\section=cmsy10 scaled\magstep1

\def\scr#1{{\fam\eusmfam\relax#1}}
\def\scrA{{\scr A}}
\def\scrB{{\scr B}}
\def\scrC{{\scr C}}
\def\scrD{{\scr D}}

\def\scrF{{\scr F}}
\def\scrG{{\scr G}}
\def\scrH{{\scr H}}
\def\scrI{{\scr I}}

\def\scrK{{\scr K}}
\def\scrJ{{\scr J}}

\def\scrN{{\scr N}}

\def\scrS{{\scr S}}

\def\scrR{{\scr R}}

\def\scrT{{\scr T}}
\def\scrV{{\scr V}}

\def\scrW{{\scr W}}
\def\gr#1{{\fam\eufmfam\relax#1}}

\def\grC{{\gr C}}   
\def\grD{{\gr D}}

   \def\grg{{\gr g}}
   
\def\grI{{\gr I}}

\def\grL{{\gr L}}

\def\grM{{\gr M}}   
   \def\grn{{\gr n}}
   
\def\grP{{\gr P}}   
   
\def\grR{{\gr R}}   
\def\grS{{\gr S}}

\def\db#1{{\fam\msbfam\relax#1}}

\def\dbA{{\db A}} 
\def\dbC{{\db C}} 
 \def\dbF{{\db F}}
\def\dbG{{\db G}} \def\dbH{{\db H}}

 \def\dbN{{\db N}}
 
\def\dbQ{{\db Q}} \def\dbR{{\db R}}
\def\dbS{{\db S}}

 \def\dbZ{{\db Z}}

\def\eps{{\varepsilon}}

\def\Gtil{\widetilde{G}}

\def\dbZhat{\widehat{\dbZ}}

\def\Ker{\text{Ker}}
\def\der{\text{der}}
\def\Sh{\hbox{\rm Sh}}

\def\sc{\text{sc}}
\def\Res{\text{Res}}
\def\ab{\text{ab}}
\def\ad{\text{ad}}
\def\Ad{\text{Ad}}
\def\Gal{\text{Gal}}
\def\Hom{\text{Hom}}
\def\End{\text{End}}
\def\Spec{\text{Spec}}

\def\Lie{\text{Lie}}

\def\leaderfill{\leaders\hbox to 1em
     {\hss.\hss}\hfill}
\def\nspace{\lineskip=1pt\baselineskip=12pt\lineskiplimit=0pt}

     %the way to use this is "\Proclaim{Theorem 1.1.}" for instance.
\def\finishproclaim{\par\rm
     \ifdim\lastskip<\medskipamount\removelastskip
     \penalty55\medskip\fi}
\def\endproof{$\hfill \square$}
\def\proof{\par\noindent {\it Proof:}\enspace}
\def\references#1{\par
  \centerline{\boldsectionfont References}\smallskip
     \parindent=#1pt\nspace}
\def\Ref[#1]{\par\hang\indent\llap{\hbox to\parindent
     {[#1]\hfil\enspace}}\ignorespaces}
\def\Item#1{\par\smallskip\hang\indent\llap{\hbox to\parindent
     {#1\hfill$\,\,$}}\ignorespaces}
\def\ItemItem#1{\par\indent\hangindent2\parindent
     \hbox to \parindent{#1\hfill\enspace}\ignorespaces}

\def\Le{{\mathchoice{\,{\scriptstyle\le}\,}
  {\,{\scriptstyle\le}\,}
  {\,{\scriptscriptstyle\le}\,}{\,{\scriptscriptstyle\le}\,}}}
\def\Ge{{\mathchoice{\,{\scriptstyle\ge}\,}
  {\,{\scriptstyle\ge}\,}
  {\,{\scriptscriptstyle\ge}\,}{\,{\scriptscriptstyle\ge}\,}}}

\def\arrowsim{\,\smash{\mathop{\to}\limits^{\lower1.5pt
  \hbox{$\scriptstyle\sim$}}}\,}

\def\doublemaprights#1#2#3#4{\raise3pt\hbox{$\mathop{\,\,\hbox to
     #1pt{\rightarrowfill}\kern-30pt\lower3.95pt\hbox to
     #2pt{\rightarrowfill}\,\,}\limits_{#3}^{#4}$}}

\def\rightcapdownarrow{\raise9pt\hbox{$\ssize\cap$}\kern-7.75pt
     \Big\downarrow}

\def\rcapmapdown#1{\rightcapdownarrow\kern-1.0pt\vcenter{
     \hbox{$\scriptstyle#1$}}}

\def\rmapdown#1{\Big\downarrow\kern-1.0pt\vcenter{
     \hbox{$\scriptstyle#1$}}}
\def\rightsubsetarrow#1{{\ssize\subset}\kern-4.5pt\lower2.85pt
     \hbox to #1pt{\rightarrowfill}}
\def\longtwoheadedrightarrow#1{\raise2.2pt\hbox to #1pt{\hrulefill}
     \!\!\!\twoheadrightarrow}

\def\Gal{\operatorname{\hbox{Gal}}}
\def\Hom{\operatorname{\hbox{Hom}}}

\def\im{\hbox{Im}}

\NoBlackBoxes
\parindent=25pt
\document
\footline={\hfil}

\null 
\vskip 0.4 cm 
\centerline{\BIGtitle Some Cases of the Mumford--Tate Conjecture} 
\centerline{\BIGtitle and Shimura Varieties} 
\vskip 0.4 in 
\centerline{\bigtitle Adrian Vasiu, Binghamton University} 
\vskip 0.6 cm
\centerline{December 18, 2007}
\centerline{Final version, to appear in Indiana Univ. Math. J.}
\vskip 0.6 cm 

\centerline{Dedicated to Jean-Pierre Serre, on his
$81^{\text{th}}$ anniversary} 
\footline={\hfill} 
\vskip 0.6 cm

\noindent {\bf ABSTRACT}. We prove the Mumford--Tate conjecture for
those abelian varieties over number fields whose extensions to
$\dbC$ have attached adjoint Shimura varieties that are products of
simple, adjoint Shimura varieties of certain Shimura types. In
particular, we prove the conjecture for the orthogonal case (i.e.,
for the $B_n$ and $D_n^{\dbR}$ Shimura types). As a main tool, we
construct embeddings of Shimura varieties (whose adjoints are) of
prescribed abelian type into unitary Shimura varieties of PEL type.
These constructions implicitly classify the adjoints of Shimura
varieties of PEL type.

\vskip 0.2 in
\centerline{\bigtitle Contents}

{\nspace{

\bigskip
\line{\item{1.}{Introduction}\leaderfill 1}

\smallskip
\line{\item{2.}{Some complements on Shimura pairs} \leaderfill 11}

\smallskip
\line{\item{3.}{Basic techniques and the proof of Theorem 1.3.1} \leaderfill 17}

\smallskip
\line{\item{4.}{Injective maps into unitary Shimura pairs of PEL type} \leaderfill 22}

\smallskip
\line{\item{5.}{Applications to Frobenius tori} \leaderfill 34}

\smallskip
\line{\item{6.}{Non-special $A_n$ types} \leaderfill 42}

\smallskip
\line{\item{7.}{The proof of the Main Theorem} \leaderfill 48}

\smallskip
\line{\item{}{References}\leaderfill 63}

}}

\footline={\hss\tenrm \folio\hss}
\pageno=1

\bigskip\smallskip
\noindent
{\boldsectionfont 1. Introduction}
\bigskip

If $K$ is a field, let $\overline{K}$ be an algebraic closure of $K$. If $M$ is a free module of finite rank over a commutative ring $R$ with unit, let $M^*:=\Hom_R(M,R)$ and let $\text{\bf GL}_M$ be the group scheme over $R$ of linear automorphisms of $M$. If $*_R$ or $*$ is either an object or a morphism of the category of $\Spec(R)$-schemes, let $*_U$ be its pull back via an affine morphism $m:\Spec(U)\to\Spec(R)$. If there exist several ways to take $m$, we mention the homomorphism $R\to U$ that defines $m$ but often we do not add it to the notation $*_U$.

\bigskip\noindent
{\bf 1.1. The conjecture.}
Let $E$ be a number field. We fix an embedding $i_E:E\hookrightarrow\dbC$. To fix the ideas, we will take $\overline{E}=\overline{\dbQ}$ to be the algebraic closures of $E=i_E(E)$ and $\dbQ$ in $\dbC$. Let $A$ be an abelian variety over $E$. For the simplicity of notations, let
$$L_A:=H_1(A({\dbC}),\dbZ)$$
be the first homology group of the complex manifold $A(\dbC)$ with coefficients in $\dbZ$, let
$$\mu_A\colon\dbG_{m,\dbC}\to \text{\bf GL}_{L_A\otimes_{\dbZ} \dbC}$$
be the {\it Hodge cocharacter}, and let $H_A$ be the {\it Mumford--Tate group} of $A_{\dbC}$. We recall that if $L_A\otimes_{\dbZ} \dbC=F^{-1,0}\oplus F^{0,-1}$ is the classical Hodge decomposition, then $\beta\in\dbG_{m,\dbC}(\dbC)$ acts through $\mu_A$ trivially on $F^{0,-1}$ and as the multiplication with $\beta$ on $F^{-1,0}$. We also recall that $H_A$ is a {\it reductive group} over $\dbQ$ and that $H_A$ is the smallest subgroup of $\text{\bf GL}_{L_A\otimes_{\dbZ} \dbQ}$ with the property that $\mu_A$ factors through $H_{A,\dbC}$, cf. [De3, Props. 3.4 and 3.6].

Let $p\in\dbN$ be a prime. Let $T_p(A)$ be the {\it Tate-module} of $A$. As a $\dbZ_p$-module (resp. as a $\Gal(E)$-module) we identify it canonically with $L_A\otimes_{\dbZ} \dbZ_p$ (resp. $\Hom_{\dbZ_p}(H^1_{\acute{et}}(A_{\overline{E}},\dbZ_p),\dbZ_p)$). Let $G_{\dbQ_p}$ be the {\it identity component of the algebraic envelope} of the $p$-adic Galois representation
$$\rho\colon\Gal(E)\to \text{\bf GL}_{T_p(A)\otimes_{\dbZ_p} \dbQ_p}(\dbQ_p)=\text{\bf GL}_{L_A\otimes_{\dbZ} \dbQ_p}(\dbQ_p),$$
cf. [Bog] and [Se3, Sect. 135]. Let $E^{\text{conn}}$ be the smallest finite field extension of $E$ that is contained in $\overline{\dbQ}$ and that has the property that $\rho(\Gal(E^{\text{conn}}))$ is a compact, open subgroup of $G_{\dbQ_p}(\dbQ_p)$ (see [Bog]). It is known that the field $E^{\text{conn}}$ does not depend on the prime $p$ (see [Se3, p. 15]). The Mumford--Tate conjecture says (see [Mu1] and [Se2, Sect. 9]):

\medskip\noindent
{\bf 1.1.1. Conjecture.}
{\it As subgroups of $\text{\bf GL}_{L_A\otimes_{\dbZ} \dbQ_p}$ we have $G_{\dbQ_p}=H_{A,\dbQ_p}$.}

\medskip
The main goal of the paper is to prove the Conjecture in many cases. We first recall fundamental results that pertain to the Conjecture and that play a key role in the paper.

\medskip\noindent
{\bf 1.1.2. Theorem (Pyatetskii-Shapiro--Deligne--Borovoi).}
{\it As subgroups of $\text{\bf GL}_{L_A\otimes_{\dbZ} \dbQ_p}$, $G_{\dbQ_p}$ is a subgroup of $H_{A,\dbQ_p}$ (for instance, cf. [De3, Prop. 2.9 and Thm. 2.11]).}

\medskip\noindent
{\bf 1.1.3. Theorem (Faltings).}
{\it The group $G_{\dbQ_p}$ is reductive (cf. [Fa]).}

\medskip\noindent
{\bf 1.1.4. Theorem (Faltings).} {\it The centralizers of $G_{\dbQ_p}$ and $H_{A,\dbQ_p}$ in $\text{\bf GL}_{L_A\otimes_{\dbZ} \dbQ_p}$ coincide; more precisely the centralizers of $G_{\dbQ_p}(\dbQ_p)$ and $H_{A,\dbQ_p}(\dbQ_p)$ in $\End(L_A\otimes_{\dbZ} \dbQ_p)$ are equal to the semisimple $\dbQ_p$-subalgebra $\End(A_{\overline{E}})\otimes_{\dbZ} \dbQ_p$ of $\End(L_A\otimes_{\dbZ} \dbQ_p)$ (cf. [Fa]).}

\medskip
Let $v$ be a prime of $E$ such that $A$ has good reduction $A_v$ with respect to it. Let $l\in\dbN$ be the prime divided by $v$. Let $T_v$ be the {\it Frobenius torus} of the abelian variety $A_v$ (see [Ch]); it is a torus over $\dbQ$ that depends only on the isogeny class of $A_v$. If $p\neq l$, then $T_{v,\dbQ_p}$ is naturally a torus of $G_{\dbQ_p}$, uniquely determined up to $\im(\rho)$-conjugation.

\medskip\noindent
{\bf 1.1.5. Theorem (Serre).} {\it Let $\grP$ be the set of primes $v$ of $E$ which are of good reduction for $A$, which are relatively prime to $p$, and for which $T_{v,\dbQ_p}$ is a maximal torus of $G_{\dbQ_p}$. If $E=E^{\text{conn}}$, then the set $\grP$ has Dirichlet density 1 (cf. [Ch, Cor. 3.8]). Thus the rank of $G_{\dbQ_p}$ does not depend on $p$ (cf. [Ch, Thm. 3.10]).}

\medskip\noindent
{\bf 1.1.6. Theorem (Larsen--Pink--Tankeev).} {\it If the Mumford--Tate conjecture for the abelian variety $A$ is true for a prime $p$, then it is true for all rational primes (cf. either [LP, Thm. 4.3] or [Ta1, Lem. 3.4]).}

\medskip
The proof of [Ta1, Lem. 3.4] relies on Theorems 1.1.2 to 1.1.5 and on the following Lemma.

\medskip\noindent
{\bf 1.1.7. Lemma (Zarhin).} {\it Let $\scrV$ be a finite dimensional vector space over a field of characteristic 0. Let $\grg_1\subseteq\grg_2\subseteq\End(\scrV)$ be Lie algebras of reductive subgroups of $\text{\bf GL}_{\scrV}$. We assume that the centralizers of $\grg_1$ and $\grg_2$ in $\End(\scrV)$ coincide and that the ranks of $\grg_1$ and $\grg_2$ are equal. Then $\grg_1=\grg_2$ (cf. [Za1, Sect. 5, Key Lemma]).}

\medskip\noindent
{\bf 1.1.8. Theorem (Pink).} {\it The group $G_{\overline{\dbQ_p}}$ is generated by cocharacters of it with the property that the resulting representations of $\dbG_{m,\overline{\dbQ_p}}$ on $L_A\otimes_{\dbZ} \overline{\dbQ_p}$ are associated to the trivial and the identity characters of $\dbG_{m,\overline{\dbQ_p}}$, the multiplicities being equal (cf. [Pi, Thms. 3.18 and 5.10]).}

\medskip\noindent
{\bf 1.1.9. Corollary (Pink).} {\it Every simple factor $\scrF_1$ of $G^{\ad}_{\overline{\dbQ_p}}$ is of classical Lie type and the highest weight of each non-trivial, simple  $\Lie(\scrF_1)$-submodule of $L_A\otimes_{\dbZ} \overline{\dbQ_p}$, is a minuscule weight (cf. [Pi, Cor. 5.11]; see also [Ta2] for a slightly less general result).}

\bigskip\noindent
{\bf 1.2. Conventions, notations, and preliminary organization.} A reductive group $H$ over a field $K$ is assumed to be connected. Let $Z(H)$, $H^{\der}$, $H^{\ad}$, and  $H^{\ab}$ denote the center, the derived group, the adjoint group, and the abelianization (respectively) of $H$. We have $Z^{\ab}(H)=H/H^{\der}$ and $H^{\ad}=H/Z(H)$. Let $Z^0(H)$ be the maximal subtorus of $Z(H)$. If $K$ is a finite, separable extension of another field $\tilde K$, let $\Res_{K/\tilde K} H$ be the reductive group over $\tilde K$ obtained from $H$ through the Weil restriction of scalars (see [BLR, Ch. 7, Subsect. 7.6] and [Va2, Subsect. 2.3]). For each commutative $\tilde K$-algebra $\tilde R$, we have a functorial group identity $\Res_{K/\tilde K} H(\tilde R)=H(\tilde R\otimes_{\tilde K} K)$. Let $\dbS:=\Res_{\dbC/\dbR} \dbG_{m,\dbR}$.

Always $n$ is a natural number. If $M$ and $R$ and as before Subsection 1.1 and if $\psi$ is a perfect alternating form on the free $R$-module $M$, then $\text{\bf Sp}(M,\psi)$ and $\text{\bf GSp}(M,\psi)$ are viewed as reductive group schemes over $R$.

A {\it Shimura pair} $(H,X)$ consists of a reductive group $H$ over $\dbQ$ and an $H(\dbR)$-conjugacy class $X$ of homomorphisms $\dbS\to H_{\dbR}$ that satisfy Deligne's axioms of [De2, Subsect. 2.1]: the Hodge $\dbQ$--structure on $\Lie(H)$ defined by any $x\in X$ is of type $\{(-1,1),(0,0),(1,-1)\}$, $\Ad\circ x(i)$ is a Cartan involution of $\Lie(H^{\ad}_{\dbR})$, and no simple factor of $H^{\ad}$ becomes compact over $\dbR$. Here $\Ad:H_{\dbR}\to \text{\bf GL}_{\Lie(H^{\ad}_{\dbR})}$ is the adjoint representation. These axioms imply that $X$ has a natural structure of a hermitian symmetric domain, cf. [De2, Cor. 1.1.17]. For generalities on (the types of) the Shimura variety $\Sh(H,X)$ defined by $(H,X)$ we refer to [De1], [De2], [Mi2], [Mi3], and [Va1, Subsects. 2.2 to 2.5]. Shimura pairs of the form $(\text{\bf GSp}(W,\psi),S)$ define {\it Siegel modular varieties}. The adjoint (resp. the toric part) of a Shimura pair $(H,X)$ is denoted by $(H^{\ad},X^{\ad})$ (resp. $(H^{\ab},X^{\ab})$), cf. [Va1, Subsubsect. 2.4.1].

We say $(H,X)$ (or $\Sh(H,X)$) is {\it unitary} if the group $H^{\ad}$ is non-trivial and each simple factor of $H^{\ad}_{\overline{\dbQ}}$ is of some $A_n$ Lie type.

\medskip\noindent
{\bf 1.2.1. On techniques.}
In [Va1] we developed some techniques to study Shimura varieties and local deformations of abelian varieties endowed with de Rham tensors. With this paper we begin to exploit them to achieve progress in the proof of Conjecture 1.1.1. We mostly use techniques from reductive groups (see Sections 3, 6, and 7) and from Shimura varieties (see Sections 2, 4, 5, 6, and 7). See Subsection 1.3 for what the techniques of Sections 3 to 7 can achieve in connection to Conjecture 1.1.1.

For more details on how this work is structured see Subsections 1.5 and 1.6.

\bigskip\noindent
{\bf 1.3. Our results on Conjecture 1.1.1.}
The group of homotheties of $L_A\otimes_{\dbZ}\dbQ_p$ is a torus of $G_{\dbQ_p}$ (see [Bog]). The following theorem refines this result.

\medskip\noindent
{\bf 1.3.1. Theorem.} {\it The two tori $Z^0(H_{A,\dbQ_p})$ and $Z^0(G_{\dbQ_p})$ of $\text{\bf GL}_{L_A\otimes_{\dbZ}\dbQ_p}$ coincide i.e., we have $Z^0(H_{A,\dbQ_p})=Z^0(G_{\dbQ_p})$.}

\medskip
Let $(H_A,X_A)$ be the {\it Shimura pair attached to} $A_{\dbC}$. Thus $X_A$ is the $H_A(\dbR)$-conjugacy class of the monomorphism $h_A:\dbS\hookrightarrow H_{A,\dbR}$ that defines the classical Hodge $\dbZ$-structure on $L_A$. We refer to $(H_A^{\ad},X_A^{\ad})$ as the {\it adjoint Shimura pair attached to} $A_{\dbC}$ and we write it
$$(H_A^{\ad},X_A^{\ad})=\prod_{j\in J} (H_j,X_j)$$
as a finite product of simple, adjoint Shimura pairs indexed by a set $J$ which is assumed to be disjoint from $\dbZ$ (in order not to create confusion with homology groups).

Each $H_j$ is a simple group over $\dbQ$. The simple factors of $H_{j,\overline{\dbQ}}$ have the same Lie type $\grL_j$ which is (cf. [Sa1], [Sa2], or [De2]) a classical Lie type. If $\grL_j$ is the $A_n$, $B_n$, or $C_n$ Lie type, then $(H_j,X_j)$ is said to be of $A_n$, $B_n$, or $C_n$ type (respectively). If $n\Ge 4$ and $\grL_j=D_n$, then $(H_j,X_j)$ is of either $D_n^{\dbH}$ or $D_n^{\dbR}$ type (cf. [De2, Table 2.3.8]; see Definition 2.2.1 (a) for the $D_n^{\dbH}$ and $D_n^{\dbR}$ types). See Subsection 2.1 for the classical groups $\text{\bf SO}^*(2n)_{\dbR}$ and $\text{\bf SO}(2,2n-2)_{\dbR}$ over $\dbR$. If $(H_j,X_j)$ is of $D_n^{\dbH}$ (resp. of $D_n^{\dbR}$) type, then all simple, non-compact factors of $H_{j,\dbR}$ are isomorphic to $\text{\bf SO}^*(2n)_{\dbR}^{\ad}$ (resp. to $\text{\bf SO}(2,2n-2)^{\ad}_{\dbR}$). The converse of the last sentence holds provided we have $n\Ge 5$.

We say $(H_j,X_j)$ is of {\it non-inner type} if $H_j$ is isomorphic to $\Res_{F_j/\dbQ} \tilde H_j$, where $F_j$ is a number field and where $\tilde H_j$ is an absolutely simple, adjoint group over $F_j$ which is not an inner form of a split group (see Definition 2.3 (a)).

\medskip\noindent
{\bf 1.3.2. Theorem.}
{\it We assume that $p\neq l$ and that $v$ is unramified over $l$. There exists a number $N(A)\in\dbN$ that depends only on $A$ (and not on $p$ or $l$) and such that the following statements holds: if $l\Ge N(A)$ and $j\in J$, then all images of $T_{v,\overline{\dbQ_p}}$ in simple factors of $H_{j,\overline{\dbQ_p}}$ have equal ranks that do not depend on the prime $p\neq l$.}

\medskip
If the abelian variety $A_v$ is ordinary, then Theorem 1.3.2 follows easily from [No, Thm. 2.2]. The following Corollary complements Theorem 1.1.5.

\medskip\noindent
{\bf 1.3.3. Corollary.} {\it The rank of the subgroup of $H_{A,\dbQ_p}^{\ad}$ whose extension to $\overline{\dbQ_p}$ is the product of the images of $G_{\overline{\dbQ_p}}$ in the simple factors of $H_{A,\overline{\dbQ_p}}^{\ad}$, does not depend on $p$.}

\medskip\noindent
{\bf 1.3.4. Main Theorem.}
{\it Conjecture 1.1.1 holds if for each element $j\in J$ the adjoint Shimura pair $(H_j,X_j)$ is of one of the following types:

\medskip
{\bf (a)} non-special $A_n$ type;

\smallskip
{\bf (b)} $B_n$ type;

\smallskip
{\bf (c)} $C_n$ type, with $n$ odd;

\smallskip
{\bf (d)} $D_n^{\dbH}$ type, with $n\Ge 5$ odd and such that $2n\notin\{\binom {2^{m+1}}{2^m}|m\in\dbN\}$;

\smallskip
{\bf (d')} non-inner $D_{2n}^{\dbH}$ type, with $n$ an odd prime;

\smallskip
{\bf (e)} $D_n^{\dbR}$ type, with $n\Ge 4$.}

\medskip
The types listed in Theorem 1.3.4 and the non-inner $D_4^{\dbH}$ type are the only types for which we can show that for each $j\in J$ the group $G_{\overline{\dbQ_p}}$ surjects onto a simple factor of $H_{j,\overline{\dbQ_p}}$ (see Theorem 7.1). The non-special $A_n$ types are introduced in Definition 6.1 (c). Here we list some examples, cf. Examples 6.2.4.

\medskip\noindent
{\bf 1.3.5. Examples.}
The simple, adjoint Shimura pair $(H_j,X_j)$ is of non-special $A_n$ type if any one of the following four conditions is satisfied:

\medskip
{\bf (i)} the group $H_{j,\dbR}$ has two simple factors $\text{\bf SU}(a_i,n+1-a_i)^{\ad}_{\dbR}$, $i=\overline{1,2}$, such that $1\Le a_1< a_2\le {{n+1}\over 2}$ and either $g.c.d.(a_1,n+1)=1$ or $g.c.d.(a_2,n+1)=1$;

\smallskip
{\bf (ii)} the group $H_{j,\dbR}$ has a simple factor $\text{\bf SU}(a,n+1-a)^{\ad}_{\dbR}$, with $g.c.d.(a,n+1)=1$ and with the pair $(a,n+1-a)\notin\{(\binom {r}{s},\binom {r}{s-1})|s,r\in\dbN,\;2\le s\le r-1\}$;

\smallskip
{\bf (iii)} the natural number $n+1$ is an integral power of $2$ and the group $H_{j,\dbR}$ has as a normal factor a product $\text{\bf SU}({{n+1}\over 2},{{n+1}\over 2})^{\ad}_{\dbR}\times_{\dbR} \text{\bf SU}(b,n+1-b)_{\dbR}^{\ad}$, where $b\in\dbN$ is an odd number different from ${{n+1}\over 2}$;

\smallskip
{\bf (iv)} the Shimura pair $(H_j,X_j)$ is of $A_n$ type, where $n+1$ is either $4$ or a prime.

\medskip\noindent
{\bf 1.3.6. Theorem.} {\it We assume that there exists no element $j\in J$ such that $(H_j,X_j)$ is of $D_4^{\dbH}$ type. Then no simple factor of $\Lie(G_{\overline{\dbQ_p}}^{\ad})$ surjects onto two or more distinct simple factors of $\Lie(H_{A,\overline{\dbQ_p}}^{\ad})$.}

\medskip\noindent
{\bf 1.3.7. Theorem (the independence property for abelian varieties).} {\it We assume that there exists no element $j\in J$ such that $(H_j,X_j)$ is of $D_4^{\dbH}$ type and that $A$ is a product of abelian varieties over $E$ for which the Mumford--Tate conjecture holds. Then Conjecture 1.1.1 holds (i.e., the Mumford--Tate conjecture holds also for $A$).}

\medskip
Theorems 1.3.1 and 1.3.2 are used in the proof of Theorem 1.3.4. Theorem 1.3.6 is a consequence of part of the proof of Theorem 1.3.4. Theorem 1.3.7 is an application of Theorem 1.3.6.

\bigskip\noindent
{\bf 1.4. Earlier results that pertain to Conjecture 1.1.1.}
Previous works of Serre, Chi, Tankeev, Pink, Larsen, Zarhin, and few others proved Conjecture 1.1.1 in many cases. The cases involve Shimura pairs $(H_A,X_A)$ whose adjoints are simple of $A_n$, $B_n$, $C_n$, and $D_n^{\dbH}$ types (with $n$ arbitrary large) and are worked out under the assumption that either the rank of $\End(A_{\overline{\dbQ}})$ is small (like $1$ or $4$) or $\dim(A)\in\{qr|q\;\text{is}\;\text{a}\;\text{prime},\;r\in\{1,\ldots,9\}\}$ (see [Ch], [Ta1], [Ta2], [LP], [Pi, p. 216 and Thms. 5.14 and 5.15], and [Za2]). To our knowledge, except [Pi, Thm. 5.14] and [Za2], all previous results that pertain to Conjecture 1.1.1 are encompassed by Theorem 1.3.4. The result [Pi, Thm. 5.14] proves Conjecture 1.1.1 when $\End(A_{\overline{\dbQ}})=\dbZ$ and $2\dim(A)\notin\{m^k,\binom {2k}{k}|m\in\dbN,\;k\in 1+2\dbN\}$. Results of [Za2] prove Conjecture 1.1.1 when $\End(A_{\overline{\dbQ}})=\dbZ$ and $A$ is the Jacobian of certain hyperelliptic curves. In both places [Pi, Thm. 5.14] and [Za2], the group $H_A$ is a $\text{\bf GSp}_{2\dim(A),\dbQ}$ group and therefore the adjoint Shimura pair $(H_A^{\ad},X_A^{\ad})$ is simple of $C_{\dim(A)}$ type.

We recall from [Pi, Def. 4.1] that a pair $(\scrG_0,\rho_0)$ defined by an irreducible, faithful, finite dimensional representation $\rho_0:\scrG_0\hookrightarrow \text{\bf GL}_{\scrW_0}$ over $\overline{\dbQ_p}$, is called an {\it irreducible weak} (resp. an {\it irreducible strong}) {\it Mumford--Tate pair of weight $\{0,1\}$} if $\scrG_0$ is a reductive group over $\overline{\dbQ_p}$ that is generated by cocharacters of it which act  (resp. is generated by the $\scrG_0(\overline{\dbQ_p})$-conjugates of one cocharacter of it which acts) on $\scrW_0$ via the trivial and the identity characters of $\dbG_{m,\overline{\dbQ_p}}$. We say that the pair $(\scrG_0,\rho_0)$ is non-trivial, if $\dim_{\overline{\dbQ_p}}(\scrW_0)\ge 2$.

Irreducible weak Mumford--Tate pairs of weight $\{0,1\}$ are naturally attached to the $p$-adic Galois representation $\rho$ of Subsection 1.1 as follows. We take $\scrW_0$ to be a simple $H_{A,\overline{\dbQ_p}}$-submodule of $L_A\otimes_{\dbZ} \overline{\dbQ_p}$. Let $\scrG_0$ and $\scrH_0$ be the images of $G_{\overline{\dbQ_p}}$ and $H_{A,\overline{\dbQ_p}}$ (respectively) in $\text{\bf GL}_{\scrW_0}$. The natural faithful representation $\rho_0:\scrG_0\hookrightarrow \text{\bf GL}_{\scrW_0}$ is irreducible, cf. Theorem 1.1.4. Theorem 1.1.8 implies that the pair $(\scrG_0,\rho_0)$ is an irreducible weak Mumford--Tate pair of weight $\{0,1\}$. From the very definition of $H_A$, one gets that the pair $(\scrH_0,\scrH_0\hookrightarrow\text{\bf GL}_{\scrW_0})$ is also an irreducible weak Mumford--Tate pair of weight $\{0,1\}$.

The previous works mentioned above rely on the classification of (non-trivial) irreducible weak and strong Mumford--Tate pairs of weight $\{0,1\}$ over $\overline{\dbQ_p}$ (see [Pi]; see also [Sa2] and [Se1]) and on the results 1.1.2 to 1.1.9.

\medskip\noindent
{\bf 1.4.1. What the results 1.3.4 to 1.3.7 bring new.} To our knowledge, Theorem 1.3.4 is the very first result on Conjecture 1.1.1 which is stated in terms of Shimura types (and thus which is stated independently on the structures of either $\End(A_{\overline{\dbQ}})$ or $H^{\ab}$). If $H_A^{\ad}$ is an absolutely simple, adjoint group, then only the $D_n^{\dbR}$ type case of Theorem 1.3.4 is completely new. The $A_n$ and $D_n^{\dbR}$ types cases of Theorem 1.3.4 are new general instances for the validity of Conjecture 1.1.1. For instance, the $D_n^{\dbR}$ type case of Theorem 1.3.4 was not known before even when $\End(A_{\overline{\dbQ}})=\dbZ$. Also, the ``applicability" of Theorem 1.3.4 (a) goes significantly beyond the numerical tests of Subsection 6.2 (see Remark 7.5.1 (a)). See Subsection 7.5 for the simplest cases not settled by Theorem 1.3.4 (they contain the case when $\dim(A)=4$ and $H_A$ is a $\text{\bf GSp}_{8,\dbQ}$ group). See Corollary 7.4 for a variant of Theorem 1.3.4 (d') that handles to some extent the non-inner $D_4^{\dbH}$ types. Theorem 1.3.7 is a particular case of a more general result on abelian motives over $E$ (see Theorem 7.6); we do not know any non-trivial analogue of Theorems 1.3.7 and 7.6 in the literature.

The proofs of the results 1.3.1 to 1.3.7 rely on:

\medskip
$\bullet$ results of {\it arithmetic nature} (like the results 1.1.2 to 1.1.5, 1.1.8, and 1.1.9 and the main results of [Wi] and [Zi]);

\smallskip
$\bullet$ results of {\it group theoretical nature} (like the results 1.1.6 and 1.1.7, the classification of irreducible weak and strong Mumford--Tate pairs of weight $\{0,1\}$, and the Satake--Deligne results on symplectic embeddings to be recalled in Subsection 1.5 below).

\medskip
Our new ideas that put together all these results are presented below in Subsections 1.5 and 1.6. The main contribution of these new ideas is that they allow us: (i) to state Theorem 1.3.4 only in terms of Shimura types and (ii) to prove Theorems 1.3.6 and 1.3.7.

\bigskip\noindent
{\bf 1.5. On injective maps into unitary Shimura pairs of PEL type.}
Let $\lambda_A$ be a polarization of $A$. Let $\psi_A$ be the non-degenerate alternating form on $W_A:=L_A\otimes_{\dbZ} \dbQ$ induced by $\lambda_A$. We get naturally an injective map
$$f_A:(H_A,X_A)\hookrightarrow (\text{\bf GSp}(W_A,\psi_A),S_A)$$
of Shimura pairs, where $S_A$ is the $\text{\bf GSp}(W_A,\psi_A)(\dbR)$-conjugacy class of $f_{A,\dbR}\circ h_A$. As a hermitian symmetric domain, $S_A$ is two disjoint copies of the Siegel domain of genus $\dim(A)$. The inclusion $X_A\hookrightarrow S_A$ is a totally geodesic and holomorphic embedding between hermitian symmetric domains. Such embeddings were classified in [Sa1].

Satake and Deligne used rational variants of [Sa1] to show that each simple, adjoint Shimura pair of $A_n$, $B_n$, $C_n$, $D_n^{\dbH}$, or $D_n^{\dbR}$ type is isomorphic to the adjoint Shimura pair attached to some complex abelian variety (cf. [Sa2, Part III] and [De2, Prop. 2.3.10]). The goal of Section 4 is to present unitary variants of [De2, Subsubsects. 2.3.10 to 2.3.13]. Starting from a simple, adjoint Shimura pair $(G,X)$ of $A_n$, $B_n$, $C_n$, $D_n^{\dbH}$, or $D_n^{\dbR}$ type, we construct explicitly in Subsections 4.1 to 4.8 ``good'' injective maps
$$f_2:(G_2,X_2)\hookrightarrow (G_4,X_4)\hookrightarrow (\text{\bf GSp}(W_1,\psi_1),S_1)$$
of Shimura pairs, where $(G_2^{\ad},X_2^{\ad})=(G,X)$, where $(G_4^{\ad},X_4^{\ad})$ is a simple, adjoint Shimura pair of $A_{n_4}$ type for some positive integer $n_4$ that depends only on the type of $(G,X)$, and where $G_4$ is the subgroup of $\text{\bf GSp}(W_1,\psi_1)$ that fixes a semisimple $\dbQ$--subalgebra of $\End(W_1)$. If $(G,X)$ is of $A_n$ type, then we can take $G_4=G_2$ (cf. Proposition 4.1) and thus we have $n_4=n$. If  $(G,X)$ is not of $A_n$ type, see Theorem 4.8 for the extra ``good'' properties enjoyed by the injective maps $f_2$ of Shimura pairs.

These constructions play key roles in the proofs of the Theorems 1.3.2, 1.3.4, 1.3.6, and 1.3.7 (see Subsection 1.6 below). We also use these constructions to classify in Corollary 4.10 all Shimura pairs that are adjoints of Shimura pairs of PEL type (their definition is recalled in Subsection 2.4; for instance, the Shimura pair $(G_4,X_4)$ is of PEL type).

The constructions rely on a modification of the proof of [De2, Prop. 2.3.10] (see proof of Proposition 4.1) and on a natural {\it twisting process} of certain homomorphisms between semisimple groups over number fields (see Subsection 4.5). The embeddings of hermitian symmetric domains of classical Lie type into Siegel domains constructed in [Sa1], are obtained by first achieving embeddings into hermitian symmetric domains associated to unitary groups over $\dbR$. This idea is not used in [De2] and thus the twisting process of Subsection 4.5 is the main new ingredient needed to get constructions that are explicit rational versions of [Sa1] and that are more practical than the ones of [Sa2].

\bigskip\noindent
{\bf 1.6. On the proofs of the results 1.3.1 to 1.3.7.}
The proof of Theorem 1.3.1 relies on an elementary Lemma 3.3 and appeals to canonical split cocharacters of filtered $F$-isocrystals over finite fields defined in [Wi, p. 512] (see Subsection 3.4). Lemma 3.3 is due in essence to Bogomolov and  Serre, cf. [Bog] and [Se3, pp. 31, 43, and 44]. Though Theorem 1.3.1 was known to  Serre (and to others), we could not trace it in the literature.

We call an abelian variety $B$ over $E$ {\it adjoint-isogenous} to $A$ if the adjoints of the Mumford--Tate groups of $A_{\dbC}$, $B_{\dbC}$, and $A_{\dbC}\times_{\dbC} B_{\dbC}$ are the same (cf. Definition 5.1). Due to Theorems 1.3.1, 1.1.2, 1.1.3, and 1.1.4, to prove Conjecture 1.1.1 one only has to show that $G^{\ad}_{\dbQ_p}=H_{A,\dbQ_p}^{\ad}$. But to check that we have $G^{\ad}_{\dbQ_p}=H_{A,\dbQ_p}^{\ad}$, we can replace $A$ by any other abelian variety $B$ over $E$ that is adjoint-isogenous to $A$ (cf. the shifting process expressed by Proposition 5.4.1). In other words, (as the prime $p$ is arbitrary) we have:

\medskip\noindent
{\bf 1.6.1. Proposition.} {\it We assume that $A$ and $B$ are adjoint-isogenous abelian varieties over a number field. Then the Mumford--Tate conjecture holds for $A$ if and only if it holds for $B$.}

\medskip\noindent
{\bf 1.6.2. Example.} We assume that the adjoint group $H_A^{\ad}$ is absolutely simple and that the simple, adjoint Shimura pair $(H_A^{\ad},X_A^{\ad})$ is of $D_n^{\dbR}$ type with $n\ge 5$.
Let $H_A^{\sc}$ be the simply connected semisimple group cover of $H_A^{\ad}$ (or of $H_A^{\der}$). Let $K$ be a number field such that the group $H_{A,K}^{\sc}$ is split (i.e., it is a $\text{\bf Spin}_{2n,K}$ group). Let $H_{A,K}^{\sc}\hookrightarrow \text{\bf GL}_{2^n,K}$ be the spin representation (i.e., the direct sum of the two half spin representations of $H_{A,K}^{\sc}$). We have a natural sequence of monomorphisms
$$H_A^{\sc}\hookrightarrow \Res_{K/\dbQ} H_{A,K}^{\sc}\hookrightarrow \Res_{K/\dbQ}\text{\bf GL}_{2^n,K}\hookrightarrow  \text{\bf GL}_{2^n[K:\dbQ],\dbQ}\hookrightarrow \text{\bf GSp}_{2^{n+1}[K:\dbQ],\dbQ}$$
between reductive groups over $\dbQ$. We write $\text{\bf GSp}_{2^{n+1}[K:\dbQ],\dbQ}=\text{\bf GSp}(W_B,\tilde\psi_B)$, where $W_B$ is a rational vector space of dimension $2^{n+1}[K:\dbQ]$ and where $\tilde\psi_B$ is a non-degenerate alternating form on $W_B$. Let $H_B$ be the reductive subgroup of $\text{\bf GL}_{W_B}$ generated by the subgroups $H_A^{\sc}$ and $Z(\text{\bf GSp}(W_B,\psi_B))$; we have natural identifications $H_B^{\der}=H_A^{\sc}$ and $H_B^{\ad}=H_A^{\ad}$. The centralizer of $H_{B,\overline{\dbQ_p}}$ in $\text{\bf GL}_{W_B\otimes_{\dbQ} \overline{\dbQ_p}}$ is isomorphic to $\text{\bf GL}_{[K:\dbQ],\overline{\dbQ_p}}^2$ (each copy of $\text{\bf GL}_{[K:\dbQ],\overline{\dbQ_p}}$ corresponds to one of the two half spin representations).

Let $h_A^{\ad}=h_B^{\ad}:\dbS\to H_{A,\dbR}^{\ad}=H_{B,\dbR}^{\ad}$ be the composite of the monomorphism $h_A:\dbS\hookrightarrow H_{A,\dbR}$ with the natural epimorphism $H_{A,\dbR}\twoheadrightarrow H_{A,\dbR}^{\ad}$. It is easy to see that there exists a unique monomorphism $h_B:\dbS\hookrightarrow H_{B,\dbR}$ which lifts $h_B^{\ad}$ and which defines a Hodge $\dbQ$--structure on $W_B$ of type $\{(-1,0),(0,-1)\}$, cf. either the proof of [De2, Prop. 2.3.10] or [Pi, Table 4.2]. Let $\psi_B$ be a non-degenerate alternating form on $W_B$ that is fixed by $H_B^{\der}=H_A^{\sc}$ and such that $2\pi i\psi_B$ is a polarization of this Hodge $\dbQ$--structure on $W_B$, cf. [De2, Cor. 2.3.3]. Let $X_B$ be the $H_B(\dbR)$-conjugacy class of $h_B$ and let $S_B$ be the $\text{\bf GSp}(W_B,\psi_B)(\dbR)$-conjugacy class of the composite of $h_B$ with the monomorphism $H_{B,\dbR}\hookrightarrow \text{\bf GSp}(W_B\otimes_{\dbQ} \dbR,\psi_B)$. By replacing $E$ with a finite field extension of it, there exists a principally polarized abelian variety $(B,\lambda_B)$ over $E$ such that the notations match (i.e., the analogues of $h_A$ and $f_A$ for $(B,\lambda_B)$ are precisely $h_B$ and $f_B:(H_B,X_B)\hookrightarrow (\text{\bf GSp}(W_B,\psi_B),S_B)$) and such that the abelian varieties $A$ and $B$ are adjoint-isogenous.

Let $\tilde G_{\dbQ_p}$ be the identity component of the algebraic envelope of the $p$-adic Galois representation attached to $B$; it is a subgroup of $H_{B,\dbQ_p}$. We show that the assumption that $\tilde G_{\dbQ_p}\neq H_{B,\dbQ_p}$ leads to a contradiction. Due to Theorem 1.1.8, this assumption implies that $\tilde G_{\overline{\dbQ_p}}^{\der}$ is isogenous to a product of two simple, adjoint groups of $B_{n_1}$ and $B_{n-n_1-1}$ Lie type (respectively), where $n_1\in\{0,\ldots,n-1\}$ (cf. [Pi, Prop. 4.3]; here the $B_0$ Lie type corresponds to the trivial group). Moreover, loc. cit. implies that the representation $\Lie(\tilde G_{\overline{\dbQ_p}}^{\der})\hookrightarrow \End(W_B\otimes_{\dbQ} \overline{\dbQ_p})$ is a direct sum of $2[K:\dbQ]$ irreducible representations that are isomorphic to the tensor product of the spin representations of the two simple factors of $\Lie(\tilde G_{\overline{\dbQ_p}}^{\der})$ of $B_{n_1}$ and $B_{n-n_1-1}$ Lie type (in other words, the restriction to $\tilde G_{\overline{\dbQ_p}}^{\der}$ of either one of the two half spin representations of $H_{B,\overline{\dbQ_p}}$ is the same irreducible representation of $\tilde G_{\overline{\dbQ_p}}^{\der}$). Thus the centralizer of $\tilde G_{\overline{\dbQ_p}}$ in $\text{\bf GL}_{W_B\otimes_{\dbQ} \overline{\dbQ_p}}$ is isomorphic to $\text{\bf GL}_{2[K:\dbQ],\overline{\dbQ_p}}$ and therefore it is distinct from the centralizer of $H_{B,\overline{\dbQ_p}}$ in $\text{\bf GL}_{W_B\otimes_{\dbQ} \overline{\dbQ_p}}$. This contradicts Theorem 1.1.4.

Therefore we have $\tilde G_{\dbQ_p}=H_{B,\dbQ_p}$. Thus the Mumford--Tate conjecture holds for $B$ (cf. Theorem 1.1.6) and thus also for $A$ (cf. Proposition 1.6.1).
This holds even if $n=4$, provided we define the two half spin representations in such away that $h_B$ does exist. If either $n$ is odd or $n$ is even and the simple, adjoint Shimura pair $(H_A^{\ad},X_A^{\ad})$ is of non-inner $D_n^{\dbR}$ type, then the representation $H_{A,\overline{\dbQ_p}}^{\der}\hookrightarrow \text{\bf GL}_{W_A\otimes_{\dbQ} \overline{\dbQ_p}}$ is a direct sum of spin representations and therefore the above replacement of $A$ by $B$ is not truly necessary.

\medskip\noindent
{\bf 1.6.3. The role of the shifting process.}
The main role of the shifting process (i.e., of Proposition 5.4.1) is to replace the ``potentially complicated" injective map $f_A$ of Shimura pairs by a ``simpler" one that is a ``direct sum" (i.e., is a Hodge quasi product in the sense of Subsection 2.4) of composite injective maps $f_2:(G_2,X_2)\hookrightarrow (G_4,X_4)\hookrightarrow(\text{\bf GSp}(W_1,\psi_1),S_1)$ of Shimura pairs that are as in Subsection 1.5. See Proposition 5.10 for three advantages one gets through such replacements; they are due to Deligne. We emphasize that:

\medskip
{\bf (i)} If $B$ is an abelian variety over $E$ which is adjoint-isogenous to $A$, then in general $B$ is unrelated to $A$ in any usual way (like via isogenies, Weil restrictions of scalars, products of abelian subvarieties of $A$, etc.). For instance, $B$ could be the product between $A$ and any other abelian variety over $E$ which has complex multiplication over $\overline{E}$.

\smallskip
{\bf (ii)} To perform the replacement of $f_A$ by a ``direct sum'' of composite injective maps $f_2$ of Shimura pairs, we often have to pass to a finite field extension of $E$.

\smallskip
{\bf (iii)} If the adjoint group $H_A^{\ad}$ is $\dbQ$--simple, then by replacing $f_A$ with an $f_2$, the Mumford--Tate group $H_A$ gets replaced by a reductive, normal subgroup of $G_2$ which contains $G_2^{\der}$ (more precisely, the injective map $f_A:(H_A,X_A)\hookrightarrow (\text{\bf GSp}(W_A,\psi_A),S_A)$ of Shimura pairs gets replaced by a composite injective map
$$(H_A,X_A)\hookrightarrow (G_2,X_2)\hookrightarrow (G_4,X_4)\hookrightarrow (\text{\bf GSp}(W_1,\psi_1),S_1)=(\text{\bf GSp}(W_A,\psi_A),S_A)$$
of Shimura pairs for which we have $H_A^{\der}=G_2^{\der}$ and which is defined by an $f_2$ as in Subsection 1.5).

\medskip\noindent
{\bf 1.6.4. On Theorem 1.3.2.} The proof of Theorem 1.3.2 goes as follows (see Subsection 5.8). Using the shifting process, we can assume that $H_A^{\ad}$ is a simple $\dbQ$--group and we can replace $f_A$ by an injective map $f_2$ of Shimura pairs as in Subsection 1.5. Using this and the main result of [Zi], we get that for $l\gg 0$ and up to a $G_4(\dbQ_p)$-conjugation of the torus $T_{v,\dbQ_p}$, $T_v$ is naturally isomorphic to a torus $T_4^\prime$ of an inner form $G_4^\prime$ of $G_4$ that is locally isomorphic. But from the constructions of $f_2$ we get that each simple factor of $\Lie(G_{4,\overline{\dbQ_p}}^{\ad})$ contains the non-zero image of a unique simple factor of $\Lie(G_{2,\overline{\dbQ_p}}^{\ad})$ (see property 5.7 (a)). Thus we get Theorem 1.3.2 by considering the ranks of the images of $T^\prime_{4,\overline{\dbQ_p}}$ in simple factors of $G^{\ad}_{4,\overline{\dbQ_p}}=G^{\prime,\ad}_{4,\overline{\dbQ_p}}$. The passage from Theorem 1.3.2 to Corollary 1.3.3 is standard (see Subsection 5.8).

\medskip\noindent
{\bf 1.6.5. On Theorems 1.3.4, 1.3.6, and 1.3.7.} The proof of Theorem 1.3.4 involves two main steps (see Theorems 7.1 and 7.2) and it ends in Subsection 7.3. Theorem 7.1 shows that if for each element $j\in J$ the simple, adjoint Shimura pair $(H_j,X_j)$ is of either non-inner $D_4^{\dbH}$ type or a type listed in Theorem 1.3.4, then for each $j\in J$ the group $G_{\overline{\dbQ_p}}$ surjects onto a simple factor of $H_{j,\overline{\dbQ_p}}$. The greatest part of the proof of Theorem 7.1 is a standard application of the classification of irreducible weak and strong Mumford--Tate pairs of weight $\{0,1\}$ and of the advantages offered by the shifting process.

Theorem 7.2 shows that under some hypotheses (like the ones of Theorem 1.3.4), we have $G_{\dbQ_p}^{\ad}=H_{A,\dbQ_p}^{\ad}$. Its proof is technical and lengthy. To explain the main point of the proofs of Theorems 7.2, 1.3.6, and 1.3.7, we will assume in this and the next paragraph that there exists no element $j\in J$ such that $(H_j,X_j)$ is of $D_4^{\dbH}$ type and that (based on the shifting process) the injective map $f_A$ has been replaced by a ``direct sum'' of injective map $f_2$ as in Subsection 1.5. Even more, to ease the notations we will also assume that the adjoint group $H_A^{\ad}$ is $\dbQ$--simple and thus that $f_A$ has been replaced (as in (iii)) by a single injective map $f_2$. The torus $Z^0(H_A)$ is a subtorus of $Z^0(G_2)$ (cf. (iii)) and its dimension is very much controlled by $G_{\dbQ_p}$. If there exists a simple factor of $\Lie(G_{\overline{\dbQ_p}}^{\ad})$ that surjects onto two or more distinct simple factors of $\Lie(H_{A,\overline{\dbQ_p}}^{\ad}$), then directly from the construction of the injective map $f_2$ of Shimura pairs we get that the centralizers of $G^{\der}_{\dbQ_p}$ and $H_{A,\dbQ_p}^{\der}$ in $\text{\bf GL}_{T_p(A)\otimes_{\dbZ_p} \dbQ_p}$ are distinct and thus one hopes to reach a contradiction with Theorem 1.1.4. But in order to reach the desired contradiction, one has to check that the subtorus $Z^0(H_A)$ of $Z^0(G_2)$ is small enough (so that the centralizers of $G_{\dbQ_p}$ and $H_{A,\dbQ_p}$ in $\text{\bf GL}_{T_p(A)\otimes_{\dbZ_p} \dbQ_p}$ are distinct). Thus the main philosophy of the proof of Theorem 7.2 is:

\medskip
{\bf (iv)} the smaller the subgroup $G_{\dbQ_p}$ of $H_{A,\dbQ_p}\leqslant G_{2,\dbQ_p}$ is known to be (in terms of several simple factors of $\Lie(G_{\overline{\dbQ_p}}^{\ad})$ surjecting onto two or more distinct simple factors of $\Lie(H_{A,\overline{\dbQ_p}}^{\ad})$), the smaller subtori of $Z^0(G_2)$ we get which contain $Z^0(H_A)$ (and thus the better upper bounds we get for $\dim(Z^0(H_A))$).

\medskip
Thus if there exists a simple factor of $\Lie(G_{\overline{\dbQ_p}}^{\ad})$ that surjects onto two or more distinct simple factors of $\Lie(H_{A,\overline{\dbQ_p}}^{\ad}$), then the proof of Theorem 7.2 shows that the torus $Z^0(H_A)$ is small enough so that we reach a contradiction with Theorem 1.1.4. In the proof of Theorem 7.2 we will use very intimately the construction of the injective map $f_2$ of Shimura pairs and the proof of Theorem 1.3.2, in order to get the desired control on the subtorus $Z^0(H_A)$ of $Z^0(G_2)$.

The essence of the proofs of Theorems 1.3.6 and 1.3.7 is the same main point mentioned above. The proof of Theorem 1.3.6 is presented in Subsection 7.3. Theorem 1.3.7 is a particular case of Theorem 7.6.

\medskip\noindent
{\bf Acknowledgments.}
We would like to thank UC at Berkeley, U of Utah, U of Arizona, and Binghamton U for providing us with good conditions with which to write this work. We would like to thank Tankeev for making us aware of Theorem 1.1.6, Serre and Deligne for some comments, Borovoi for mentioning the reference [Sa2], and the referee and  Larsen for many valuable comments and suggestions. We are obliged to [Pi]. This research was partially supported by the NSF grant DMS 97-05376.

\bigskip\smallskip
\noindent
{\boldsectionfont 2. Some complements on Shimura pairs}
\bigskip

In Subsection 2.1 we gather extra conventions and notations to be used in Sections 2 to 6. In Subsections 2.2 and 2.3 we include complements on reflex fields and types of Shimura pairs. In Subsection 2.4 we include complements on injective maps into Shimura pairs that define Siegel modular varieties. The notations of Section 2 will be independent from the ones of Subsections 1.1, 1.3, and 1.5.

\bigskip\noindent
{\bf 2.1. Extra conventions and notations.} Let $p\in\dbN$ be a prime. For each field $K$ and every cocharacter $\mu_0$ of a reductive group $H$ over $K$, let $[\mu_0]$ be the $H(K)$-conjugacy class of $\mu_0$. A pair of the form $(H,[\mu_1])$, with $\mu_1$ a cocharacter of $H_{\overline K}$, is called a {\it group pair} over $K$. If $D$ is a semisimple $K$-algebra, let $\text{\bf GL}_1(D)$ be the reductive group over $K$ of invertible elements of $D$. If $K$ is a $p$-adic field, we call $H$ {\it unramified} if it is quasi-split and splits over an unramified extension of $K$. If $K$ is a number field, let $K_w$ be the completion of $K$ with respect to a finite prime $w$ of $K$. Let $k(w)$ be the residue field of $w$. Let $B(k)$ be the field of fractions of the ring $W(k)$ of Witt vectors with coefficients in a perfect field $k$. Let $\dbF:=\overline{\dbF_p}$. Let $\dbA_f:=\dbZhat\otimes_{\dbZ} \dbQ$ be the $\dbQ$--algebra of finite ad\`eles of $\dbQ$. If an abstract group $C$ acts on a set $V$, let $V^C:=\{y\in V|y\,\, \text{fixed}\,\,\text{by}\,\, C\}$.

Let $(G,X)$ be a Shimura pair. Let $E(G,X)$ be the reflex field of $(G,X)$. The Shimura variety $\Sh(G,X)$ is identified with the canonical model over $E(G,X)$ of the complex Shimura variety $\Sh(G,X)_{\dbC}$ (see [De1], [De2], [Mi1], [Mi3], and [Mi4]).

Let $a$, $b\in\dbN\cup\{0\}$ with $a+b>0$. Let $\text{\bf SU}(a,b)$ be the simply connected semisimple group over $\dbQ$ whose $\dbQ$--valued points are the $\dbQ(i)$--valued points of $\text{\bf SL}_{a+b,\dbQ}$ that leave invariant the hermitian form $-z_1\overline z_1-\cdots-z_a\overline z_a+z_{a+1}\overline z_{a+1}+\cdots+z_{a+b}\overline z_{a+b}$ over $\dbQ(i)$. Let $\text{\bf SO}(a,b)$ be the semisimple group over $\dbQ$ of $a+b$ by $a+b$ matrices of determinant 1 that leave invariant the quadratic form $-x_1^2-\cdots-x_a^2+x_{a+1}^2+\cdots+x_{a+b}^2$ on $\dbQ^{a+b}$. Let $\text{\bf SO}(a):=\text{\bf SO}(0,a)$. Let $\text{\bf SO}^*(2a)$ be the semisimple group over $\dbQ$ whose group of $\dbQ$--valued points is the subgroup of $\text{\bf SO}(2n)(\dbQ(i))$ that leaves invariant the skew hermitian form $-z_1\overline{z}_{n+1}+z_{n+1}\overline{z}_1-\cdots-z_n\overline{z}_{2n}+z_{2n}\overline{z}_n$ over $\dbQ(i)$ ($z_i$'s and $x_i$'s are related here over $\dbQ(i)$ via $z_i=x_i$).

The {\it minuscule weights} of a split, simple Lie algebra of classical Lie type $\grL$  over a field of characteristic $0$ are: $\varpi_i$ with $i\in\{1,\ldots,n\}$ if $\grL=A_n$, $\varpi_n$ if $\grL=B_n$, $\varpi_1$ if $\grL=C_n$, and $\varpi_1$, $\varpi_{n-1}$, and $\varpi_n$ if $\grL=D_n$ with $n\Ge 3$ (see [Sa1], [Bou2, pp. 127--129], [De2], [Se1]).

\bigskip\noindent
{\bf 2.2. On reflex fields and Shimura types.}
Let $(G,X)$ be a Shimura pair. We have $\dbS(\dbR)=\dbG_{m,\dbR}(\dbC)$. We identify $\dbS(\dbC)=\dbG_{m,\dbR}(\dbC)\times\dbG_{m,\dbR}(\dbC)$ in such a way that the monomorphism $\dbS(\dbR)\hookrightarrow\dbS(\dbC)$ induces the map $z\to (z,\overline z)$ (here $z\in \dbG_{m,\dbR}(\dbC)$). Let $x\in X$; it is a homomorphism $x:\dbS\to G_{\dbR}$. Let $\mu_x:\dbG_{m,\dbC}\to G_{\dbC}$ be the cocharacter given on complex points by the rule $z\to x_{\dbC}(z,1)$. Let $\mu_0$ be a cocharacter of $G_{\overline{\dbQ}}$ such that we have $[\mu_{0,\dbC}]=[\mu_x]$. If $K$ is a field of characteristic $0$, the $G(\overline K)$-conjugacy class $[\mu_{0,\overline K}]$ of $\mu_{0,\overline K}$ will be denoted simply by $\mu_K$. If $K$ is a $p$-adic field, let $[\mu_p]:=[\mu_K]=[\mu_{0,\overline{K}}]$.

Let $\grC(G,X)$ be the set of cocharacters of $G_{\dbC}$ that are $G(\dbC)$-conjugates of the extensions to $\dbC$ of all Galois conjugates of $\mu_0$ (i.e., of all cocharacters of the form $\gamma(\mu_0):\dbG_{m,\overline{\dbQ}}\to G_{\overline{\dbQ}}$ for some element $\gamma\in\Gal(\overline{\dbQ})$). As the notation suggests, the set $\grC(G,X)$ does not depend on the choice of the cocharacter $\mu_0\in [\mu_{\dbQ}]=[\mu_0]$.

The {\it reflex field} $E(G,X)$ of $(G,X)$ is the field of definition of $[\mu_{\dbC}]=[\mu_x]$ (see [De1], [De2], [Mi2], and [Mi3]). It is the composite field of $E(G^{\ad},X^{\ad})$ and $E(G^{\ab},X^{\ab})$, cf. [De1, Prop. 3.8 (i)]. Thus, we now recall how one computes $E(G,X)$, if $G$ is either a torus or a simple, adjoint group. If $G$ is a torus, then $X=\{x\}$ and $E(G,X)$ is the field of definition of $\mu_x$ itself. From now on until Subsection 2.4 we will assume that $G$ is a simple, adjoint group over $\dbQ$; thus $(G,X)$ is a simple, adjoint Shimura pair.

Let $F(G,X)$ be a totally real number subfield of $\overline{\dbQ}\subseteq\dbC$ such that we have
$$G=\Res_{F(G,X)/\dbQ} G[F(G,X)],$$
with $G[F(G,X)]$ as an absolutely simple, adjoint group over $F(G,X)$ (cf. [De2, Subsubsect. 2.3.4]); the subfield $F(G,X)$ of $\overline{\dbQ}$ is unique up to $\Gal(\dbQ)$-conjugation. Let $I_0$ be the split, simple, adjoint group over $\dbQ$ of the same Lie type as every simple factor of $G_{\overline{\dbQ}}$. Let $Aut(I_0)$ be the group over $\dbQ$ of automorphisms of $I_0$. The group $I_0$ is the group of inner automorphisms of $I_0$ and thus it is a normal subgroup of $Aut(I_0)$. The quotient group $Aut(I_0)/I_0$ is either trivial or $\dbZ/2\dbZ$ or the symmetric group $S_3$. Let $I(G,X)$ be the smallest field extension of $F(G,X)$ such that the group $G[F(G,X)]_{I(G,X)}$ is an inner form of $I_{0,I(G,X)}$; the degree $[I(G,X):F(G,X)]$ divides the order of the finite group $Aut(I_0)/I_0$.

Let $T$ be a maximal torus of $G$. Let $B$ be a Borel subgroup of $G_{\overline{\dbQ}}$ that contains $T_{\overline{\dbQ}}$. Let $\grD$ be the Dynkin diagram of $\Lie(G_{\overline{\dbQ}})$ with respect to $T_{\overline{\dbQ}}$ and $B$. As $\break G_{\overline{\dbQ}}=\prod_{i\in\Hom(F(G,X),\dbR)} G[F(G,X)]\times_{F(G,X),i} \overline{\dbQ}$, we have a disjoint union
$$\grD=\cup_{i\in\Hom(F(G,X),\dbR)} \grD_i,$$
where $\grD_i$ is the Dynkin diagram of $\Lie(G[F(G,X)]\times_{F(G,X),i} \overline{\dbQ})$. For a vertex $\grn$ of $\grD$, let $\grg_{\grn}$ be the 1 dimensional Lie subalgebra of $\Lie(B)$ that corresponds to $\grn$.

The Galois group $\Gal(\dbQ)$ acts on $\grD$ as follows. If $\gamma\in\Gal(\dbQ)$, then $\gamma(\grn)$ is the vertex of $\grD$ defined by the identity $\grg_{\gamma(\grn)}=i_{g_{\gamma}}(\gamma(\grg_{\grn}))$, where $i_{g_{\gamma}}$ is the inner conjugation of $\Lie(G_{\overline{\dbQ}})$ defined by an element $g_{\gamma}\in G(\overline{\dbQ})$ that normalizes $T_{\overline{\dbQ}}$ and that has the property that $g_{\gamma}\gamma(B) g_{\gamma}^{-1}=B$. Here $\gamma(B)$ is the Galois conjugate of $B$ through $\gamma$.

Let $\grn_X$ be the set of vertices of $\grD$ such that the cocharacter of $T_{\overline{\dbQ}}$ that acts on $\grg_{\grn}$ trivially if $\grn\notin \grn_X$ and via the identical character of $\dbG_{m,\overline{\dbQ}}$ if $\grn\in\grn_X$, belongs to $[\mu_0]$. The image of $x:\dbS\to G_{\dbR}$ in a simple factor $\scrF$ of $G_{\dbR}$ is trivial if and only if the group $\scrF$ is compact (this is so as the centralizer of $\text{Im}(x)$ in $G_{\dbR}$ is a maximal compact subgroup of $G_{\dbR}$, cf. [De2, p. 259]). Thus, as the group $G_{\dbR}$ is not compact, the set $\grn_X$ is non-empty.

As the set $\grn_X$ is non-empty and as the Hodge $\dbQ$--structure on $\Lie(G)$ defined by each $x\in X$ is of type $\{(-1,1),(0,0),(1,-1)\}$, the group $I_0$ is of either classical or $E_6$ or $E_7$ Lie type (see [De2, Table 1.3.9]). Moreover, if $\grn_X$ contains a vertex of $\grD_i$, then (with the standard notations of [Bou1, Planches I to VI]) this vertex is: an arbitrary vertex if $I_0$ is of $A_n$ Lie type, vertex $1$ if $I_0$ is of $B_n$ Lie type, vertex $n$ if $I_0$ is of $C_n$ Lie type, an extremal vertex if $I_0$ is of $D_n$ Lie type, vertex $1$ or $6$ if $I_0$ is of $E_6$ Lie type, and vertex $7$ if $I_0$ is of $E_7$ Lie type (cf. [De2, Table 1.3.9]). The reflex field $E(G,X)$ is the fixed field of the open subgroup of $\Gal(\dbQ)$ that is the stabilizer of $\grn_X$, cf. [De2, Prop. 2.3.6].

The following three definitions conform with [De2] and [Mi3].

\medskip\noindent
{\bf 2.2.1. Definitions.}  {\bf (a)} If $I_0$ is of $A_n$, $B_n$, $C_n$, $E_6$, or $E_7$ Lie type, then we say $(G,X)$ is of $A_n$, $B_n$, $C_n$, $E_6$, or $E_7$ type (respectively).

Let now $I_0$ be of $D_n$ Lie type, with $n\Ge 4$. If $n=4$, we say $(G,X)$ is of $D_4^{\dbR}$ (resp. of $D_4^{\dbH}$) type if for each embedding $i:F(G,X)\hookrightarrow\dbR$, the orbit of $\grn_X$ under $\Gal(\dbQ)$ contains only one vertex (resp. only two vertices) of $\grD_i$. If $n\ge 5$, we say $(G,X)$ is of $D_n^{\dbR}$ type if for each embedding $i:F(G,X)\hookrightarrow\dbR$, the orbit of $\grn_X$ under $\Gal(\dbQ)$ contains only the vertex $1$ of $\grD_i$. If $n\ge 5$, we say $(G,X)$ is of $D_n^{\dbH}$ type if for each embedding $i:F(G,X)\hookrightarrow\dbR$, the orbit of $\grn_X$ under $\Gal(\dbQ)$ contains only vertices of $\grD_i$ that belong to the set $\{n-1,n\}$.  If $n\ge 4$, we say $(G,X)$ is of $D_n^{\text{mixed}}$ type if it is neither of $D_n^{\dbH}$ type nor of $D_n^{\dbR}$ type.

\smallskip
{\bf (b)} An adjoint Shimura pair is said to be of {\it abelian type} if it is a product of simple, adjoint Shimura pairs of $A_n$, $B_n$, $C_n$, $D_n^{\dbH}$, or $D_n^{\dbR}$ type.

\smallskip
{\bf (c)} Let $(G,X)$ be of $D_n^{\dbR}$ type, $n\Ge 4$. If $n=4$ we index the vertices of each $\grD_i$ in such a way that the $\Gal(\dbQ)$-orbit of $\grn_X$ contains only the vertex $1$ of $\grD_i$. Let $*$ be a field of characteristic $0$ such that the group $G_*$ is split. By the two half spin representations of a simple factor of $\Lie(G_*)$, we mean the two representations associated to the highest weights $\varpi_{n-1}$ and $\varpi_n$. By the spin representation of a simple factor of $\Lie(G_*)$, we mean the direct sum of the two half spin representations of the mentioned simple factor. Similarly, we speak about the two half spin representations and about the spin representation of the simply connected, semisimple group cover of a simple factor of $G_*$.

\medskip\noindent
{\bf 2.2.2. Extra fields.} There exists a unique totally real number field $E_1(G,X)$ such that $E(G,X)$ is either it or a totally imaginary quadratic extension of it, cf. [De1, Prop. 3.8 (ii)]. We say $(G,X)$ is: (i) {\it without involution} if $E(G,X)=E_1(G,X)$, and (ii) {\it with involution} if $E(G,X)\neq E_1(G,X)$. We have $E(G,X)=E_1(G,X)$ if and only if the set $\grn_X$ is fixed by the automorphism of $\overline{\dbQ}$ defined by the complex conjugation, cf. loc. cit. The last automorphism of $\overline{\dbQ}$ acts on each $\scrD_i$ via the standard opposition (it is either trivial or an involution) defined by the element $-w_0$ of [Bou1, item (XI), Planches I to IV], cf. [De1, Prop. 3.8 (ii)]. Thus, if $(G,X)$ is not of some $A_{2n+1}$ type, then from the description of elements $w_0$ we get that $(G,X)$ is with involution if and only if it is of one of the following types: $A_{2n}$, $D_{2n+3}^{\dbH}$, $D_{2n+3}^{\text{mixed}}$, or $E_6$. If $(G,X)$ is of $A_{2n+1}$ type, then it is without involution if and only if all simple, non-compact factors of $G_{\dbR}$ are $\text{\bf SU}(n+1,n+1)^{\ad}_{\dbR}$ groups (i.e., if and only if for each element $i\in\Hom(F(G,X),\dbR)$, the set $\grn_X$ contains either no vertex of $\grD_i$ or only the vertex $n+1$).

Let $EF_1(G,X)$ (resp. $EF(G,X)$) be the composite field of $F(G,X)$ with $E_1(G,X)$ (resp. with $E(G,X)$) inside $\overline{\dbQ}$. Let $G_1(G,X)$ (resp. $G(G,X)$) be the Galois extension of $\dbQ$ generated by $EF_1(G,X)$ (resp. by $EF(G,X)$). Let $G_2(G,X)$ be the Galois extension of $\dbQ$ generated by $F(G,X)$. Both number fields $G_1(G,X)$ and $G_2(G,X)$ are totally real.

\medskip\noindent
{\bf 2.2.3. A construction.}
For a totally real number field $F_1$ that contains $F(G,X)$ let
$$G^{F_1}:=\Res_{F_1/\dbQ} G[F(G,X)]_{F_1}.$$
We have a natural monomorphism $G[F(G,X)]\hookrightarrow \Res_{F_1/F(G,X)} G[F(G,X)]_{F_1}$ over $F(G,X)$ which at the level of $F(G,X)$-valued points is defined by the inclusion
$$G[F(G,X)](F(G,X))\hookrightarrow G[F(G,X)](F_1)=\Res_{F_1/F(G,X)} G[F(G,X)]_{F_1}(F(G,X)).$$
The $\Res_{F(G,X)/\dbQ}$ of this natural monomorphism is a monomorphism $G\hookrightarrow G^{F_1}$ that extends to an injective map
$$(G,X)\hookrightarrow (G^{F_1},X^{F_1})\leqno (1)$$
 of Shimura pairs. Here $X^{F_1}$ is the $G^{F_1}(\dbR)$-conjugacy class of the composite of any $x\in X$ with the monomorphism $G_{\dbR}\hookrightarrow G^{F_1}_{\dbR}$. As $G^{F_1}_{\dbR}\arrowsim G_{\dbR}^{[F_1:F(G,X)]}$, the hermitian symmetric domain $X^{F_1}$ is isomorphic to $X^{[F_1:F(G,X)]}$. We have $F(G^{F_1},X^{F_1})=F_1$.

Let $(\grD^{F_1},\grn_X^{F_1})$ be the analogue of the pair $(\grD,\grn_X)$  but for the Shimura pair $(G^{F_1},X^{F_1})$ and the pair $(T^{F_1},B^{F_1})$, where $T^{F_1}$ is the maximal torus of $G^{F_1}$ which is the centralizer of $T$ in $G^{F_1}$ and where $B^{F_1}$ is the unique Borel subgroup of $G^{F_1}_{\overline{\dbQ}}$ that contains both $T^{F_1}_{\overline{\dbQ}}$ and $B$. Let $s^{F_1}:\grD^{F_1}\twoheadrightarrow \grD$ be the $\Gal(\dbQ)$-invariant surjective map of Dynkin diagrams which is naturally associated to the monomorphism $G_{\overline{\dbQ}}\hookrightarrow G^{F_1}_{\overline{\dbQ}}\arrowsim G_{\overline{\dbQ}}^{[F_1:F(G,X)]}$; if $i:F(G,X)\hookrightarrow \dbR$ is an embedding, then $(s^{F_1})^{-1}(\grD_{i})$ is the disjoint union of the connected Dynkin diagrams of those simple factors of $\Lie(G^{F_1}_{\overline{\dbQ}})=\oplus_{i_1\in\Hom(F_1,\dbR)} \Lie(G^{F_1}\otimes_{F_1,i_1} \overline{\dbQ})$ that correspond to embeddings $i_1:F_1\hookrightarrow \dbR$ that extend $i$. We have $(s^{F_1})^{-1}(\grn_X)=\grn_X^{F_1}$, cf. the very definition of $X^{F_1}$.

 The identity $(s^{F_1})^{-1}(\grn_X)=\grn_X^{F_1}$ implies that the fixed fields of the stabilizers of $\grn_{X}$ and $\grn_X^{F_1}$ in $\Gal(\dbQ)$ are equal; thus $E(G,X)=E(G^{F_1},X^{F_1})$. The existence of $s^{F_1}$ also implies that the types of $(G,X)$ and $(G^{F_1},X^{F_1})$ are equal, cf. Definition 2.2.1 (a).

\medskip\noindent
{\bf 2.3. Definitions.} Let $(G,X)$ be a simple, adjoint Shimura pair. We say $(G,X)$ is:

\smallskip
{\bf (a)} of {\it inner} (resp. {\it non-inner}) {\it type}, if $I(G,X)$ is (resp. is not) $F(G,X)$;

\smallskip
{\bf (b)} with {\it apparent $\dbR$-involution}, if $EF(G,X)\neq F(G,X)$;

\smallskip
{\bf (c)} with {\it $\dbR$-involution}, if $E(G,X)$ is not a subfield of $G_2(G,X)$;

\smallskip
{\bf (d)} with {\it $\dbQ_p$-involution}, if there exists a finite field extension $F_{p,1}$ of $\dbQ_p$ with the property that the group pair $(G_{F_{p,1}},[\mu_p])$ has a simple factor $(G_1,[\mu_1])$ such that $G_1$ is a non-split, absolutely simple, unramified group over $F_{p,1}$ and $\mu_1$ is the extension to $\overline{\dbQ_p}$ of a cocharacter of $G_{1,F_{p,2}}$ whose $G_1(F_{p,2})$-conjugacy class is not fixed by $\Gal(F_{p,2}/F_{p,1})$; here $F_{p,2}$ is an unramified field extension of $F_{p,1}$ of degree 2 if $(G,X)$ is not of $D_4^{\text{mixed}}$ type and of some degree in the set $\{2,3\}$ if $(G,X)$ is of $D_4^{\text{mixed}}$ type;

\smallskip
{\bf (e)} of {\it compact type}, if the $F(G,X)$-rank of $G[F(G,X)]$ (i.e.,  the $\dbQ$--rank of $G$) is 0;

\smallskip
{\bf (f)} of {\it strong compact type}, if the $I(G,X)$-rank of $G[F(G,X)]_{I(G,X)}$ is 0;

\smallskip
{\bf (g)} of {\it $p$-compact type}, if the $\dbQ_p$-rank of $G_{\dbQ_p}$ is 0.

\medskip\noindent
{\bf 2.3.1. Remark.} In Definitions 2.3 (b) and (c) the use of ``$\dbR$-'' emphasizes that we are dealing with totally real number fields $F(G,X)$ and $G_2(G,X)$ while the word involution points out to the fact that suitable group pairs are in the natural sense with involution. The group pairs are not necessarily over $F(G,X)$ or $G_2(G,X)$; often we have to be in a $\dbQ_p$-context that is similar to the one of Definition 2.3 (d) (cf. Subsubsection 2.3.3 below).

\medskip\noindent
{\bf 2.3.2. On involutions.} We assume that $(G,X)$ is with involution i.e., $[E(G,X):E_1(G,X)]=2$. The field $G(G,X)$ is a Galois extension of $\dbQ$ that is a totally imaginary quadratic extension of $G_1(G,X)$. Thus $(G,X)$ is also with $\dbR$-involution. Standard facts on number fields (for instance, see [La, p. 168]) imply that the set of primes $p\in\dbN$ which split in $G_1(G,X)$, which do not split in $G(G,X)$, and which have the property that $G_{\dbQ_p}$ is unramified, has positive Dirichlet density. If $p$ is such a prime, then $G_{\dbQ_p}$ is a product of absolutely simple, unramified, adjoint groups and each prime of $E_1(G,X)$ (resp. of $E(G,X)$) that divides $p$ has residue field $\dbF_p$ (resp. $\dbF_{p^2}$). From this and [Mi3, Prop. 4.6 and Cor. 4.7] we easily get that $(G_{\dbQ_p},[\mu_p])$ has a simple factor $(G_1,[\mu_1])$ for which the following two properties hold:

\medskip
{\bf (a)} the group $G_1$ is non-split, absolutely simple, unramified, and splits over $B(\dbF_{p^2})$;

\smallskip
{\bf (b)} the cocharacter $\mu_1$ is the extension to $\overline{B(\dbF_{p^2})}$ of a cocharacter $\mu_{1,2}$ of $G_{1,B(\dbF_{p^2})}$ with the property that the $G_1(B(\dbF_{p^2}))$-conjugacy class $[\mu_{1,2}]$ is not fixed by $\Gal(B(\dbF_{p^2})/\dbQ_p)$.

\medskip
By taking $F_{p,1}$ and $F_{p,2}$ to be $\dbQ_p$ and $B(\dbF_{p^2})$ (respectively), from the property (b) we get that the Shimura pair $(G,X)$ is also with $\dbQ_p$-involution.

\medskip\noindent
{\bf 2.3.3. On $\dbR$-involutions.} We assume that $(G,X)$ is with $\dbR$-involution. As in Subsubsection 2.3.2 we argue that the set of primes $p\in\dbN$ which split in $G_2(G,X)$, which do not split in $G(G,X)$, and which have the property that $G_{\dbQ_p}$ is unramified, has positive Dirichlet density. For such a prime $p$, if $(G,X)$ is not (resp. is) of $D_4^{\text{mixed}}$ type we get the existence of a simple factor $(G_1,[\mu_1])$ of $(G_{\dbQ_p},[\mu_p])$ such that properties 2.3.2 (a) and (b) hold (resp. hold with $\dbF_{p^2}$ replaced by $\dbF_{p^s}$ for some $s\in\{2,3\}$). Thus $(G,X)$ is with $\dbQ_p$-involution. From the property 2.3.2 (a) we get that $G_1$ is not an inner form of $I_{0,\dbQ_p}$ and thus $(G,X)$ is also of non-inner type.

\medskip\noindent
{\bf 2.3.4. On the non-inner type and the $D_{2n}$ Lie type.} If $(G,X)$ is of $B_n$, $C_n$, or $E_7$ type, then $I(G,X)=F(G,X)$. If $n>1$ and $a\in\{0,\ldots,n+1\}$, then the group $\text{\bf SU}(a,n+1-a)_{\dbR}^{\ad}$ is not an inner form of $\text{\bf PGL}_{n+1,\dbR}$. Thus, as $F(G,X)$ is a totally real number field, $(G,X)$ is of non-inner type if it is of $A_{n}$ type with $n>1$. If $(G,X)$ is of $E_6$ type, then it is of non-inner type (cf. Subsubsections 2.2.2, 2.3.2, and 2.3.3). Thus the notion non-inner type is of interest only when the simple, adjoint Shimura pair $(G,X)$ is of $D_{2n}^{\dbH}$, $D_{2n}^{\dbR}$, or $D_{2n}^{\text{mixed}}$ type with $n\Ge 2$.

Let now $G$ be such that $I_0$ is of $D_{2n}$ Lie type, with $n\Ge 2$. We have $E(G,X)=E_1(G,X)$, cf. Subsubsection 2.2.2. Each simple factor of $G^{\ad}_{\dbR}$ is an inner form of $\text{\bf SO}(4n)^{\ad}_{\dbR}$ (cf. [De2, Subsubsect. 1.3.6]) and thus also of $I_{0,\dbR}$. Therefore the field $I(G,X)$ is totally real. If $i_0:F(G,X)\hookrightarrow\dbR$ is the embedding defined naturally by the inclusion $F(G,X)\hookrightarrow\dbC$, then $(G,X)$ is of inner type if and only if the action of $\Gal(F(G,X))$ on $\grD_{i_0}$ is trivial. Thus $(G,X)$ is of inner type if and only if the action of $\Gal(\dbQ)$ on the set of extremal vertices of $\grD$ has exactly three orbits.

If $(G,X)$ is of $D_{2n}^{\dbR}$ type, then the notions of apparent $\dbR$-involution and of non-inner type are unrelated (i.e., no one implies the other one). Also, if $(G,X)$ is of $D_{2n}^{\dbR}$ type, then it is never with $\dbR$-involution (cf. [De2, Rm. 2.3.12]) or with $\dbQ_p$-involution (this is trivial). If $(G,X)$ is of  non-inner $D_{2n}^{\dbH}$ or $D_{2n}^{\text{mixed}}$ type, then each element of $\Gal(F(G,X))$ that does not fix $I(G,X)$ will also not fix $E(G,X)$ and therefore $(G,X)$ is with apparent $\dbR$-involution.

\medskip\noindent
{\bf 2.3.5. Remark.} The Shimura pair $(G,X)$ is of compact type if and only if $\Sh(G,X)_{\dbC}$ is a pro-\'etale cover of a normal, projective scheme over $\dbC$; this is implied by [BHC, Thm. 12.3 and Cor. 12.4]. If $(G,X)$ is of $p$-compact type, then it is also of compact type. If $(G,X)$ is of $p$-compact type and if all primes of $F(G,X)$ that divide $p$ split in $I(G,X)$, then $(G,X)$ is also of strong compact type.

\medskip\noindent
{\bf 2.3.6. Example.} We assume that $[I(G,X):F(G,X)]=2$ and that $G[F(G,X)]$ is the adjoint group of the special unitary group $\text{\bf SU}_1(D,h_D)$, where $D$ is a central division algebra over $I(G,X)$ of dimension $(n+1)^2$ that is equipped with an involution $*$ which is of the second kind and for which we have an identity $F(G,X)=\{x\in I(G,X)|x^*=x\}$ and where $h_D$ is a non-degenerate hermitian form of index $0$ relative to $*$. Then $G[F(G,X)]_{I(G,X)}$ is $\text{\bf GL}_1(D)^{\ad}$ and thus has rank $0$. Thus the Shimura pair $(G,X)$ is of strong compact type.

\medskip\noindent
{\bf 2.3.7. Example.} We assume that $G_{\dbR}$ has at least one simple, compact factor and $(G,X)$ is of strong compact type. Let $F_1$ be a totally real finite field extension of $F(G,X)$ such that the $F_{1,I}$-rank of $G[F(G,X)]_{F_{1,I}}$ is positive, where $F_{1,I}$ is the composite field of $F_1$ and $I(G,X)$. We have $I(G^{F_1},X^{F_1})=F_{1,I}$. The Shimura pair $(G^{F_1},X^{F_1})$ is of compact type (as $G^{F_1}_{\dbR}\arrowsim G_{\dbR}^{[F_1:F(G,X)]}$ has at least one simple, compact factor) but is not of strong compact type (as the $F_{1,I}$-rank of $G[F(G,X)]_{F_{1,I}}$ is positive).

\bigskip\noindent
{\bf 2.4. Hodge quasi products.}
Let $I:=\{1,\ldots,n\}$. For each element $i\in I$ we consider an injective map
$$f_i:(G_i,X_i)\hookrightarrow (\text{\bf GSp}(W_i,\psi_i),S_i)$$
of Shimura pairs. Let $\scrB_i:=\End(W_i)^{G_i(\dbQ)}$. We say $f_i$ is a {\it PEL type embedding} if the group $G_i$ is the identity component of the subgroup of $\text{\bf GSp}(W_i,\psi_i)$ that fixes $\scrB_i$. A Shimura pair that admits a PEL type embedding into a Shimura pair that defines a Siegel modular variety, is called a {\it Shimura pair of PEL type} (cf. [Sh], [De1], [Zi], [Ko2]).

What follows next extends [Va1, Subsect. 2.5, Ex. 5] (loc. cit. dealt with $n=2$). Let
$$(W,\psi):=\oplus_{i\in I} (W_i,\eps_i\psi_i),$$
where $\eps_i\in\{-1,1\}$. Let $G_I:=\prod_{i\in I} G_i$. It is a subgroup of both $\prod_{i\in I} \text{\bf GSp}(W_i,\psi_i)$ and $\text{\bf GL}_W$. Let $G_0$ be the identity component of the intersection $G_I\cap \text{\bf GSp}(W,\psi)$. Let $x:\dbS\hookrightarrow G_{I,\dbR}$ be a monomorphism that defines an element of $\prod_{i\in I} X_i$; it factors through $G_{0,\dbR}$. We choose the $\eps_i$'s (with $i\in I$) such that either $+2\pi i\psi$ or $-2\pi i\psi$ is a polarization of the Hodge $\dbQ$--structure on $W$ defined by $x$. Thus, if $X_0$ (resp. if $S$) is the $G_0(\dbR)$-conjugacy (resp. the $\text{\bf GSp}(W,\psi)(\dbR)$-conjugacy) class of $x$, we get an injective map
$$f_0:(G_0,X_0)\hookrightarrow (\text{\bf GSp}(W,\psi),S)$$
 of Shimura pairs which we call a Hodge quasi product of $f_i$'s. We emphasize that both $\pm\psi$ and $f_0$ depend uniquely on the choice of $x$. We denote each such $f_0$ by $\times^{\scrH}_{i\in I} f_i$. We refer to the Shimura pair $(G_0,X_0)$ as a {\it Hodge quasi product} of $(G_i,X_i)$'s and we denote it by $\times^{\scrH}_{i\in I} (G_i,X_i)$. If $n=2$, then we also denote $f_0$ by $f_1\times^{\scrH} f_2$ and $(G_0,X_0)$ by $(G_1,X_1)\times^{\scrH} (G_2,X_2)$.

We call a $\dbZ$-lattice $L$ of $W$ well adapted for $f_0$ if it is of the form $\oplus_{i\in I} L_i$, where for each $i\in I$ the $\dbZ$-lattice $L_i$ of $W_i$ is such that we get a perfect alternating form $\psi_i:L_i\otimes_{\dbZ} L_i\to\dbZ$.

\medskip\noindent
{\bf 2.4.1. Facts. (a)} {\it The group $G_I$ is the extension of $\dbG_{m,\dbQ}^{n-1}$ by $G_0$. We have $Z^0(G_0)=\dbG_{m,\dbQ}$ if and only if for all $i\in I$ we have $Z^0(G_i)=\dbG_{m,\dbQ}$.}

\smallskip
{\bf (b)} {\it We assume that for each $i\in I$ the injective map $f_i$ is a PEL type embedding. Then $f_0$ is also a PEL type embedding.}

\medskip
\proof
As $\prod_{i\in I} \text{\bf GSp}(W_i,\psi_i)$ is the extension of $\dbG_{m,\dbQ}^{n-1}$ by $\text{\bf GSp}(W,\psi)\cap (\prod_{i\in I} \text{\bf GSp}(W_i,\psi_i))$, the first part of (a) follows. The image of each $x_i\in X_i$ contains $Z(\text{\bf GL}_{W_i\otimes_{\dbQ} \dbR})$. Thus $Z^0(G_i)$ contains the $1$ dimensional split torus $Z(\text{\bf GL}_{W_i})$. The first part of (a) implies that we have a short exact sequence $0\to Z^0(G_0)\to\prod_{i\in I} Z^0(G_i)\to \dbG_{m,\dbQ}^{n-1}\to 0$. From the last two sentences we get that the second part of (a) also holds.

Part (b) is implied by the fact that $G_0$ is the identity component of the subgroup of $\text{\bf GSp}(W,\psi)$ that fixes the semisimple $\dbQ$--subalgebra $\prod_{i\in I} \scrB_i$ of $\End(W)$.\endproof

\bigskip\smallskip
\noindent
{\boldsectionfont 3. Basic techniques and the proof of Theorem 1.3.1}
\bigskip

In Subsections 3.1 to 3.3 we review basic facts about Conjecture 1.1.1. In Subsection 3.4 we prove Theorem 1.3.1. Until Section 4 we use the notations of Subsection 1.1. Thus $G_{\dbQ_p}$ is a connected subgroup of $\text{\bf GL}_{L_A\otimes_{\dbZ} \dbQ_p}$.

\bigskip\noindent
{\bf 3.1. Simple properties.} {\bf (a)} The Mumford--Tate conjecture holds for $A$ if and only if it holds for an abelian variety isogenous to $A_{E^\prime}$, where $E^{\prime}$ is a finite field extension of $E$.

\smallskip
{\bf (b)} The following three assertions are equivalent:

\medskip
{\bf (i)} The abelian variety $A_{\overline{E}}$ has complex multiplication (i.e., the semisimple $\dbQ$--algebra $\End(A_{\overline{E}})\otimes_{\dbZ} \dbQ$ contains a commutative, semisimple $\dbQ$--algebra of dimension $2\dim(A)$).

\smallskip
{\bf (ii)} The group $H_A$ is a torus.

\smallskip
{\bf (iii)} The group $G_{\dbQ_p}$ is a torus.

\medskip
The equivalence of (i) (resp. (iii)) and (ii) follows from  definitions (resp. from Theorem 1.1.4).

\smallskip
{\bf (c)} As $Z(G_{\dbQ_p})$ is a subgroup of $Z(H_{A,\dbQ_p})$ (cf. Theorems 1.1.3 and 1.1.4), the monomorphism $G_{\dbQ_p}\hookrightarrow H_{A,\dbQ_p}$ gives birth naturally to a monomorphism $G^{\ad}_{\dbQ_p}\hookrightarrow H_{A,\dbQ_p}^{\ad}$.

\bigskip\noindent
{\bf 3.2. Frobenius tori.}
For the constructions reviewed here we refer to [Ch]. Let $v$, $A_v$, and $l$ be as before Theorem 1.1.5. Until Subsection 3.4 we will assume that $p\neq l$. Let $\overline{v}$ be a prime of $\overline{E}$ that divides $v$. Let $F(\overline{v})\in\Gal(E)$ be a Frobenius automorphism of $\overline{v}$. The pair $(\overline{v},F(\overline{v}))$ is uniquely determined up to $\Gal(E)$-conjugation. As $v$ is a prime of good reduction for $A$, the restriction of $\rho$ to the inertia group of $\overline{v}$ is trivial (cf. [BLR, Ch. 7, 7.4, Thm. 5]). Thus the element $\rho(F(\overline{v}))\in \text{\bf GL}_{L_A\otimes_{\dbZ} \dbQ_p}(\dbQ_p)$ is well defined regardless of the choice of $F(\overline{v})$; it and its positive, integral powers will be referred as {\it Frobenius elements}. As $\dbZ_p$-modules, via $\overline{v}$ we identify $L_A\otimes_{\dbZ} \dbZ_p=T_p(A)$ with the Tate module $T_p(A_{v,\overline{k(v)}})$ of $A_{v,\overline{k(v)}}$; thus the Frobenius endomorphism $Fr_v$ of $A_v$ acts on $L_A\otimes_{\dbZ} \dbQ_p$ as $\rho(F(\overline{v}))$ does.

Let $\scrI_v:=\text{\bf GL}_1(\End(A_v)\otimes_{\dbZ}\dbQ)$. Let $T_v$ be the identity component of the smallest subgroup of $\scrI_v$ that has $Fr_v$ as a $\dbQ$--valued point. It is a $\dbQ$--model of the identity component of the Zariski closure in $\text{\bf GL}_{L_A\otimes_{\dbZ} \dbQ_p}$ of the subgroup $\{\rho(F(\overline{v})^s)|s\in\dbZ\}$ of $\text{\bf GL}_{L_A\otimes_{\dbZ} \dbQ_p}(\dbQ_p)$. As $Fr_v$ is a semisimple endomorphism of $L_A\otimes_{\dbZ} \dbQ_p$ (see [Tat]), $T_v$ is a torus. Thus $T_{v,\dbQ_p}$ is a torus of $G_{\dbQ_p}$, uniquely determined up to $\im(\rho)$-conjugation. The torus $T_v$ depends only on some positive, integral power $\pi_v$ of $Fr_v$, cf. [Ch, Subsects. 3.a and 3.b]; here we can replace $Fr_v$ by the Frobenius endomorphism of any other abelian variety $B_{v^\prime}$ over a finite field extension $k(v^\prime)$ of $k(v)$ that is isogenous to $A_{v,k(v^\prime)}$.

If the finite field $k(v^\prime)$ is large enough, then we can choose $B_{v^\prime}$ such that there exists an abelian scheme $B_R$ over a complete discrete valuation ring $R$ of mixed characteristic $(0,l)$ and of residue field $k(v^\prime)$, which lifts $B_{v^\prime}$ and has complex multiplication. Let $B:=B_R\times_R R[{1\over l}]$. We can assume that $E$ is a subfield of $R[{1\over l}]$. We fix an embedding $R[{1\over l}]\hookrightarrow\dbC$ that extends the embedding $i_E$ of Subsection 1.1. Let $H_B$ be the Mumford--Tate group of $B_{\dbC}$; as $B_R$ has complex multiplication, $H_B$ is a torus. Let $C_B$ be the centralizer of $H_B$ in $\text{\bf GL}_{H_1(B(\dbC),\dbQ)}$. Let $D_B$ be the centralizer of $C_B$ (i.e., the double centralizer of $H_B$) in $\text{\bf GL}_{H_1(B(\dbC),\dbQ)}$. As $B_R$ has complex multiplication and as $\Lie(C_B)$ is naturally identified with the Lie algebra of $\dbQ$--endomorphisms of $B_{\dbC}$, we have a natural Lie monomorphism  $\Lie(C_B)\hookrightarrow \End(B_{v^\prime})\otimes_{\dbZ} \dbQ$. Each invertible element of the $\dbQ$--subalgebra of $\End(B_{v^\prime})\otimes_{\dbZ} \dbQ$ generated by the Frobenius endomorphism of $B_{v^\prime}$, is naturally identified with an element of the center of $\Lie(C_B)$ and thus also with a $\dbQ$--valued point of $D_B$; thus $T_{v}$ is naturally identified with a subtorus of $D_B$. If $B_R$ is equipped with a polarization $\lambda_{B_R}$ and if $\scrB$ is a semisimple $\dbQ$--subalgebra of $\End(H_1(B(\dbC),\dbQ))^{H_B(\dbQ)}$, then $T_v$ is in fact a torus of the connected component of the intersection of $\text{\bf GSp}(H_1(B(\dbC),\dbQ),\psi_{B})$ with the centralizer of $\scrB$ in $\text{\bf GL}_{H_1(B(\dbC),\dbQ)}$. Here $\psi_B$ is the alternating form on $H_1(B(\dbC),\dbQ)$ induced by $\lambda_{B_R}$.

The following Lemma is due in essence to Bogomolov and Serre, cf. [Bog] and [Se3, pp. 31, 43, and 44].

\bigskip\noindent
{\bf 3.3. Lemma.} {\it There exists a subtorus $T^0$ of $Z^0(H_A)$ that does not depend on $p$ and such that $G_{\dbQ_p}$ is the subgroup of $H_{A,\dbQ_p}$ generated by $G^{\der}_{\dbQ_p}$ and $T^0_{\dbQ_p}$.}

\medskip
\proof
By replacing $E$ with a finite field extension of it, we can assume that $\End(A)=\End(A_{\overline{\dbQ}})$. Let $C_A$ be the centralizer of $Z^0(H_A)$ in $\text{\bf GL}_{L_A\otimes_{\dbZ} \dbQ}$. The group $C_A$ is reductive and $Z(C_A)$ is a torus. The commutative, semisimple $\dbQ$--algebra $Z(C_A)(\dbQ)$ is a $\dbQ$--subalgebra of
$$\End(L_A\otimes_{\dbZ} \dbQ_p)=\End(T_p(A)\otimes_{\dbZ_p} \dbQ_p)=\End(T_p(A_{v,\overline{k(v)}}))$$
which is generated by $\dbQ$--linear combinations of \'etale $\dbQ_p$-realizations of endomorphisms of $A$ and thus also of $A_v$. Therefore $Z(C_A)$ is naturally a torus of $\scrI_v$. Let $C_v$ be the centralizer of $Z(C_A)$ in $\scrI_v$; it is a reductive group over $\dbQ$ (cf. [Bor, Ch. IV, 13.17, Cor. 2 (a)]). We identify $\scrI_{v,\dbQ_p}$ and $C_{v,\dbQ_p}$ with subgroups of $\text{\bf GL}_{L_A\otimes_{\dbZ} \dbQ_p}$ that contain $Z(C_A)_{\dbQ_p}$.

We have $\pi_v\in C_v(\dbQ)$. We denote also by $\pi_v$ its \'etale $\dbQ_p$-realization (i.e., the corresponding integral power of $\rho(F(\overline{v}))$). Thus $\pi_v\in G_{\dbQ_p}(\dbQ_p)\leqslant H_A(\dbQ_p)\leqslant C_A(\dbQ_p)$. After replacing $\pi_v$ with a positive, integral power of it that divides the order of the finite group $C_A^{\der}\cap Z(C_A)$, we can write $\pi_v\in C_A(\dbQ_p)$ as a product
$$\pi_v=\pi_0\pi_1,$$
where $\pi_0\in Z(C_A)(\dbQ_p)$ and $\pi_1\in C_A^{\der}(\dbQ_p)$.
As $\pi_v\in H_A(\dbQ_p)$ and as $H_A$ is a reductive subgroup of $C_A$ such that $Z^0(H_A)$ is a subtorus of $Z(C_A)$, after a similar replacement of $\pi_v$ we get that $\pi_0\in Z^0(H_A)(\dbQ_p)$ and $\pi_1\in H_A^{\der}(\dbQ_p)$. Let $C_v^\prime$ be the unique reductive subgroup of $C_v$ such that we have a natural isogeny $C_v^\prime\times_{\dbQ} Z(C_A)\to C_v$. As $\pi_v\in C_v(\dbQ)$, after a similar replacement of $\pi_v$ we get that $\pi_0\in Z(C_A)(\dbQ)$ and $\pi_1\in C_v^\prime(\dbQ)$.

Thus $\pi_0\in Z^0(H_A)(\dbQ_p)\cap Z(C_A)(\dbQ)=Z^0(H_A)(\dbQ)$. Let $T^0$ be the identity component of the smallest subgroup of $Z^0(H_A)$ which has $\pi_0$ as a $\dbQ$--valued point; it is a torus.

We now take $v$ such that $T_{v,\dbQ_p}$ is a maximal torus of $G_{\dbQ_p}$, cf. Theorem 1.1.5. The images of $Z^0(G_{\dbQ_p})$, $G_{\dbQ_p}$, and $T_{v,\dbQ_p}$ in $H_{A,\dbQ_p}^{\ab}$ are the same torus $T^1_{\dbQ_p}$. As the images of $\pi_v$ and $\pi_0$ in $H_{A,\dbQ_p}^{\ab}(\dbQ_p)$ coincide, $T^1_{\dbQ_p}$ is the extension to $\dbQ_p$ of $T^1:=\text{Im}(T^0\to H_A^{\ab})$. Thus $Z^0(G_{\dbQ_p})$ is a subtorus of $Z^0(H_A)_{\dbQ_p}$ (cf. Theorems 1.1.2 and 1.1.4) that has the same image in $H_{A,\dbQ_p}^{\ab}$ as $T^0_{\dbQ_p}$. This implies that $Z^0(G_{\dbQ_p})=T^0_{\dbQ_p}$. Thus $G_{\dbQ_p}$ is the subgroup of $H_{A,\dbQ_p}$ generated by $G^{\der}_{\dbQ_p}$ and $T^0_{\dbQ_p}$.

As $\pi_0\in Z(C_A)(\dbQ)$ depends only on $v$, the torus $T^0$ does not depend on the prime $p\neq l$. The fact that $T^0$ also does not depend on $l$, is checked in the standard way by considering another finite prime $v_1$ of $E$ which is of good reduction for $A$, whose characteristic is not $l$, and whose Frobenius torus $T_{v_1}$ has the same rank as $T_v$. \endproof

\medskip\noindent
{\bf 3.3.1. Remark.}
If $T_{v,\dbQ_p}$ is not a maximal torus of $G_{\dbQ_p}$, then the same argument shows that the image of $T_{v,\dbQ_p}$ in $H_{A,\dbQ_p}^{\ab}$ is the extension to $\dbQ_p$ of a subtorus of $H_A^{\ab}$.

\bigskip\noindent
{\bf 3.4. Proof of  1.3.1.}
We prove  Theorem 1.3.1. Let $H_A^0$ be the subgroup of $H_A$ generated by $H_A^{\der}$ and the torus $T^0$ of Lemma 3.3. As $Z^0(G_{\dbQ_p})=T^0_{\dbQ_p}$, to prove Theorem 1.3.1 it suffices to show that we have an identity $T^0=Z^0(H_A)$. We include a proof of the identity $T^0=Z^0(H_A)$ that is also used to get several other properties which play key roles in the subsequent Sections. In this Subsection we do not assume that $p\neq l$.

For Hodge cycles on abelian varieties over reduced $\dbQ$--schemes we refer to [De3]. It suffices to prove the identity $T^0=Z^0(H_A)$ under the extra assumptions that $E^{\text{conn}}=E$ and that all Hodge cycles on $A_{\dbC}$ are also Hodge cycles on $A$, cf. [De3, Prop. 2.9 and Thm. 2.11]. As the torus $T^0$ does not depend on $p$, until Section 4 we will choose the prime $v$ to be unramified over $l$. We identity $E_v$ with $B(k(v))$. Let $A_{W(k(v))}$ be the abelian scheme over $W(k(v))$ whose generic fibre is $A_{B(k(v))}$. Let $(M,\phi)$ be the $F$-crystal (i.e., the contravariant Dieudonn\'e module) of the special fibre $A_v$ of $A_{W(k(v))}$. Thus $M$ is a free $W(k(v))$-module of rank $2\dim(A)$ and $\phi$ is a Frobenius-linear monomorphism of $M$. Let $F^1$ be the Hodge filtration of $M=H^1_{\text{crys}}(A_v/W(k(v)))=H^1_{\text{dR}}(A_{W(k(v))}/W(k(v)))$ defined by $A_{W(k(v))}$. We identify $L_A^*\otimes_{\dbZ} \dbQ_l$ with $H^1_{\acute et}(A_{\overline{B(k(v))}},\dbQ_l)$ as $\dbQ_l$-vector spaces.

For a vector space $V$ over a field $K$ let $\scrT(V):=\oplus_{s,t\in\dbN\cup\{0\}} V^{\otimes s}\otimes_K V^{*\otimes t}$. We identify $\End(V)=V\otimes_K V^*\subseteq\scrT(V)$. The action of $\phi$ on $M[{1\over l}]$ extends in the standard (tensorial) way to an action of $\phi$ on $\scrT(M[{1\over l}])$. Let $(v_{\alpha})_{\alpha\in\scrJ}$ be a family of tensors of $\scrT(L_A^*\otimes_{\dbZ} \dbQ_l)$ such that the reductive group $G_{\dbQ_l}$ is the subgroup of $\text{\bf GL}_{L_A\otimes_{\dbZ} \dbQ_l}$ (equivalently of $\text{\bf GL}_{L_A^*\otimes_{\dbZ} \dbQ_l}$) that fixes $v_{\alpha}$ for all $\alpha\in\scrJ$, cf. [De3, Prop. 3.1 (c)].

Let $t_{\alpha}$ be the tensor of $\scrT(M[{1\over l}])$ that corresponds to $v_{\alpha}$ via Fontaine comparison theory (see [Fo] and [Ts]). It is fixed by the action of $\phi$ on $\scrT(M[{1\over l}])$ and moreover it belongs to the $F^0$-filtration of the tensor product filtration of $\scrT(M[{1\over l}])$ defined naturally by the filtration $(F^i[{1\over l}])_{i\in\{0,1,2\}}$ of $M[{1\over l}]$, where $F^2:=0$ and $F^0:=M$, and by the filtration of $M[{1\over l}]^*$ that is the dual of $(F^i[{1\over l}])_{i\in\{0,1,2\}}$.

Let $\tilde G_{B(k(v))}$ be the subgroup of $\text{\bf GL}_{M[{1\over l}]}$ that fixes $t_{\alpha}$ for all $\alpha\in\scrJ$; it is a form of $G_{B(k(v))}$ and thus it is a reductive group. Based on Theorem 1.1.2, we can assume that there exists a subset $\scrJ_1$ of $\scrJ$ such that $(v_{\alpha})_{\alpha\in\scrJ_1}$ is the family of all $\dbQ_l$-\'etale realizations of Hodge cycles on $A_{B(k(v))}$ that belong to $\scrT(L_A^*\otimes_{\dbZ} \dbQ_l)$ (i.e., it is the family of all tensors of $\scrT(L_A^*\otimes_{\dbZ} \dbQ)$ fixed by $H_A$). From this and [Bl, Thm. (0.3)] we get that $(t_{\alpha})_{\alpha\in\scrJ_1}$ is the family of all de Rham realizations of Hodge cycles on $A_{B(k(v))}$ that belong to $\scrT(M[{1\over l}])$. Let $\tilde H_{A,B(k(v))}$ be the subgroup of $\text{\bf GL}_{M[{1\over l}]}$ that fixes $t_{\alpha}$ for all $\alpha\in\scrJ_1$; it is a ``de Rham'' form of $H_{A,B(k(v))}$. The existence of an isomorphism between the two pairs $(L_A^*\otimes_{\dbZ} B(k(v)),(v_{\alpha})_{\alpha\in\scrJ})$ and $(M[{1\over l}],(t_{\alpha})_{\alpha\in\scrJ})$ is measured by a class $\delta_v$ in the pointed set $H^1(B(k(v)),G_{B(k(v))})$, cf. Fontaine comparison theory. As $G_{B(k(v))}$ acts trivially on $T^0_{B(k(v))}=Z^0(G_{B(k(v))})$ (via inner conjugation), the image of $\delta_v$ in $H^1(B(k(v)),Aut(T^0_{B(k(v))}))$ is trivial. Thus we can identify $T^0_{B(k(v))}=Z^0(G_{B(k(v))})$ with the subgroup $Z^0(\tilde G_{B(k(v))})$ of either $\tilde G_{B(k(v))}$ or $\tilde H_{A,B(k(v))}$.

Let $\mu_v:\dbG_{m,W(k(v))}\to \text{\bf GL}_M$ be the inverse of the canonical split cocharacter of $(M,F^1,\phi)$ defined in [Wi, p. 512]. We have a direct sum decomposition $M=F^1\oplus \tilde F^0$ such that $\tilde F^0$ is fixed by $\mu_v$ and $\beta\in\dbG_{m,W(k(v))}(W(k(v)))$ acts via $\mu_v$ on $F^1$ as the multiplication with $\beta^{-1}$. For $\alpha\in\scrJ$, the tensor $t_{\alpha}$ is fixed by $\mu_v$ (cf. the functorial aspects of [Wi, p. 513]). Thus $\mu_{v,B(k(v))}$ factors through $\tilde G_{B(k(v))}$.

Let $i_{\overline{\dbQ_l}}:\overline{\dbQ_l}\hookrightarrow\dbC$ be an embedding that extends the embedding $i_E$ of Subsection 1.1. We view $\dbC$ as a $B(k(v))$-algebra via the restriction $i_{B(k(v))}$ of $i_{\overline{\dbQ_l}}$ to $B(k(v))$. The canonical identification $H^1(A(\dbC),\dbZ)=L_A^*$ allows us to view $H_A$ and $H_{A,\dbC}$ naturally as subgroups of the linear groups $\text{\bf GL}_{H^1(A(\dbC),\dbQ)}$ and $\text{\bf GL}_{H^1(A(\dbC),\dbC)}$ (respectively). Let $P$ be the parabolic subgroup of $\tilde H_{A,B(k(v))}$ that is the normalizer of $F^1[{1\over l}]$. Under the standard identification of $M_{\dbC}:=M[{1\over l}]\otimes_{B(k(v))} \dbC=H^1_{\text{dR}}(A_{\dbC}/\dbC)$ with $H^1(A(\dbC),\dbC)$, the groups $\tilde H_{A,\dbC}$ and $P_{\dbC}$ get identified with $H_{A,\dbC}$ and respectively with the parabolic subgroup of $H_{A,\dbC}$ that normalizes the direct summand $F^{1,0}=F^1\otimes_{W(k(v))} \dbC$ of the Hodge decomposition $M_{\dbC}=F^{1,0}\oplus F^{0,1}$ of $A_{\dbC}$. Both $\mu_A$ and $\mu_{v,\dbC}$ are (identified with) cocharacters of $P_{\dbC}$.

Let $g\in P(\dbC)$ be such that $g\mu_A g^{-1}$ commutes with $\mu_{v,\dbC}$. As the commuting cocharacters $g\mu_A g^{-1}$ and $\mu_{v,\dbC}$ act in the same way on $F^{1,0}$ and on $M_{\dbC}/F^{1,0}$, they coincide. As $\mu_{v,B(k(v))}$ factors through $\tilde G_{B(k(v))}$ and as $T^0_{B(k(v))}=Z^0(\tilde G_{B(k(v))})$, the cocharacter $\mu_{v,\dbC}=g\mu_A g^{-1}$ factors through the normal subgroup $H_{A,\dbC}^0$ of $H_{A,\dbC}$. From the last two sentences, we get that both $g\mu_Ag^{-1}$ and $\mu_A$ factor through $H_{A,\dbC}^0$. Due to the smallest property of $H_A$ (see Subsection 1.1) and the fact that $\mu_A$ factors through $H_{A,\dbC}^0$, we get that $H_A=H_A^0$. Thus $T^0=Z^0(H_A)$. This ends the proof of Theorem 1.3.1.\endproof

\medskip\noindent
{\bf 3.4.1. Simple properties.}
{\bf (a)} Conjecture 1.1.1 holds if $A$ has complex multiplication over $\overline{\dbQ}$. This well known fact also follows from the property 3.1 (b) and Theorems  1.1.2 and 1.3.1.

\smallskip
{\bf (b)} The monomorphism $\tilde G_{\overline{\dbQ_l}}\hookrightarrow\tilde H_{A,\overline{\dbQ_l}}$ is naturally identified up to $H_A(\overline{\dbQ_l})$-conjugation with the monomorphism $G_{\overline{\dbQ_l}}\hookrightarrow H_{A,\overline{\dbQ_l}}$, cf. the existence of the class $\delta_v$. As $\mu_v$ factors through $\tilde G_{B(k(v))}$ and as the cocharacters $\mu_{v,\dbC}$ and $\mu_A$ of $H_{A,\dbC}$ are $H_A(\dbC)$-conjugate, there exists a cocharacter $\mu^{\acute et}_{v,\overline{\dbQ_l}}:\dbG_{m,\overline{\dbQ_l}}\to G_{\overline{\dbQ_l}}$ whose extension to $\dbC$ (via $i_{\overline{\dbQ_l}}:\overline{\dbQ_l}\hookrightarrow\dbC$) is $H_A(\dbC)$-conjugate to $\mu_A$.

\smallskip
{\bf (c)} We use the notations of Subsection 3.4. The element $\phi^{[k(v):\dbF_l]}\in \tilde G_{B(k(v))}(B(k(v)))$ is the crystalline realization of $Fr_v$ and thus it is a semisimple element. This implies that the torus $T_{v,B(k(v))}$ is naturally identified with the smallest subtorus of $\tilde G_{B(k(v))}$ that has $\phi^{s[k(v):\dbF_l]}$ as a $B(k(v))$-valued point, where $s\in\dbN$ is such that $\pi_v=Fr_v^s$. Thus $T_{v,\overline{\dbQ_l}}$ is naturally identified with a torus of $\tilde G_{\overline{\dbQ_l}}$ and therefore (cf. (b)) also of $G_{\overline{\dbQ_l}}$.

\medskip\noindent
{\bf 3.4.2. Corollary.} {\it We assume that $E^{\text{conn}}=E$, that $v$ is unramified over $l$, that $k(v)=\dbF_l$, and that $T_v$ has the same rank as $G_{\dbQ_l}$. We also assume that $T_v$ is (isomorphic to) a torus of $H_A$ in such a way that there exists an element $g_v\in H_A(\overline{\dbQ_l})$ with the property that the torus $g_vT_{v,\overline{\dbQ_l}}g_v^{-1}$ of $H_{A,\overline{\dbQ_l}}$ is a maximal torus of $G_{\overline{\dbQ_l}}$. Let $i_{\overline{\dbQ_l}}:\overline{\dbQ_l}\hookrightarrow \dbC$ be an embedding that extends the embedding $i_E:E\hookrightarrow\dbC$ of Subsection 1.1. Let $\grC_l(H_A,X_A)$ be the set of cocharacters of $H_{A,\overline{\dbQ_l}}$ whose extensions to $\dbC$ via $i_{\overline{\dbQ_l}}$ belong to the set $\grC(H_A,X_A)$ of cocharacters of $H_{A,\dbC}$ defined in Subsection 2.2. Then the torus $T_{v,\overline{\dbQ_l}}$ of $H_{A,\overline{\dbQ_l}}$ is generated by images of certain cocharacters in the set $\grC_l(H_A,X_A)$.}

\medskip
\proof
We can choose the element $g_v$ such that the cocharacter $\mu^{\acute et}_{v,\overline{\dbQ_l}}$ of the property 3.4.1 (b) factors through $g_vT_{v,\overline{\dbQ_l}} g_v^{-1}$. As the cocharacters $\mu_{v,\dbC}$ and $\mu_A$ are $P(\dbC)$-conjugate (see Subsection 3.4) and thus also $H_A(\dbC)$-conjugate, the set $\grC_l(H_A,X_A)$ is the set of $H_A(\overline{\dbQ_l})$-conjugates of the $\Gal(\dbQ)$-conjugates of the cocharacter $g_v^{-1}\mu^{\acute et}_{v,\overline{\dbQ_l}}g_v$ of the extension of $T_v$ to $\overline{\dbQ_l}$ (equivalently to $\overline{\dbQ}$). But up to $G_{\dbQ_l}(\overline{\dbQ_l})$-conjugation, we can choose the cocharacter $\mu^{\acute et}_{v,\overline{\dbQ_l}}:\dbG_{m,\overline{\dbQ_l}}\to G_{\overline{\dbQ_l}}$ such that $g_v^{-1}\mu^{\acute et}_{v,\overline{\dbQ_l}}g_v$ is a strong Hodge cocharacter of $T_{v,\overline{\dbQ_l}}$ in the sense of [Pi, Def. 3.14] (see [Pi, Prop. 3.13 and paragraph before Def. 3.14]). Thus the fact that $T_{v,\overline{\dbQ_l}}$ is generated by images of certain cocharacters in the set $\grC_l(H_A,X_A)$ is only a weaker form of [Pi, Thm. 3.16]. \endproof

\bigskip\smallskip
\noindent
{\boldsectionfont 4. Injective maps into unitary Shimura pairs of PEL type}
\bigskip

In Subsections 4.1 to 4.9 we include different rational versions of the embedding results of [Sa1]; the approach is close in spirit to [Sa2], [De2, Prop. 2.3.10], and [Va1, Subsects. 6.5 and 6.6]. Corollary 4.10 classifies all adjoint Shimura pairs that are adjoints of Shimura pairs of PEL type. Until Section 5 the notations are independent of the ones of Subsections 1.1, 1.3, or 1.5.

Until Subsection 4.9 let $(G,X)$ be a simple, adjoint Shimura pair of abelian type. Thus $G$ is a simple, adjoint group over $\dbQ$. Let $F:=F(G,X)$. We write $G=\Res_{F/\dbQ} G[F]$, cf. Subsection 2.2.
For a totally real number field $F_1$ that contains $F$, let $G^{F_1}:=\Res_{F_1/\dbQ} G[F]_{F_1}$ and $(G,X)\hookrightarrow (G^{F_1},X^{F_1})$ be as in Subsubsection 2.2.3.

\bigskip\noindent
{\bf 4.1. Proposition.} {\it We assume that the simple, adjoint Shimura pair $(G,X)$ is of $A_n$ type. Then there exists an injective map
$$f_2\colon (G_2,X_2)\hookrightarrow (\text{\bf GSp}(W,\psi),S)$$
of Shimura pairs such that the following two properties hold:

\medskip
{\bf (i)} the adjoint Shimura pair $(G_2^{\ad},X_2^{\ad})$ is isomorphic to $(G,X)$;

\smallskip
{\bf (ii)} the group $G_2$ is the subgroup of $\text{\bf GSp}(W,\psi)$ that fixes $\scrB:=\End(W)^{G_2(\dbQ)}$ (thus $f_2$ is a PEL type embedding).}

\medskip
\proof
Let $G[F]^{\sc}$ be the simply connected semisimple group cover of $G[F]$. The group $G^{\sc}:=\Res_{F/\dbQ} G[F]^{\sc}$ is the simply connected semisimple group cover of $G$. Let $T$ be a maximal torus of $G^{\sc}$.
Let $F_0$ be a finite field extension of $F$ such that the group $G[F]^{\sc}_{F_0}$ is split (i.e., it is an $\text{\bf SL}_{n+1,F_0}$ group). For instance, we can take $F_0$ to be the smallest Galois extension of $\dbQ$ such that the torus $T_{F_0}$ is split. Let $K$ be a totally imaginary quadratic extension of $F$ that is linearly disjoint from $F_0$; if $n=1$ and $G[F]$ is non-split (i.e., if $(G,X)$ is of compact $A_1$ type), then we choose $K$ such that $G[F]_K$ is also non-split. If $n>1$ (resp. $n=1$) let $E_0$ be $F_0$ (resp. be $F_0\otimes_F K$). Let $W_0$ be an $F_0$-vector space of dimension $n+1$. Let $W_2:=W_0\otimes_{F_0} E_{0}$; when we view $W_2$ as a $\dbQ$--vector space we denote it by $W$. We identify $G[F]^{\sc}_{F_0}=\text{\bf GL}_{W_0}^{\der}$ and thus also $G[F]^{\sc}_{E_0}=\text{\bf GL}_{W_2}^{\der}$. We have natural monomorphisms $G^{\sc}\hookrightarrow \Res_{E_0/\dbQ} G[F]^{\sc}_{E_0}\hookrightarrow \text{\bf GL}_{W}^{\der}$ and thus we can view $G^{\sc}$ as a subgroup of $\text{\bf GL}_W$. Next we follow the proof of [De2, Prop. 2.3.10].

Let $\scrS$ be the set of extremal vertices of the Dynkin diagram $\grD$ of $\Lie(G_{\overline{\dbQ}})$ with respect to $T_{\overline{\dbQ}}$ and a Borel subgroup of $G_{\overline{\dbQ}}$ that contains $T_{\overline{\dbQ}}$. The Galois group $\Gal(\dbQ)$ acts on $\scrS$ (see Subsection 2.2). If $n>1$ we can identify $\scrS$ with the $\Gal(\dbQ)$-set $\Hom_{\dbQ}(\scrK,\overline{\dbQ})$, where $\scrK:=I(G,X)$ is a totally imaginary quadratic extension of $F$. If $n=1$ let $\scrK:=K$. Always $\scrK$ is a subfield of $E_0$ and thus $W_2$ is naturally a $\scrK$-vector space. Therefore the torus $\scrT:=\Res_{\scrK/\dbQ} \dbG_{m,\scrK}$ acts naturally on $W$. If $n>1$ this action over $\overline{\dbQ}$ can be described as follows: if $V_0$ is a simple $G^{\sc}_{\overline{\dbQ}}$-submodule of $W\otimes_{\dbQ} \overline{\dbQ}$ whose highest weight is the fundamental weight associated to a vertex $\grn\in\scrS$, then $\scrT_{\overline{\dbQ}}$ acts on $V_0$ via its character that corresponds naturally to $\grn$. As $\scrK$ is a totally imaginary quadratic extension of $F$, for $n\Ge 1$ the extension to $\dbR$ of the quotient torus $\scrT/\Res_{F/\dbQ}\dbG_{m,F}$ is compact. This implies that the maximal subtorus $\scrT_c$ of $\scrT$ which over $\dbR$ is compact, is isomorphic to $\scrT/\Res_{F/\dbQ}\dbG_{m,F}$ and therefore it has rank $[F:\dbQ]$.

The centralizer $C$ of $G^{\sc}$ in $\text{\bf GL}_W$ is a reductive group and $Z(C)$ is a torus of rank equal to $2[F:\dbQ]$ if $n>1$ (resp. to $[F:\dbQ]$ if $n=1$). If $n>1$ (resp. $n=1$) we view $\scrT$ and $\scrT_c$ as subtori of $Z(C)$ (resp. of $C$). Let $\scrT_{c,+}$ be the torus of $C$ generated by $\scrT_c$ and $Z(\text{\bf GL}_W)$; its rank is $[F:\dbQ]+1$ and it commutes with $G^{\sc}$. Let $G_2$ be the reductive subgroup of $\text{\bf GL}_W$ that is generated by $G^{\sc}$ and $\scrT_{c,+}$; thus we have $G_2^{\der}=G^{\sc}$ and $Z^0(G_2)=\scrT_{c,+}$. The centralizer $C_2$ of $G_2$ in $\text{\bf GL}_W$ is also a reductive group. Let $\scrB$ be the semisimple $\dbQ$--subalgebra of $\End(W)$ defined by the elements of $\Lie(C_2)$. We have $\scrB=\End(W)^{G_2(\dbQ)}$.

We fix an element $x\in X$. We will construct a monomorphism $x_2:\dbS\hookrightarrow G_{2,\dbR}$ such that it defines a Hodge $\dbQ$--structure on $W$ of type $\{(-1,0),(0,-1)\}$ and it lifts the homomorphism $x:\dbS\to G^{\ad}_{2,\dbR}=G_{\dbR}$.

Let $\grR:=\Hom_{\dbQ}(F,\dbR)$. We have $G_{\dbR}^{\sc}=\prod_{i\in\grR} G[F]^{\sc}\times_{F,i} \dbR$. For $i\in \grR$ let $V_{i}$ be the real vector subspace of $W\otimes_{\dbQ} \dbR$ generated by its simple $G[F]^{\sc}\times_{F,i} \dbR$-submodules. We have a product decomposition $\scrT_{\dbR}=\prod_{i\in \grR} T_{i}$ in 2 dimensional tori such that for each pair $(i,i^\prime)\in\grR^2$, the torus $T_{i}$ acts trivially on $V_{i^\prime}$ if and only if $i\neq i^\prime$. Each $T_{i}$ is isomorphic to $\dbS$ and acts faithfully on $V_{i}$. We have
$$\scrT_{c,\dbR}=\prod_{i\in \grR} T_{c,i},$$
where $T_{c,i}$ is the $1$ dimensional compact subtorus of $T_{i}$.
The direct sum decomposition
$$W\otimes_{\dbQ} \dbR=\oplus_{i\in \grR} V_{i}$$
is normalized by $G_{2,\dbR}$. We have a unique direct sum decomposition
$$V_{i}\otimes_{\dbR} \dbC=V^+_{i}\oplus V^{-}_{i}\leqno (2)$$
such that $\scrT_{\dbC}$ acts on each $V^u_{i}$ via a unique character that is uniquely determined by the pair $(i,u)$ (here $u\in\{-,+\}$). Each $G^{\der}_{2,\dbC}$-module $V_{i}^u$ is isotypic (i.e., is a direct sum of isomorphic simple $G^{\der}_{2,\dbC}$-modules). More precisely, the highest weight of the representation of the factor $G[F]^{\sc}\times_{F,i} \dbC$ of $G^{\der}_{2,\dbC}$ on $V_{i}^u$ is $\varpi_{s_i(u)}$, where $s_i:\{+,-\}\to \{1,n\}$ is a surjective function. Thus $I_{i}:=\im(T_{i}\to \text{\bf GL}_{V_{i}})$ is the center of the centralizer of $C_{2,\dbR}$ in $\text{\bf GL}_{V_i}$.

We will take $x_2:\dbS\hookrightarrow G_{2,\dbR}$ such that its restriction to the split subtorus $\dbG_{m,\dbR}$ of $\dbS$ induces an isomorphism $\dbG_{m,\dbR}\arrowsim Z(\text{\bf GL}_{W\otimes_{\dbQ} \dbR})$ and the following two properties hold:

\medskip
\item{{\bf (a)}} if $i\in\grR$ is such that $G[F]\times_{F,i} \dbR$ is non-compact, then the homomorphism $\dbS\to \text{\bf GL}_{V_i}$ defined by $x_2$ is constructed as in the proof of [De2, Prop. 2.3.10] (thus this homomorphism factors through $\text{Im}(G_{2,\dbR}\to\text{\bf GL}_{V_i})$, defines a Hodge $\dbR$-structure on $V_i$ of type $\{(-1,0),(0,-1)\}$, and it is uniquely determined by these two properties and the fact that it lifts the homomorphism $\dbS\to G[F]\times_{F,i} \dbR$ defined by $x$);

\smallskip
\item{{\bf (b)}} if $i\in\grR$ is such that $G[F]\times_{F,i} \dbR$ is compact, then the homomorphism $\dbS\to \text{\bf GL}_{V_i}$ defined by $x_2$ is a monomorphism whose image is naturally identified with $I_{i}$ (and thus whose image in $G[F]\times_{F,i} \dbR$ is trivial).

\medskip
For the rest of the proof it is irrelevant which one of the two possible natural identifications of (b) we choose; each such choice is obtained naturally from the other one via the standard non-trivial automorphism of the compact subtorus of $\dbS$.

Let $X_2$ be the $G_2(\dbR)$-conjugacy class of $x_2$. The homomorphism $x_2$ lifts $x\in X$, cf. (a) and (b) and the fact (see Subsection 2.2) that the image of $x$ in each simple, compact factor of $G_{\dbR}$ is trivial. Therefore the pair $(G_2,X_2)$ is a Shimura pair whose adjoint is naturally isomorphic to $(G,X)$. Thus (i) holds.
The torus $\scrT$ splits over the field of CM type which is the Galois extension of $\dbQ$ generated by $\scrK$. Due to (a) and (b), the homomorphism $x_2$ defines a Hodge $\dbQ$--structure on $W$ of type $\{(-1,0),(0,-1)\}$. From the last two sentences and [De2, Cor. 2.3.3], we get that there exists a non-degenerate alternating form $\psi$ on $W$ such that we have an injective map $f_2:(G_2,X_2)\hookrightarrow (\text{\bf GSp}(W,\psi),S)$ of Shimura pairs.

The representation of $G_{2,\overline{\dbQ}}$ on $W\otimes_{\dbQ} \overline{\dbQ}$ is a direct sum of irreducible representations of dimension $n+1$ which (as $\overline{\dbQ}$--vector spaces) are isotropic with respect to $\psi$. For all $n\Ge 1$, the number of non-isomorphic such irreducible representations is $2[F:\dbQ]$. Thus the subgroup $G^\prime_{2,\overline{\dbQ}}$ of $\text{\bf GSp}(W\otimes_{\dbQ} \overline{\dbQ},\psi)$ that centralizes $\scrB\otimes_{\dbQ} \overline{\dbQ}$ (equivalently $C_{2,\overline{\dbQ}}$) is a reductive group that has the following properties: it contains $G_{2,\overline{\dbQ}}$, $G^{\prime,\der}_{2,\overline{\dbQ}}\arrowsim\text{\bf SL}_{n+1,\overline{\dbQ}}^{[F:\dbQ]}$, and $Z^0(G^\prime_{2,\overline{\dbQ}})\arrowsim\dbG_{m,\overline{\dbQ}}^{[F:\dbQ]+1}$. By reasons of dimensions we have $G_{2,\overline{\dbQ}}=G^\prime_{2,\overline{\dbQ}}$. Therefore the subgroup of $\text{\bf GSp}(W,\psi)$ that centralizes $\scrB$ (equivalently $C_2$) is $G_2$. Thus (ii) also holds.\endproof

\medskip\noindent
{\bf 4.1.1. Two properties.} {\bf (a)} The double centralizer $G_5$ of $G_2$ in $\text{\bf GL}_W$ is such that $G_{5,\overline{\dbQ}}$ is isomorphic to $\text{\bf GL}_{n+1,\overline{\dbQ}}^{2[F:\dbQ]}$ and $G_5^{\ad}$ is $\Res_{\scrK/\dbQ} G[F]_{\scrK}$. The semisimple $\dbQ$--algebra $\scrB$ has $\scrK$ as its center and thus it is simple. Moreover, if $n>1$ (resp. $n=1$) and $(G,X)$ is of strong compact type (resp. of compact type), then from the description of $\scrK$ in the proof of Proposition 4.1 we get that the $\scrK$-rank of $G[F]_{\scrK}$ is $0$; thus the $\dbQ$--rank of $G_5^{\ad}$ is also $0$.

\smallskip
{\bf (b)} We assume that there exists a number $a\in\dbN$ such that $n=2a-1$ and $G^{\sc}_{\dbR}=G_{2,\dbR}^{\der}$ is a product of $\text{\bf SU}(a,a)_{\dbR}$ groups. Then the monomorphism $x_2:\dbS\hookrightarrow G_{2,\dbR}$ factors through the subgroup of $G_{2,\dbR}$ generated by $G_{2,\dbR}^{\der}$ and $Z(\text{\bf GL}_{W\otimes_{\dbQ} \dbR})$, cf. [De2, Rms. 2.3.12 and 2.3.13]. In other words, the homomorphism $\dbS\to\text{\bf GL}_{V_i}$ of the property 4.1 (a) factors through the subgroup of $\text{\bf GL}_{V_i}$ generated by $\text{Im}(G^{\der}_{2,\dbR}\to\text{\bf GL}_{V_i})$ and $Z(\text{\bf GL}_{V_i})$.

\medskip\noindent
{\bf 4.1.2. A variant of 4.1.} We refer to the proof of Proposition 4.1 with $n>1$ and with $(G,X)$ without involution. We take $E_0$ to be $F_0\otimes_F K$, with $K$ a totally imaginary quadratic extension of $F$ that is linearly disjoint from $\scrK=I(G,X)$. Let $\scrT:=\Res_{\scrK\otimes_F K/\dbQ}\dbG_{m,\scrK\otimes_F K}$. As in the proof of Proposition 4.1, $W$ is naturally a $\scrK\otimes_F K$-module and thus a $\scrT$-module. Let $\scrT_c$ and $\scrT_{c,+}$ be the subtori of $\scrT$ obtained as in the proof of Proposition 4.1; their ranks are $2[F:\dbQ]$ and $2[F:\dbQ]+1$ (respectively). We take $G_2$ to be defined as in the proof of Proposition 4.1. We have a unique direct sum decomposition in $G_{2,\dbC}$-modules
$$V_i\otimes_{\dbR} \dbC=V^{+,+}_{i}\oplus V^{+,-}_{i}\oplus V^{-,+}_{i}\oplus V^{-,-}_{i}$$
such that all irreducible subrepresentations of the representation of $G_{2,\dbC}^{\der}$ (resp. of $Z^0(G_{2,\dbC})$) on $V^{u,+}_{i}\oplus V^{u,-}_{i}$ (resp. on $V^{+,u}_{i}\oplus V^{-,u}_{i}$) are isomorphic; here $u\in\{-,+\}$.

The proof of [De2, Prop. 2.3.10] guarantees the existence of a monomorphism $x_2:\dbS\hookrightarrow G_{2,\dbR}$ such that it defines a Hodge $\dbQ$--structure on $W$ of type $\{(-1,0),(0,-1)\}$ and it lifts the homomorphism $x:\dbS\to G^{\ad}_{2,\dbR}=G_{\dbR}$. As in the proof of Proposition 4.1 we argue that there exists a non-degenerate alternating form $\psi$ on $W$ such that we have an injective map $f_2:(G_2,X_2)\hookrightarrow (\text{\bf GSp}(W,\psi),S)$ of Shimura pairs. The representation of $G_{2,\overline{\dbQ}}$ on $W\otimes_{\dbQ} \overline{\dbQ}$ is a direct sum of irreducible representations which (as $\overline{\dbQ}$--vector spaces) are isotropic with respect to $\psi$. But this time the number of non-isomorphic such irreducible representations is $4[F:\dbQ]$. Thus the subgroup $G_4$ of $\text{\bf GSp}(W,\psi)$ that centralizes $\scrB=\End(W)^{G_2(\dbQ)}$ is not $G_2$ but a reductive group with the property that $G^{\der}_{4,\overline{\dbQ}}$ is isomorphic to $(G^{\der}_{2,\overline{\dbQ}})^2$.

As $G_2\leqslant G_4\leqslant \text{\bf GSp}(W,\psi)$ and as obviously there exists no simple factor of $G_4^{\ad}$ that over $\dbR$ is compact, the $G_4(\dbR)$-conjugacy class $X_4$ of $x_2$ has the property that we have injective maps $(G_2,X_2)\hookrightarrow (G_4,X_4)\hookrightarrow (\text{\bf GSp}(W,\psi),S)$ of Shimura pairs. The injective map $(G_4,X_4)\hookrightarrow (\text{\bf GSp}(W,\psi),S)$ of Shimura pairs is a PEL type embedding, cf. the very definition of $G_4$. As $F_0\otimes_F K$ is a field (cf. our linearly disjoint assumption), the center of the $\dbQ$--algebra $\scrB$ is a $\dbQ$--subalgebra of $F_0\otimes_F K$ and thus it is a field. Therefore the $\dbQ$--algebra $\scrB$ is simple. Thus the injective map $(G,X)=(G_2^{\ad},X_2^{\ad})\hookrightarrow (G_4^{\ad},X_4^{\ad})$ of Shimura pairs is as in (1) of Subsubsection 2.2.3, with $F_1$ as the unique totally real quadratic extension of $F$ contained in $F_0\otimes_F K$. We have identities $E(G_2^{\ad},X_2^{\ad})=E(G,X)=E(G^{F_1},X^{F_1})=E(G_4^{\ad},X_4^{\ad})$, cf. Subsubsection 2.2.3.

\medskip\noindent
{\bf 4.1.3. Remark.} If in Proposition 4.1 (resp. in Subsubsection 4.1.2) we have $n>1$ and we can take $G_2$ (resp. $G_4$) such that $Z^0(G_2)=\dbG_{m,\dbQ}$ (resp. such that $Z^0(G_4)=\dbG_{m,\dbQ}$), then the property 4.1 (ii) (resp. the part of Subsubsection 4.1.2 that involves PEL type embeddings) does not hold anymore.

\bigskip\noindent
{\bf 4.2. The group $\tilde G$.} Until Subsection 4.9 we will assume that $(G,X)$ is of $B_n$, $C_n$, $D_n^{\dbH}$, or $D_n^{\dbR}$ type. All simple, non-compact factors of $G_{\dbR}$ are isomorphic. We now introduce the ``simplest" absolutely simple, adjoint group over $\dbQ$ with the property that $\tilde G_{\dbR}$ is isomorphic to each simple, non-compact factor of $G_{\dbR}$. Depending on the type of $(G,X)$ we take $\tilde G$ as follows (the values of $n$ are such that we do not get repetitions among the Shimura types).

\medskip
$\bullet$ If $(G,X)$ is of $B_n$ type ($n>1$), then $\tilde G:=\text{\bf SO}(2,2n-1)$.

\smallskip
$\bullet$ If $(G,X)$ is of $C_n$ type ($n>2$), then $\tilde G:=\text{\bf Sp}^{\ad}_{2n,\dbQ}$.

\smallskip
$\bullet$ If $(G,X)$ is of $D_n^{\dbH}$ type ($n\Ge 4$), then $\tilde G:=\text{\bf SO}^*(2n)^{\ad}$.

\smallskip
$\bullet$ If $(G,X)$ is of $D_n^{\dbR}$ type ($n\Ge 4$), then $\tilde G:=\text{\bf SO}(2,2n-2)^{\ad}$.

\medskip
Let $\tilde G^{\sc}$ be the simply connected semisimple group cover of $\tilde G$. Let $\tilde x:\dbS\to\tilde G_{\dbR}$ be a homomorphism whose $\tilde G(\dbR)$-conjugacy class $\tilde X$ has the property that the pair $(\tilde G,\tilde X)$ is a Shimura pair. We know that $\tilde x$ is unique up to $\tilde G(\dbR)$-conjugation and possibly up to a replacement of it by its inverse, cf. [De2, Prop. 1.2.7 and Cor. 1.2.8].

\bigskip\noindent
{\bf 4.3. The group $G(A)$.} Until Subsection 4.9 we will have the same restrictions on $n$ as in Subsection 4.2. Let the group $G(A)$ over $\dbQ$ be defined by the rules:

\medskip
$\bullet$  if $(G,X)$ is of $B_n$ or $D_n^{\dbR}$ type, then $G(A):=\text{\bf SU}(2^{n-1},2^{n-1})$;

\smallskip
$\bullet$  if $(G,X)$ is of $C_n$ or $D_n^{\dbH}$ type, then $G(A):=\text{\bf SU}(n,n)$.

\medskip
These rules conform with [Sa1, Subsects. 3.3 to 3.7 or p. 461]. Subsections [Sa1, 1.1 and 3.3 to 3.7] are stated in terms of simple Lie algebras over $\dbR$ and show the existence of a Lie monomorphism $dh_{\dbR}:\Lie(\tilde G_{\dbR}^{\sc})\hookrightarrow\Lie(G(A)_{\dbR})$ that satisfies the condition $(H_2)$ of [Sa1, end of Subsect. 1.1]. This condition says that we have an identity
$$(dh_{\dbR}\circ d\tilde x)(\sqrt{-1})=d\tilde y(\sqrt{-1})\in\Lie(G(A)_{\dbR}^{\ad}),\leqno (3)$$
where $\sqrt{-1}\in\dbC=\Lie(\dbS)$, where $\tilde x$ is as in the end of Subsection 4.2, and where $\tilde y:\dbS\to G(A)^{\ad}_{\dbR}$ is a homomorphism that induces a monomorphism $\dbS/\dbG_{m,\dbR}\hookrightarrow G(A)^{\ad}_{\dbR}$ whose image is the identity component of the center of a maximal compact subgroup of $G(A)^{\ad}_{\dbR}$. Here $d\tilde x$ and $d\tilde y$ are the Lie homomorphisms induced naturally by $\tilde x$ and $\tilde y$ (respectively). Let
$$h_{\dbR}\colon\tilde G^{\sc}_{\dbR}\to G(A)_{\dbR}$$
be the homomorphism of finite kernel whose differential at the level of Lie algebras is $dh_{\dbR}$.

The only difference here from [Sa1, Subsects. 3.3 to 3.7 or p. 461] is the following one. See [He, Ch. X, Sect. 6, Table V] for the classification of hermitian symmetric domains. If $(G,X)$ is of $D_n^{\dbR}$ type, then $\tilde X$ is two copies of the hermitian symmetric domain $BD\, I(2,2n-2)$ associated to the group $\text{\bf SO}(2,2n-2)^{\ad}_{\dbR}$. In [Sa1, Subsect. 3.5] this hermitian symmetric domain is embedded into the hermitian symmetric domain $A\, III(2^{n-2},2^{n-2})$ associated to the group $\text{\bf SU}(2^{n-2},2^{n-2})_{\dbR}^{\ad}$. We have two such embeddings: they correspond to the two half spin representations of $\tilde G^{\sc}_{\dbC}$. We get our $h_{\dbR}$ by working with the spin representation of $\tilde G^{\sc}_{\dbC}$ i.e., by composing the homomorphism $\tilde h_{\dbR}:\tilde G_{\dbR}^{\sc}\to \text{\bf SU}(2^{n-2},2^{n-2})_{\dbR}\times_{\dbR} \text{\bf SU}(2^{n-2},2^{n-2})_{\dbR}$ defined by the mentioned two embeddings with a standard monomorphism
$$\Delta:\text{\bf SU}(2^{n-2},2^{n-2})_{\dbR}\times_{\dbR} \text{\bf SU}(2^{n-2},2^{n-2})_{\dbR}\hookrightarrow \text{\bf SU}(2^{n-1},2^{n-1})_{\dbR}.$$
It is easy to see that $\Delta$ is uniquely determined up to $\text{\bf SU}(2^{n-1},2^{n-1})_{\dbR}(\dbR)$-conjugation.

Let $V_A$ be a complex vector space such that $G(A)_{\dbC}$ is (isomorphic to) $\text{\bf GL}_{V_A}^{\der}$. If $(G,X)$ is of $C_n$ or $D_n^{\dbH}$ (resp. $B_n$ or $D_n^{\dbR}$) type, then the representation of $\tilde G_{\dbC}^{\sc}$ on $V_A$ is the standard (resp. is the spin) representation of dimension $2n$ (resp. $2^n$).

\medskip\noindent
{\bf 4.3.1. Simple properties.} {\bf (a)} We assume that $(G,X)$ is of $D_n^{\dbR}$ type. Thus  $\tilde G^{\sc}$ is $\text{\bf Spin}(2,2n-2)$. We check that $h_{\dbR}$ is the composite $h_{\dbR}^{(1)}:\text{\bf Spin}(2,2n-2)_{\dbR}\to G(A)_{\dbR}$ of a standard monomorphism $\text{\bf Spin}(2,2n-2)_{\dbR}\hookrightarrow \text{\bf Spin}(2,2n-1)_{\dbR}$ with the monomorphism $\text{\bf Spin}(2,2n-1)_{\dbR}\hookrightarrow G(A)_{\dbR}$ we mentioned above but for the $B_n$ type case. The representation of $\tilde G^{\sc}_{\dbC}$ on $V_A$ is the spin representation, regardless if this representation is obtained via $h_{\dbR}$ or via $h_{\dbR}^{(1)}$. Thus the double centralizer of $\text{Im}(h_{\dbR}^{(1)})$ in $G(A)_{\dbR}$ has a derived group $D$ whose extension to $\dbC$ is isomorphic to $\text{\bf SL}_{2^{n-1},\dbC}^2$. As $h_{\dbR}^{(1)}$ satisfies at the level of Lie algebras the condition $(H_2)$ and as $\text{Im}(h_{\dbR}^{(1)})\leqslant D\leqslant G(A)_{\dbR}$, $D$ has a compact torus of the same rank as $D$. Thus each simple factor of $D^{\ad}$ is absolutely simple and of the form $\text{\bf SU}(2^{n-1}-a,a)_{\dbR}^{\ad}$ for some $a\in\{0,\ldots,2^{n-1}\}$ (see [He, Ch. X, Sect. 6, Table V]). The homomorphism $\tilde G_{\dbR}^{\sc}\to D$ satisfies at the level of Lie algebras the condition $(H_2)$ (as $h_{\dbR}^{(1)}$ has this property). Thus from [Sa1, Subsect. 3.3] we get that $D$ is a product of two copies of $\text{\bf SU}(2^{n-2},2^{n-2})_{\dbR}$ (i.e., $a$ is $2^{n-2}$) and that the homomorphism $\tilde G_{\dbR}^{\sc}\to D$ can be identified with $\tilde h_{\dbR}$ of Subsection 4.3. Also $\Delta$ can be identified with the monomorphism $D\hookrightarrow G(A)_{\dbR}$. We conclude that $h_{\dbR}$ can be identified with $h_{\dbR}^{(1)}$.

\smallskip
{\bf (b)} The finite group $\text{Ker}(h_{\dbR})$ is trivial if and only if $(G,X)$ is of $B_n$, $C_n$, or $D_n^{\dbR}$ type (cf. the description of the $\tilde G^{\sc}_{\dbC}$-module $V_A$).

\smallskip
{\bf (c)} From now on we take $h_{\dbR}$ such that it is the extension to $\dbR$ of a homomorphism
$$h\colon\tilde G^{\sc}\to G(A)$$
over $\dbQ$ (thus the notation $h_{\dbR}$ is justified). The fact that such a choice of $h_{\dbR}$ is possible follows from the constructions in [Sa1, Subsects. 3.3 to 3.7]: one only has to take $\Delta$ to be defined over $\dbQ$ and to replace in loc. cit. $\dbR$, $\dbC$, and $\dbH$ by $\dbQ$, $\dbQ(i)$, and respectively by the standard quaternion algebra $\dbQ(i,j)$ over $\dbQ$. Let
$$\tilde G^d:=\tilde G^{\sc}/\Ker(h).$$

\noindent
{\bf 4.3.2. Fact.} {\it Let $K$ be a field of characteristic $0$. Let $\tilde G^d_K\hookrightarrow G(A)_K$ be the monomorphism defined by $h_K:\tilde G^{\sc}_K\to G(A)_K$. Then the centralizer of a maximal torus of $\tilde G^d_K$ in $G(A)_K$, is a maximal torus of $G(A)_K$.}

\medskip
\proof
We can assume that $K$ is $\dbC$. Each such centralizer is a reductive group, cf. [Bor, Ch. IV, 13.17, Cor. 2 (a)]. Thus to prove the Fact, it suffices to show that the centralizer of a maximal torus $\tilde T$ of $\tilde G^d_{\dbC}$ in $\text{\bf GL}_{V_A}$ is a torus. It suffices to check that all weight spaces of the representation of $\tilde T$ on $V_A$ have dimension $1$. The fact that this holds is a basic property of minuscule weights (cf. [Bou2, pp. 127--129]) but for readers' convenience we recall one simple way to argue it.

If $(\tilde G,\tilde X)$ is of $C_n$ (resp. $D_n^{\dbH}$) type, then the monomorphism $\tilde G^d_{\dbC}\hookrightarrow \text{\bf GL}_{V_A}$ is a standard $\text{\bf Sp}_{2n,\dbC}\hookrightarrow \text{\bf GL}_{2n,\dbC}$ (resp. $\text{\bf SO}_{2n,\dbC}\hookrightarrow \text{\bf GL}_{2n,\dbC}$) monomorphism and it is well known that all the mentioned weight spaces have dimension $1$. If $(\tilde G,\tilde X)$ is of $B_n$ or $D_n^{\dbR}$ type, then the representation of $\tilde G^d_{\dbC}$ on $V_A$ is a spin representation and thus (cf. [FH, Props. 20.15 and 20.20]) all the mentioned weight spaces have dimension $1$.\endproof

\bigskip\noindent
{\bf 4.4. Lemma.} {\it There exists a totally real finite field extension $F_1$ of $F$ of degree a divisor of $2$ and such that the group $G[F]_{F_1}$ is an inner form of $\tilde G_{F_1}$.}

\medskip
\proof
The group $G[F]$ is a form of $\tilde G_F$. If $(G,X)$ is of $B_n$ or $C_n$ type, then this form is inner and thus we can take $F_1=F$. Let now $(G,X)$ be of $D_n^{\dbH}$ or $D_n^{\dbR}$ type ($n\Ge 4$). Let $Aut(\tilde G_F)$ be the group scheme of automorphisms of $\tilde G_F$. The group $\tilde G_F$ is a normal subgroup of $Aut(\tilde G_F)$ and the quotient group $Q:=Aut(\tilde G_F)/\tilde G_F$ is $\dbZ/2\dbZ$ if $n\Ge 5$ and is the symmetric group $S_3$ if $n=4$. We have a short exact sequence $1\to Q_1\to Q\to\dbZ/2\dbZ\to 1$, where $Q_1$ is $\dbZ/3\dbZ$ if $n=4$ and is the trivial group if $n>5$. Let $Aut_1(\tilde G_F)$ be the normal subgroup of $Aut(\tilde G_F)$ which contains $\tilde G_F$ and whose group of connected components is $Q_1$.

Let $\delta_0\in H^1(F,Aut(\tilde G_F))$ be the class that defines the form $G[F]$ of $\tilde G_{F}$. In the \'etale topology of $\Spec(F)$, the non-trivial torsors of $\dbZ/2\dbZ$ correspond to quadratic extensions of $F$. For each embedding $F\hookrightarrow\dbR$, the group $G[F]_{\dbR}$ is either $\tilde G_{\dbR}$ or the compact form of $\tilde G_{\dbR}$; thus $G[F]_{\dbR}$ is an inner form of $\tilde G_{\dbR}$ (cf. [De2, Subsubsect. 2.3.4]). This implies that there exists a smallest totally real finite field extension $F_1$ of $F$ of degree a divisor of 2 and such that the image of $\delta_0$ in $H^1(F_1,Aut(\tilde G_F)_{F_1})$ is the image of a class $\delta_1\in H^1(F_1,Aut_1(\tilde G_F)_{F_1})$. The Shimura pairs $(G,X)$ and $(G^{F_1},X^{F_1})$ have the same type that is either $D_n^{\dbH}$ or $D_n^{\dbR}$, cf. Subsubsection 2.2.3. If $n=4$, then from the very definition of the $D_4^{\dbR}$ and $D_4^{\dbH}$ types we get that the image of $\delta_1$ in $H^1(F_1,Q_1)$ is trivial. We conclude that for all $n\Ge 4$, $\delta_1$ is the image of a class in $H^1(F_1,\tilde G_{F_1})$. Thus $G[F]_{F_1}$ is an inner form of $\tilde G_{F_1}$ for all $n\Ge 4$. \endproof

\bigskip\noindent
{\bf 4.5. The twisting process.} Let $\delta\in H^1(F_1,\tilde G_{F_1})$ be the class that defines the inner form $G[F]_{F_1}$ of $\tilde G_{F_1}$. Let $\tilde G^e$ be the image of $\tilde G^d$ in $G(A)^{\ad}$. If $(G,X)$ is not of $D_n^{\dbR}$ type, then the representation of $\tilde G^d_{\dbC}$ on $V_A$ is irreducible and this implies that $Z(\tilde G^d)$ is a subgroup of $Z(G(A))$; thus $\tilde G^e$ is an adjoint group.

In this paragraph we assume that $(G,X)$ is of $D_n^{\dbR}$ type. The group $\tilde G^e_{\dbR}$ is the image of $\text{\bf Spin}(2,2n-2)_{\dbR}$ in $\text{\bf SO}(2,2n-1)_{\dbR}$ (cf. property 4.3.1 (a)) and thus it is $\text{\bf SO}(2,2n-2)_{\dbR}$. This implies that $\tilde G_e$ is $\text{\bf SO}(2,2n-2)$ and thus $\pmb{\mu}_{2,\dbQ}\arrowsim\ker(\tilde G^d\to \tilde G^e)$. The compact inner form $\text{\bf SO}(2n)_{\dbR}$ of $\text{\bf SO}(2,2n-2)_{\dbR}$ is obtained by twisting through a class in the set $H^1(\dbR,\text{\bf SO}(2,2n-2)_{\dbR})$. Thus the obstruction of lifting $\delta$ to a class $\delta_e\in H^1(F_1,\tilde G^e_{F_1})$, is measured by a class $\delta_2\in H^2(F_1,\pmb{\mu}_{2,F_1})$ whose image in $H^2(\dbR,\pmb{\mu}_{2,\dbR})$ under each embedding $F_1\hookrightarrow\dbR$ is trivial. Let $\scrD_2$ be the division algebra of center $F_1$ that corresponds to $\delta_2$. Based on [Ha, Lem. 5.5.3] applied to the reductive group $\text{\bf GL}_1(\scrD_2)$, we get the existence of a totally real finite field extension $F_2$ of $F_1$ which has a degree that divides $2$ and which is a maximal commutative $F_1$-subalgebra of $\scrD_2$. By replacing $F_1$ with $F_2$ we can assume that $\delta_2$ is the trivial class.

Thus, regardless of the type of $(G,X)$, by requesting only that the degree $[F_1:F]$ divides $4$ we can assume that $\delta$ is the image of a class $\delta_e\in H^1(F_1,\tilde G^e_{F_1})$. Let $h_e:\tilde G^e\hookrightarrow G(A)^{\ad}$ be the monomorphism defined by $h$. We use $\delta_e$ to twist $h_{e,F_1}$. By performing such a twist and by lifting the resulting monomorphism via isogenies, we get a monomorphism
$$h^t_{F_1}\colon G[F]^d_{F_1}\hookrightarrow G(A)^t_{F_1}.$$
Here $G[F]^d_{F_1}$ is the isogeny cover of $G[F]_{F_1}$ that corresponds to the isogeny cover $\tilde G^d$ of $\tilde G$ and $G(A)^t_{F_1}$ is the resulting inner form of $G(A)_{F_1}$ (the upper index $t$ stands for twist).

\bigskip\noindent
{\bf 4.6. Lemma.} {\it If by extension of scalars under an embedding $F_1\hookrightarrow\dbR$, the group $G[F]_{F_1}$ becomes a compact group (resp. it becomes $\tilde G_{\dbR}$), then under the same extension of scalars $G(A)^t_{F_1}$ becomes a compact group (resp. it becomes $G(A)_{\dbR}$).}

\medskip
\proof
We have to prove that if we twist $h_{e,\dbR}:\tilde G^e_{\dbR}\hookrightarrow G(A)^{\ad}_{\dbR}$ through the compact inner form of $\tilde G^{e}_{\dbR}$, then the inner form $\scrI^{\ad}$ of $G(A)^{\ad}_{\dbR}$ we get is the compact form. This can be proved using Cartan decompositions (see [He]). Here we include a proof of this that relies on Lemma 1.1.7. Based on the property 4.3.1 (a), the case of $D_n^{\dbR}$ type gets reduced to the $B_n$ type case. Thus it suffices to consider the cases when $(G,X)$ is of $B_n$, $C_n$, or $D_n^{\dbH}$ type. Let $\scrI$ be the simply connected semisimple group cover of $\scrI^{\ad}$.

Let $\scrC$ be a maximal compact subgroup of $\scrI$. Let $V_A$ be as before Subsubsection 4.3.1. We identify $\scrI_{\dbC}=\text{\bf GL}_{V_A}^{\der}$. If $T_{\dbR}$ is a maximal compact torus of $\tilde G_{\dbR}^e$, then the map $H^1(\dbR,T_{\dbR})\to H^1(\dbR,\tilde G^e_{\dbR})$ is surjective (see [Ko1, Lem. 10.2]). Thus, as $T_{\dbR}$ is contained in a maximal compact torus $T(A)$ of $G(A)^{\ad}_{\dbR}$, we get that $T(A)$ is isomorphic to a torus of $\scrI$. As the maximal compact Lie subgroups of $\scrI(\dbR)$ are $\scrI(\dbR)$-conjugate (see [He, Ch. VI, Sect. 2, Thm. 2.2 (ii)]), $\scrC$ has a subgroup isomorphic to $T(A)$. Thus $\scrC$ has the same rank as $\scrI$.

The representation of $\tilde G^d_{\dbC}$ on $V_A$ is irreducible. Thus the $\tilde G^d_{\dbC}$-module $V_A$ is simple. But up to $\scrI(\dbC)$-conjugation, we can assume that  $\tilde G^d_{\dbC}$ is a subgroup of $\scrC_{\dbC}$. Thus the $\scrC_{\dbC}$-module $V_A$ is also simple. From Lemma 1.1.7 applied to the inclusions $\Lie(\scrC_{\dbC})\hookrightarrow \Lie(\scrI_{\dbC})\hookrightarrow \End(V_A)$, we get that $\Lie(\scrI_{\dbC})=\Lie(\scrC_{\dbC})$. Thus $\scrI=\scrC$ i.e., the group $\scrI$ is compact. \endproof

\bigskip\noindent
{\bf 4.7. The Shimura pair $(G_1,X_1)$.} Let $G_1$ be the adjoint group of $G_1^{\sc}:=\Res_{F_1/\dbQ} G(A)_{F_1}^t$. Each simple factor of $G_{1,\dbR}$ is either compact or an $\text{\bf SU}(a,a)^{\ad}_{\dbR}$ group for some $a\in\dbN$, cf. Lemma 4.6. This implies the existence of an adjoint Shimura pair of the form $(G_1,X_1)$; it is uniquely determined (cf. [De2, Prop. 1.2.7 and Cor. 1.2.8 (ii)]).

\bigskip\noindent
{\bf 4.8. Theorem.} {\it We recall that $(G,X)$ is a simple, adjoint Shimura pair of $B_n$, $C_n$, $D_n^{\dbH}$, or $D_n^{\dbR}$ type and that $F_1$ is a totally real finite field extension of $F$ of degree a divisor of $4$. There exist injective maps
$$f_2:(G_2,X_2)\hookrightarrow (G_3,X_3)\hookrightarrow (G_4,X_4)\operatornamewithlimits{\hookrightarrow}\limits^{f_4} (\text{\bf GSp}(W_1,\psi_1),S_1)$$
of Shimura pairs such that the following seven properties hold:

\medskip
{\bf (i)} we have $(G_2^{\ad},X_2^{\ad})=(G,X)$, $(G_3^{\ad},X_3^{\ad})=(G^{F_1},X^{F_1})$, $(G_4^{\ad},X_4^{\ad})=(G_1,X_1)$, and $G_4^{\der}=G_1^{\sc}$;

\smallskip
{\bf (ii)} we have a natural monomorphism $G_2^{\ad}\hookrightarrow G_3^{\ad}$ which under the identifications of (i) is nothing else but the natural monomorphism $G\hookrightarrow G^{F_1}$;

\smallskip
{\bf (iii)} we have $G_3^{\der}=\Res_{F_1/\dbQ} G[F]^d_{F_1}$ and under this identification and of the ones of (i), the monomorphism $G_3^{\der}\hookrightarrow G_4^{\der}$ becomes the $\Res_{F_1/\dbQ}$ of $h^t_{F_1}$ (of Subsection 4.5);

\smallskip
{\bf (iv)} if $(G,X)$ is not (resp. is) of $D_n^{\dbH}$ type, then we have $E(G_2,X_2)=E(G_3,X_3)=E(G_4,X_4)$ (resp. we have $E(G_2,X_2)=E(G_3,X_3)$);

\smallskip
{\bf (v)} both $Z^0(G_2)$ and $Z^0(G_3)$ are subtori of $Z^0(G_4)$;

\smallskip
{\bf (vi)} if $(G,X)$ is not of $D_n^{\dbH}$ type, then both $G_2^{\der}$ and $G_3^{\der}$ are simply connected;

\smallskip
{\bf (vii)} the group $G_4$ is the subgroup of $\text{\bf GSp}(W_1,\psi_1)$ that fixes all elements of the semisimple $\dbQ$--algebra $\scrB:=\End(W_1)^{G_4(\dbQ)}$ (thus $f_4$ is a PEL type embedding).}

\medskip
\proof
Let $f_4:(G_4,X_4)\hookrightarrow (\text{\bf GSp}(W_1,\psi_1),S_1)$ be obtained as in Proposition 4.1 but working with $(G_1,X_1)$ (instead of $(G,X)$). Thus (vii) holds (cf. property 4.1 (ii)) and we have $G_4^{\der}=G_1^{\sc}$ and $(G_4^{\ad},X_4^{\ad})=(G_1,X_1)$. We take $G_3$ to be generated by $Z^0(G_4)$ and by the semisimple subgroup $\Res_{F_1/\dbQ} G[F]_{F_1}^d$ of $G_4^{\der}$. This takes care of (iii) and we have $Z^0(G_3)=Z^0(G_4)$. We take $G_2$ to be the reductive subgroup of $G_3$ such that $Z^0(G_2)=Z^0(G_4)$ and (ii) holds. Thus (v) also holds.

The monomorphism $G_{3,\dbR}^{\der}\hookrightarrow G_{4,\dbR}^{\der}$ is a product of $[F_1:\dbQ]$ monomorphisms that can be identified either with the monomorphism $\tilde G^d_{\dbR}\hookrightarrow G(A)_{\dbR}$ defined naturally by $h_{\dbR}$ or with a monomorphism between compact groups, cf. Lemma 4.6. Thus from (3) and the uniqueness part of Subsection 4.7, we get the existence of elements $x\in X$ and $x_1\in X_1$ such that the image of $dx(\sqrt{-1})\in\Lie(G_{\dbR})$ in $\Lie(G_{1,\dbR})$ is $dx_1(\sqrt{-1})\in\Lie(G_{1,\dbR})$. But $G_4^{\ad}(\dbQ)$ is dense in $G_4^{\ad}(\dbR)$ and therefore it permutes transitively the connected components of $X_1=X_4^{\ad}$. Thus by composing the monomorphism $G_4\hookrightarrow \text{\bf GSp}(W_1,\psi_1)$ with an isomorphism $G_4\arrowsim G_4$ defined by an element of $G_4^{\ad}(\dbQ)$, we can choose $x$ and $x_1$ such that $x_1$ is the image in $X_1=X_4^{\ad}$ of a point $x_4\in X_4$. Therefore the monomorphism $x_4:\dbS\hookrightarrow G_{4,\dbR}$ factors through $G_{2,\dbR}$ in such a way that the resulting homomorphism $\dbS\to G_{2,\dbR}^{\ad}=G_{\dbR}$ is an element of $X$. Thus by taking $X_2$ (resp. $X_3$) to be the $G_2(\dbR)$-conjugacy (resp. the $G_3(\dbR)$-conjugacy) class of the factorization of $x_4$ through $G_{2,\dbR}$ (resp. through $G_{3,\dbR}$), we get that (i) holds.

To check (iv), let $*\in\{2,3,4\}$. The reflex field of $E(G_*,X_*)$ is the composite field of $E(G_*^{\ad},X_*^{\ad})$ and $E(G_*^{\ab},X_*^{\ab})$, cf. Subsection 2.2. But $(G_*^{\ab},X_*^{\ab})$ and thus also $E(G_*^{\ab},X_*^{\ab})$ does not depend on $*$. Moreover, we have $E(G_2^{\ad},X_2^{\ad})=E(G_3^{\ad},X_3^{\ad})$ (cf. Subsubsection 2.2.3). We now assume that $(G,X)$ is not of $D_n^{\dbH}$ type. Let $(\grD_*,\grn_{X_*})$ be the analogue of $(\grD,\grn_X)$ of Subsection 2.2 but for the Shimura pair $E(G_*^{\ad},X_*^{\ad})$. We recall from [De2, Rm. 2.3.12] that $E(G_3^{\ad},X_3^{\ad})$ (resp. $E(G_4^{\ad},X_4^{\ad})$) is the fixed field of the subgroup of $\Gal(\dbQ)$ that stabilizes the union of the connected components of $\grD_3$ (resp. of $\grD_4$) that correspond to compact factors of $G_{3,\dbR}^{\ad}$ (resp. of $G_{4,\dbR}^{\ad}$). From this and Lemma 4.6 we get $E(G_3^{\ad},X_3^{\ad})=E(G_4^{\ad},X_4^{\ad})$. Thus (iv) holds.

Condition (vi) follows from the construction of $h^t_{F_1}$ and the property 4.3.1 (b). \endproof

\medskip\noindent
{\bf 4.8.1. Remark.} We assume that $n\in 2+2\dbN$ and that $(G,X)$ is of inner $D_n^{\dbR}$ type. Then we have a variant of Subsections 4.3 to 4.8 in which we work with only one fixed half spin representation of $\tilde G^{\sc}_{\dbC}$. Thus the group $G(A)$ becomes $\text{\bf SU}(2^{n-2},2^{n-2})$, we have $[F_1:F]\le 2$ (i.e., in Lemma 4.4 we have $F_1=F$ and in Subsection 4.5 we do not have to replace $F_1$ by $F_2$), and the groups $G_{2,\overline{\dbQ}}^{\der}$ and $G_{3,\overline{\dbQ}}^{\der}$ are products of $\text{\bf Spin}_{2n,\overline{\dbQ}}/\mu_{2,\overline{\dbQ}}$ groups which for $n\neq 4$ are not $\text{\bf SO}_{2n,\overline{\dbQ}}$ groups.

Also if $(G,X)$ is of non-inner $D_4^{\dbR}$ type, then we have a variant of Subsections 4.3 to 4.8 in which we construct $h:\tilde G^{\sc}\to G(A)$ as for the $D_4^{\dbH}$ type. Thus $G(A)$ is $\text{\bf SU}(4,4)$, $F_1\neq F$, and due to the non-inner assumption the groups $G_{2,\overline{\dbQ}}^{\der}$ and $G_{3,\overline{\dbQ}}^{\der}$ are products of $\text{\bf Spin}_{8,\overline{\dbQ}}$ groups.

\medskip\noindent
{\bf 4.8.2. Remarks.} {\bf (a)} In the proof of` Proposition 4.1, the degree $[F_0:\dbQ]$ is bounded in terms of $\dim(T_{\dbQ})$ and thus also of the rank of $G_{\dbC}$. Moreover, in Theorem 4.8 we have $[F_1:F]\Le 4$. Thus there exists an effectively computable number $M(G,X)\in\dbN$ such that in Proposition 4.1 (resp. Theorem 4.8), we can take the injective map $f_2:(G_2,X_2)\hookrightarrow (\text{\bf GSp}(W,\psi),S)$ (resp. $f_2:(G_2,X_2)\hookrightarrow (\text{\bf GSp}(W_1,\psi_1),S_1)$) of Shimura pairs in such a way that the following inequality $\dim_{\dbQ}(W_1)\Le M(G,X)$ holds.

\smallskip
{\bf (b)} We assume that $(G,X)$ is of $D_n^{\dbH}$ type. We denote also by $h:\tilde G^d\hookrightarrow G(A)$ the monomorphism defined by $h$. Then there exist monomorphisms over $\dbR$
$$\tilde G_{\dbR}^{d}\operatornamewithlimits{\hookrightarrow}^{h_{\dbR}} G(A)_{\dbR}\operatornamewithlimits{\hookrightarrow}^{h_{1,\dbR}} \text{\bf Sp}_{4n,\dbR}\operatornamewithlimits{\hookrightarrow}^{h_{2,\dbR}} \text{\bf SU}(2n,2n)_{\dbR}\leqno (4)$$
which at the level of Lie algebras satisfy the condition $(H_2)$. Argument: see [Sa1, Subsect. 1.4] for $h_{1,\dbR}$; as $h_{2,\dbR}$ we can take the analogue of $h_{\dbR}$ but for the $C_{2n}$ type.

By reasons of ranks of maximal compact subgroups, $\text{\bf SU}(2n,0)_{\dbR}$ is not a subgroup of a non-compact form of $\text{\bf Sp}_{4n,\dbR}$. Thus Lemma 4.6 also holds in the context of $h_{1,\dbR}$. As in the property 4.3.1 (c) we argue that (4) is the extension to $\dbR$ of monomorphisms $h$, $h_1$, and $h_2$ over $\dbQ$. Therefore Subsections 4.3 to 4.8 can be adapted in the context of the monomorphisms $h$, $h_1$, and $h_2$ to get that we can choose $f_2$ and $f_4$ such that the injective map $f_4:(G_4,X_4)\hookrightarrow (\text{\bf GSp}(W_1,\psi_1),S_1)$ of Shimura pairs factors as injective maps $(G_4,X_4)\hookrightarrow (G_5,X_5)\hookrightarrow (G_6,X_6)\hookrightarrow (\text{\bf GSp}(W_1,\psi_1),S_1)$ of Shimura pairs and moreover, besides properties 4.8 (i) to (v), we also have the following four extra properties:

\medskip
{\bf (i)} the group $G_6$ is the subgroup of $\text{\bf GSp}(W_1,\psi_1)$ that centralizes the semisimple $\dbQ$--algebra $\scrB_6:=\End(W_1)^{G_6(\dbQ)}$;

\smallskip
{\bf (ii)} the monomorphisms $G_3^{\der}\hookrightarrow G_4^{\der}$, $G_4^{\der}\hookrightarrow G_5^{\der}$, and $G_5^{\der}\hookrightarrow G_6^{\der}$ are $\Res_{F_1/\dbQ}$ of inner twists of the monomorphisms $h_{F_1}$, $h_{1,F_1}$, and $h_{2,F_1}$ (respectively);

\smallskip
{\bf (iii)} the torus $Z^0(G_*)$ does not depend on $*\in\{2,\ldots,6\}$;

\smallskip
{\bf (iv)} the group $G_3$ is the identity component of the subgroup of $G_5$ that centralizes the semisimple $\dbQ$--algebra $\scrB_3:=\End(W_1)^{G_3(\dbQ)}$.
\medskip

Property (iv) is an easy consequence of (iii). As $G_{5,\overline{\dbQ}}^{\der}$ is a product of $\text{\bf Sp}_{4n,\overline{\dbQ}}$ groups, the relative PEL property expressed by (iv) is a fundamental tool in showing that many properties enjoyed by $(G_5,X_5)$ are as well enjoyed by $(G_3,X_3)$ and thus implicitly by $(G_2,X_2)$ (we have in mind mainly applications to the Langlands--Rapoport conjecture mentioned in Subsubsect. 1.2.1). In other words, this relative PEL type property ``compensates" the fact that $G_2^{\der}$ and $G_3^{\der}$ are not simply connected.

\smallskip
{\bf (c)} We assume that $(G,X)$ is of $D_n^{\dbR}$ type. The monomorphism $G^{\der}_{3,\overline{\dbQ}}\hookrightarrow G^{\der}_{4,\overline{\dbQ}}$ is a product of $[F_1:\dbQ]$ copies of the spin (representation) monomorphism $\text{\bf Spin}_{2n,\overline{\dbQ}}\hookrightarrow \text{\bf SL}_{2^n,\overline{\dbQ}}$. Thus the identity component $Z_3$ of the centralizer of $G_3^{\der}$ in $G_4^{\der}$ is a torus of rank $[F_1:\dbQ]$. Let $G_{3.5}$ be the reductive subgroup of $G_4$ generated by $G_3$ and $Z_3$. Let $G_{3.75}$ be the reductive subgroup of $\text{\bf GSp}(W_1,\psi_1)$ that centralizes $\End(W_1)^{G_3(\dbQ)}$. We have monomorphisms $G_{3.5}\hookrightarrow G_{3.75}\hookrightarrow G_4$ which give birth to an identity $Z^0(G_{3.5})=Z^0(G_{3.75})$. Moreover, the group $G^{\der}_{3.75,\overline{\dbQ}}$ is isomorphic to $\text{\bf SL}_{2^{n-1},\overline{\dbQ}}^{2[F_1:\dbQ]}$.

\medskip\noindent
{\bf 4.8.3. Special cases.} We refer to Theorem 4.8. All simple factors of $G_{1,\dbR}$ are either compact or are $\text{\bf SU}(a,a)_{\dbR}$ groups, where $a\in\dbN$. Thus $(G_1,X_1)$ is without involution. For the rest of this Subsubsection we will assume that all simple factors of $G_{\dbR}$ are non-compact. From Lemma 4.6 we get that all simple factors of $G_{1,\dbR}$ and $G_{2,\dbR}^{\ad}$ are also non-compact. Thus the homomorphism $x_4:\dbS\to G_{4,\dbR}$ of the proof of Theorem 4.8 factors through the extension to $\dbR$ of the subgroup of $G_4$ generated by $G^{\der}_4$ and by $Z(\text{\bf GL}_{W_1})$, cf. property 4.1.1 (b). Thus we can choose $f_2$ such that we have $Z^0(G_2)=Z^0(G_4)=Z(\text{\bf GL}_{W_1})\arrowsim\dbG_{m,\dbQ}$.

In this paragraph we take $(G,X)$ to be of either $C_n$ or $D_n^{\dbH}$ type. The double centralizer $D_2$ of $G_2$ in $\text{\bf GL}_{W_1}$ is a reductive group whose extension to $\overline{\dbQ}$ is a product of $[F:\dbQ]$ copies of $\text{\bf GL}_{2n,\overline{\dbF}}$. The involution of $\End(W_1\otimes_{\dbQ} \overline{\dbQ})$ induces a non-trivial involution of the Lie algebra of each such copy of $\text{\bf GL}_{2n,\overline{\dbF}}$. This implies that $\Lie(D_2\cap \text{\bf GSp}(W_1,\psi_1))$ is $\Lie(G_2)$. Thus $G_2$ is the identity component of $D_2$ and therefore the injective map $f_2:(G_2,X_2)\hookrightarrow (\text{\bf GSp}(W_1,\psi_1),S_1)$ of Shimura pairs is a PEL type embedding.

In this paragraph we take $(G,X)$ to be of inner $D_4^{\dbR}$ type. If we work with the variant of Theorem 4.8 mentioned in Remark 4.8.1, then similar arguments show that we can assume that $Z^0(G_2)=Z^0(G_4)=Z(\text{\bf GL}_{W_1})\arrowsim\dbG_{m,\dbQ}$ and that $f_2:(G_2,X_2)\hookrightarrow (\text{\bf GSp}(W_1,\psi_1),S_1)$ is a PEL type embedding.

\bigskip\noindent
{\bf 4.9. Variants.} Using Fact 2.4.1 (b) we get that Proposition 4.1 and Theorem 4.8 have variants for the case when $(G,X)$ is an arbitrary Shimura pair of adjoint, abelian type. To state such variants, let $(G,X)=\prod_{i\in I} (G_i,X_i)$ be the product decomposition into simple, adjoint Shimura pairs indexed by a set $I$ which is disjoint from $\dbZ$. For $i\in I$ let $F_i:=F(G_i,X_i)$ and $\Gtil_i:=G_i[F_i]$, cf. Subsection 2.2.

{\bf (a)} There exists a product $(G_1,X_1)=\prod_{i\in I} (G_{i,1},X_{i,1})$ of simple, adjoint, unitary Shimura pairs and there exist injective maps
$$f_2:(G_2,X_2)\hookrightarrow (G_3,X_3)\hookrightarrow (G_4,X_4)\operatornamewithlimits{\hookrightarrow}\limits^{f_4} (\text{\bf GSp}(W_1,\psi_1),S_1)$$
of Shimura pairs such that the following seven properties hold:

\medskip
{\bf (i)} we have $(G_2^{\ad},X_2^{\ad})=(G,X)$ and $(G_3^{\ad},X_3^{\ad})=\prod_{i\in I} (G_i^{F_{i,1}},X_i^{F_{i,1}})$, with $F_{i,1}$ as a totally real number field that contains $F_i$;

\smallskip
{\bf (ii)} we get a monomorphism $G_2^{\ad}\hookrightarrow G_3^{\ad}$ that is a product indexed by $i\in I$ of the natural monomorphisms $G_i\hookrightarrow G_i^{F_{i,1}}$;

\smallskip
{\bf (iii)} the monomorphism $G_3^{\der}\hookrightarrow G_4^{\der}$ is a product indexed by $i\in I$ of monomorphisms that are either isomorphisms or are Weil restrictions of monomorphisms obtained as $h^t_{F_1}$ in Subsection 4.5, depending on the fact that $(G_i,X_i)$ is or is not of $A_{n_i}$ type for some $n_i\in\dbN$;

\smallskip
{\bf (iv)} we have $(G_4^{\ad},X_4^{\ad})=(G_1,X_1)$, the group $G_4^{\der}$ is a simply connected semisimple group, and $f_4$ is a PEL type embedding that is a Hodge quasi product indexed naturally by elements of $I$ of PEL type embeddings (cf. Fact 2.4.1 (b));

\smallskip
{\bf (v)} if for all $i\in I$ (resp. if there exists $i\in I$ such that) the Shimura pair $(G_i,X_i)$ is not (resp. is) of $D_n^{\dbH}$ type, then $E(G_2,X_2)=E(G_3,X_3)=E(G_4,X_4)$ (resp. $E(G_2,X_2)=E(G_3,X_3)$); also if for all $i\in I$ the pair $(G_i,X_i)$ is of some $A_{n_i}$ type, then $G_2=G_3=G_4$;

\smallskip
{\bf (vi)} both $Z^0(G_2)$ and $Z^0(G_3)$ are subtori of $Z^0(G_4)$;

\smallskip
{\bf (vii)} if none of the Shimura pairs $(G_i,X_i)$ (with $i\in I$) is of some $D_n^{\dbH}$ type, then $G_2^{\der}$ and $G_3^{\der}$ are simply connected semisimple groups.
\medskip

{\bf (b)} In Theorem 4.8 we can replace $F_1$ by any other totally real number field $F_2$ that contains it. If there is $i\in I$ such that the simple, adjoint Shimura pair $(G_i,X_i)$ is of $A_n$ type, then for each totally real number field extension $F_{i,2}$ of $F_{i,1}$, the simple, adjoint Shimura pair $(G_i^{F_{i,2}},X_i^{F_{i,2}})$ is also of $A_n$ type and therefore one can apply Proposition 4.1 to it. Thus, as in the proof of Theorem 4.8 we argue that for such an element $i\in I$ there exist injective maps $f_{i,2}:(G_{i,2},X_{i,2})=(G_{i,3},X_{i,3})\hookrightarrow (G_{i,4},X_{i,4})\operatornamewithlimits{\hookrightarrow}\limits^{f_{i,4}} (\text{\bf GSp}(W_{i,1},\psi_{i,1}),S_{i,1})$ of Shimura pairs such that the analogues of properties 4.8 (i) to (vii) hold for them (therefore $(G_{i,2}^{\ad},X_{i,2}^{\ad})=(G_i,X_i)$, $(G_{i,4},X_{i,4})=(G_i^{F_{i,2}},X_i^{F_{i,2}})$, $f_{i,4}$ is a PEL type embedding, etc.)

Based on the previous paragraph, by using composites of fields and injective maps as in (1) we get an uniformized variant of (a) in which to the property (iii) of (a) we can add the requirement that the field $F_{i,1}$ does not depend on $i\in I$ (thus in general we have no upper bounds for the degrees $[F_{i,1}:F_i]$).

\smallskip
{\bf (c)} We have variants of (a) and (b) in which for the $A_n$ type without involution and with $n>1$ we use Subsubsection 4.1.2 instead of Proposition 4.1. In what follows we will only use the variant of (a) that is modeled on Subsubsection 4.1.2. For this variant we only have to perform the following two changes to (a):

\medskip
-- in the property (iii) of (a) we have to replace ``either isomorphisms" by ``either monomorphisms which over $\overline{\dbQ}$ are products of diagonal embeddings of the form $\text{\bf SL}_{n_i+1,\overline{\dbQ}}\hookrightarrow \text{\bf SL}_{n_i+1,\overline{\dbQ}}^{s_i}$, where $s_i\in\{1,2\}$,";

-- in  the property (v) of (a) we have to replace ``some $A_{n_i}$ type" by ``some $A_1$ type or some $A_{n_i}$ type with involution".

\bigskip\noindent
{\bf 4.10. Corollary.} {\it Let $(G_0,X_0)$ be an adjoint Shimura pair. Then $(G_0,X_0)$ is the adjoint of a Shimura pair of PEL type if and only if each simple factor $(G,X)$ of $(G_0,X_0)$ is such that one of the following two conditions holds:

\medskip
{\bf (i)} the Shimura pair $(G,X)$ is of $A_n$ type, or

\smallskip
{\bf (ii)} the Shimura pair $(G,X)$ is of $C_n$, $D_n^{\dbH}$, or inner $D_4^{\dbR}$ type and the group $G_{\dbR}$ has no simple, compact factor.}

\medskip
\proof
We first prove the ``only if" part. Let $\tilde f:(\tilde G,\tilde X)\hookrightarrow (\text{\bf GSp}(\tilde W,\tilde\psi),\tilde S)$ be an embedding of PEL type. Let $(G,X)$ be a simple factor of $(\tilde G^{\ad},\tilde X^{\ad})$. Let $\tilde\scrB:=\End(\tilde W)^{\tilde G(\dbQ)}$. We show that the assumption that $\tilde\scrB$ has two distinct factors $\tilde\scrB_1$ and $\tilde\scrB_2$ which are permuted by the involution of $\End(\tilde W)$ defined by $\tilde\psi$, leads to a contradiction. Let $\tilde W=\tilde W_0\oplus\tilde W_1\oplus\tilde W_2$ be the unique direct sum decomposition in $\tilde\scrB$-modules such that $\tilde\scrB_1\oplus\scrB_2$ acts trivially on $\tilde W_0$ and for $i\in\{1,2\}$ the $\tilde\scrB_i$-module $\tilde W_i$ is a direct sum of simple $\tilde\scrB_i$-modules. The $\dbG_{m,\dbQ}$ subgroup of $\text{\bf GL}_{\tilde W}$ which fixes $\tilde W_0$, which acts on $\tilde W_1$ via the identical character of $\dbG_{m,\dbQ}$, and which acts on $\tilde W_2$ via the inverse of the identical character of $\dbG_{m,\dbQ}$, is contained in both $Z(\tilde G)$ and $\text{\bf Sp}(\tilde W,\tilde\psi)$. This contradicts the fact that the centralizer of each element $\tilde x\in\tilde X$ in $\text{\bf Sp}(\tilde W,\tilde\psi)_{\dbR}$ is a compact group.

Thus each simple factor of $\tilde\scrB$ is normalized by the involution of $\End(\tilde W)$ defined by $\tilde\psi$. Based on this, we can assume that $\tilde\scrB$ is a simple $\dbQ$--algebra; thus either $\tilde G^{\ad}=G$ or $\tilde G^{\der}_{\dbC}$ is a product of $\text{\bf SO}_{4,\dbC}$ groups. If $\tilde G^{\der}_{\dbC}$ is a product of $\text{\bf SO}_{4,\dbC}$ groups, then $(G,X)$ is of $A_1$ type. Thus we can assume that $\tilde G^{\ad}=G$ and that $(G,X)$ is not of $A_n$ type. From this and [De2, Table 2.3.8] we get that $(G,X)$ is of either $B_n$ type with $n\Ge 2$ or $C_n$ type with $n\Ge 2$ or $D_n^{\dbH}$ type with $n\Ge 4$ or $D_n^{\dbR}$ type with $n\Ge 4$. From [Sh, p. 320] we get that $\tilde G^{\der}_{\dbR}$ is a product of non-compact groups that are isomorphic to either $\text{\bf Sp}_{2n,\dbR}$ (with $n\Ge 2$) or $\text{\bf SO}^*(2n)_{\dbR}$ (with $n\Ge 4$); here the restrictions on $n$ are as per the type of $(G,X)$. Thus either (ii) holds or $(G,X)$ is of non-inner $D_4^{\dbR}$ type. But if $(G,X)$ is of $D_4^{\dbR}$ type, then $\tilde G_{\overline{\dbQ}}^{\der}$ is a product of $\text{\bf SO}_{8,\overline{\dbQ}}$ groups. If $(G,X)$ is of non-inner $D_4^{\dbR}$ type, then the representation of each simple factor of $\Lie(G_{\overline{\dbQ}})$ on $\tilde W\otimes_{\dbQ} \overline{\dbQ}$ involves both half spin representations (see Definition 2.2.1 (c)) and thus $\tilde G_{\overline{\dbQ}}^{\der}$ is a product of $\text{\bf Spin}_{8,\overline{\dbQ}}$ groups and therefore we reached a contradiction. Thus (ii) holds. This proves the ``only if" part.

Based on Fact 2.4.1 (b), it suffices to check the ``if" part only in the case when $(G_0,X_0)$ is simple. If $(G,X):=(G_0,X_0)$ is such that (i) (resp. (ii)) holds, then from Proposition 4.1 (resp. from Subsubsection 4.8.3), we get that $(G_0,X_0)$ is the adjoint of a Shimura pair of PEL type. \endproof

\bigskip\smallskip
\noindent
{\boldsectionfont 5. Applications to Frobenius tori}
\bigskip

We now apply Section 4 to the study of Frobenius tori of abelian varieties over number fields. At the very end of Subsections 5.1 to 5.8, we prove Theorem 1.3.2 and Corollary 1.3.3. In Subsections 5.9 and 5.10 we deal with complements on the shifting process mentioned in Subsection 1.6 and elaborated in Subsection 5.4; in particular, see Proposition 5.10 for three main advantages one gains in general by performing it.

We begin by introducing a setting which is independent of Section 4. Let $E$, $A$, $L_A$, $H_A$, $\rho$, $G_{\dbQ_p}$, and $E^{\text{conn}}$ be as in Subsection 1.1. Let $W_A$, $\lambda_A$, $\psi_A$, and $f_A$ be as in Subsection 1.5. In this Section we will assume that $\lambda_A$ is a principal polarization of $A$. Let
$$K(4):=\{g\in \text{\bf GSp}(L_A,\psi_A)(\dbZhat)|g\,\text{mod 4 is the identity}\}.$$
It is known that $\Sh(\text{\bf GSp}(W_A,\psi_A),S_A)/K(4)$ is the
$\dbQ$--scheme that parameterizes isomorphism classes of principally
polarized abelian schemes that are of relative dimension
$\dim_{\dbQ}(W_A)/2$ and that have level-4 symplectic similitude
structure, cf. [De1, Thm. 4.21]. Let $K_A(4):=H_A(\dbA_f)\cap K(4)$.
We have a natural finite morphism
$$\Sh(H_A,X_A)/K_A(4)\to\Sh(\text{\bf GSp}(W_A,\psi_A),S_A)/K(4)$$
of $\dbQ$--schemes (see [De1, Cor. 5.4]). Let $(\scrA,\lambda_{\scrA})$ be the pull back to $\Sh(H_A,X_A)/K_A(4)$ of the universal principally polarized abelian scheme over $\Sh(\text{\bf GSp}(W_A,\psi_A),S_A)/K(4)$. By replacing $E$ with a finite field extension of it, we can assume that $E=E^{\text{conn}}$ and that there exists a morphism
$$m\colon\Spec(E)\to\Sh(H_A,X_A)/K_A(4)$$
such that $m^*(\scrA,\lambda_{\scrA})$ is $(A,\lambda_A)$. As $\Sh(H_A,X_A)(\dbC)=H_A(\dbQ)\backslash (X_A\times H_A(\dbA_f))$
(see [De2, Prop. 2.1.10 and Cor. 2.1.11]), we can assume that the composite of the morphism $\Spec(\dbC)\to\Spec(E)$ with $m$ is the point
$$[h_A,1_{W_A}]\in \Sh(H_A,X_A)/K_A(4)(\dbC)=H_A(\dbQ)\backslash (X_A\times H_A(\dbA_f)/K_A(4));$$ here $1_{W_A}$ is the identity element of $H_A(\dbA_f)$. As $E=E^{\text{conn}}$, we have $\im(\rho)\leqslant G_{\dbQ_p}(\dbQ_p)$.

\bigskip\noindent
{\bf 5.1. Definition.} Let $B$ and $C$ be two abelian varieties over a field $K$ embeddable into $\dbC$. We fix an embedding $i_K:K\hookrightarrow\dbC$. Let $H_B$, $H_C$, and $H_{B\times_K C}$ be the Mumford--Tate groups of $B_{\dbC}$, $C_{\dbC}$, and $B_{\dbC}\times_{\dbC} C_{\dbC}$ (respectively). We say $B$ and $C$ are {\it adjoint-isogenous} if the natural projection epimorphisms $H_{B\times_K C}\twoheadrightarrow H_B$ and $H_{B\times_K C}\twoheadrightarrow H_C$ induce isomorphisms at the level of adjoint groups (thus we can identify $H_B^{\ad}=H_C^{\ad}=H_{B\times_K C}^{\ad}$).

\medskip\noindent
{\bf 5.1.1. Remark.} The choice of another embedding $K\hookrightarrow\dbC$ corresponds to the replacements of $H_B$, $H_C$, and $H_{B\times_K C}$ by suitable forms of them. Thus the fact that the abelian varieties $B$ and $C$ are or are not adjoint-isogenous does not depend on the choice of the embedding $i_K:K\hookrightarrow\dbC$.

\bigskip\noindent
{\bf 5.2. New injective maps.} Until Subsection 5.8 we will assume that the adjoint group $H^{\ad}_A$ is non-trivial. Let
$$f_2:(G_2,X_2)\hookrightarrow (G_3,X_3)\hookrightarrow (G_4,X_4)\operatornamewithlimits{\hookrightarrow}\limits^{f_4} (\text{\bf GSp}(W_1,\psi_1),S_1)$$
be injective maps of Shimura pairs constructed as in the Variants 4.9 (a), (b), or (c) but starting from the adjoint Shimura pair $(G_2^{\ad},X_2^{\ad})\operatornamewithlimits{=}\limits^{ID}(H_A^{\ad},X_A^{\ad})$.
Thus $(G_4,X_4)$ is unitary and the injective map $(G_4,X_4)\hookrightarrow (\text{\bf GSp}(W_1,\psi_1),S_1)$  of Shimura pairs is a PEL type embedding. Until Proposition 5.10 it is irrelevant if we use Variant 4.9 (a), (b), or (c) to construct these injective maps of Shimura pairs. Let
$$\scrB:=\End(W_1)^{G_4(\dbQ)}.$$
\indent
As the group $G_2^{\ad}(\dbQ)$ is dense in $G_2^{\ad}(\dbR)$, by replacing $f_2$ with its composite with an automorphism $(G_2,X_2)\arrowsim (G_2,X_2)$ defined by an element of $G_2^{\ad}(\dbQ)$, we can assume that the intersection $X_2\cap X_A$ (taken inside $X_2^{\ad}=X_A^{\ad}$) contains the element $h_A\in X_A$.

Let $L_1$ be a $\dbZ$-lattice of $W_1$ such that we get a perfect alternating form $\psi_1\colon L_1\otimes_{\dbZ} L_1\to\dbZ$. Let $K_2(4):=G_2(\dbA_f)\cap \{g\in \text{\bf GSp}(L_1,\psi_1)(\dbZhat)|g\,\text{mod 4 is the identity}\}$.
As before Definition 5.1, we get naturally a principally polarized abelian scheme $(\scrA_2,\lambda_{\scrA_2})$ over $\Sh(G_2,X_2)/K_2(4)$.

\bigskip\noindent
{\bf 5.3. Hodge twists.} Let $(W_0,\psi_0):=(W_A\oplus W_1,\psi_A\oplus\psi_1)$. Let $G_0$ be the identity component of the intersection of $\text{\bf GSp}(W_0,\psi_0)$ with the inverse image in $H_A\times_{\dbQ} G_2$ of the subgroup $H_A^{\ad}$ of $H_A^{\ad}\times_{\dbQ} G_2^{\ad}$ embedded in $H_A^{\ad}\times_{\dbQ} G_2^{\ad}$ diagonally via $ID$. Let $X_0$ be the $G_0(\dbR)$-conjugacy class of the monomorphism $\dbS\hookrightarrow G_{0,\dbR}$ defined by $(h_A,h_A)\in X_A\times X_2$ (this makes sense due to $ID$ and the relation $h_A\in X_A\cap X_2$). We get a (diagonal) injective map
$$f_0\colon (G_0,X_0)\hookrightarrow (\text{\bf GSp}(W_0,\psi_0),S_0)$$
of Shimura pairs that is the composite of the injective map $f_3\colon (G_0,X_0)\hookrightarrow (H_A,X_A)\times^{\scrH} (G_2,X_2)$
 of Shimura pairs with a Hodge quasi product $f_A\times^{\scrH} f_2$.
The composite of $f_3$ with the natural injective map $(H_A,X_A)\times^{\scrH} (G_2,X_2)\hookrightarrow (H_A,X_A)\times (G_2,X_2)$ of Shimura pairs and with the projections of $(H_A,X_A)\times (G_2,X_2)$ onto its factors, give birth at the level of adjoint Shimura pairs to identifications
$$(G_0^{\ad},X_0^{\ad})\operatornamewithlimits{=}\limits^{ID}(G_2^{\ad},X_2^{\ad})\operatornamewithlimits{=}\limits^{ID}(H_A^{\ad},X_A^{\ad}).$$
\noindent
Due to these identifications, we call $(G_0,X_0)$ a {\it Hodge twist} of $(H_A,X_A)$ and $(G_2,X_2)$.

Let $L_0:=L_A\oplus L_1$. It is a $\dbZ$-lattice of $W_0$ and we have a perfect alternating form $\psi_0\colon L_0\otimes_{\dbZ} L_0\to\dbZ$. Let $K_0(4):=G_0(\dbA_f)\cap (K(4)\times K_2(4))$. As before Definition 5.1, we get naturally (starting from $f_0$ and $L_0$) a principally polarized abelian scheme $(\scrA_0,\lambda_{\scrA_0})$ over $\Sh(G_0,X_0)/K_0(4)$.

\bigskip\noindent
{\bf 5.4. The shifting process.} Let $B$ be an abelian variety over $E$ which is adjoint-isogenous to $A$. Let $H_B$ be the Mumford--Tate group of $B_{\dbC}$; we have a natural identification $H_B^{\ad}=H_A^{\ad}$. Let $\tilde G_{\dbQ_p}$ be the identity component of the algebraic envelope of the $p$-adic Galois representation attached to $B$; it is a subgroup of $H_{B,\dbQ_p}$.

\medskip\noindent
{\bf 5.4.1. Proposition.} {\it The images of $G_{\dbQ_p}^{\der}$ and $\tilde G_{\dbQ_p}^{\der}$ in $H^{\ad}_{A,\dbQ_p}=H_{B,\dbQ_p}^{\ad}$ coincide. Thus to prove the part of the Mumford--Tate conjecture that involves derived groups for the pair $(A,p)$ is the same thing as proving it for the pair $(B,p)$ (i.e., we have $G_{\dbQ_p}^{\der}=H^{\der}_{A,\dbQ_p}$  if and only if $\tilde G_{\dbQ_p}^{\der}=H_{B,\dbQ_p}^{\der}$).}

\medskip
\proof
The images of $G_{\dbQ_p}^{\der}$ and $\tilde G_{\dbQ_p}^{\der}$ in $H^{\ad}_{A,\dbQ_p}=H_{B,\dbQ_p}^{\ad}$ coincide with the image of the identity component of the algebraic envelope of the $p$-adic Galois representation attached to $A\times_E B$ in $H_{A\times_E B,\dbQ_p}^{\ad}=H^{\ad}_{A,\dbQ_p}=H_{B,\dbQ_p}^{\ad}$. Thus we have $G_{\dbQ_p}^{\der}=H^{\der}_{A,\dbQ_p}$  if and only if $\tilde G_{\dbQ_p}^{\der}=H_{B,\dbQ_p}^{\der}$. \endproof

\medskip\noindent
{\bf 5.4.2. Choice of $B$.} As $h_A\in X_2\cap X_A$, as $(h_A,h_A)\in X_0$, and as the functorial morphism $\Sh(G_0,X_0)/K_0(4)\to \Sh(H_A,X_A)/K_A(4)$ is finite, by replacing $E$ with a finite field extension of it we can assume that there exists a morphism
$$m_0\colon\Spec(E)\to\Sh(G_0,X_0)/K_0(4)$$
such that the composite of the natural morphism $\Spec(\dbC)\to\Spec(E)$ with $m_0$ is the point $[(h_A,h_A),1_{W_0}]\in \Sh(G_0,X_0)/K_0(4)(\dbC)$.

Let $B$ be the abelian variety over $E$ obtained from $\scrA_2$ by pull back via $m_0$ and the natural morphism $\Sh(G_0,X_0)/K_0(4)\to \Sh(G_2,X_2)/K_2(4)$. We have an identity
$$m_0^*(\scrA_0)=A\times_E B.$$
The Mumford--Tate group of $A_{\dbC}\times_{\dbC} B_{\dbC}$ (resp. of $B_{\dbC}$) is a subgroup of $G_0$ (resp. of $G_2$) that surjects onto $G_0^{\ad}=G_2^{\ad}=H_A^{\ad}$. This implies that the abelian varieties $A$ and $B$ are adjoint-isogenous and that $\tilde G_{\dbQ_p}$ is a subgroup of $H_{B,\dbQ_p}$ and thus also of $G_{2,\dbQ_p}$.

\medskip\noindent
{\bf 5.4.3. Decomposing $B$.}
See Subsection 5.2 for the definition of $\scrB$. The injective map $f_4:(G_4,X_4)\hookrightarrow (\text{\bf GSp}(W_1,\psi_1),S_1)$ of Shimura pairs is a Hodge quasi product indexed by the set $J$ of Subsection 1.3, cf. property (iv) of Variant 4.9 (a). Thus from the simple $\dbQ$--algebra statements of the property 4.1.1 (a) and of Subsubsection 4.1.2, we get that the number $s_{\scrB}$ of simple factors of $\scrB$ equals to the number of elements of $J$. If we take $L_1$ to be well adapted for $f_4$ (see Subsection 2.4), then $B$ is a product indexed by $J$ of abelian varieties over $E$ that have the following two properties: (i) they are principally polarized, and (ii) their extensions to $\dbC$ have Mumford--Tate groups whose adjoints are simple.

\bigskip\noindent
{\bf 5.5. Good choices for $v$.} Let $v$, $l$, $A_v$, and $T_v$ be as before Theorem 1.1.5. We choose the prime $v$ of $E$ such that moreover we have $l>\dim(B)$, $l\neq p$, and the following two conditions hold:

\medskip
{\bf (i)} the abelian variety $B$ has good reduction $B_v$ with respect to $v$;

\smallskip
{\bf (ii)} the prime $v$ is unramified over $l$ and $\scrB_{(l)}:=\scrB\cap\End(L_1\otimes_{\dbZ} \dbZ_{(l)})$ is a semisimple $\dbZ_{(l)}$-algebra that is self dual with respect to the perfect alternating $\psi_1$ on $L_1$.

\medskip
Let $T_{v,2}$ and $T_{v,0}$ be the Frobenius tori of $B_v$ and $A_v\times_{k(v)} B_v$ (respectively). Let $\pi_v\in T_v(\dbQ_p)$, $\pi_{v,0}\in T_{v,0}(\dbQ_p)$, and $\pi_{v,2}\in T_{v,2}(\dbQ_p)$ be (the same integral power of) Frobenius elements obtained as in Subsection 3.2 via the choice of a Frobenius automorphism of a prime $\overline v$ of $\overline E=\overline{\dbQ}$ that divides $v$. We have natural monomorphisms $T_{v,\dbQ_p}\hookrightarrow G_{\dbQ_p}$, $T_{v,2,\dbQ_p}\hookrightarrow\tilde G_{\dbQ_p}\hookrightarrow G_{2,\dbQ_p}$, and $T_{v,0,\dbQ_p}\hookrightarrow G_{0,\dbQ_p}$. The images of $\pi_v$, $\pi_{v,0}$, and $\pi_{v,2}$ in $H^{\ad}_{A,\dbQ_p}(\dbQ_p)=G_{0,\dbQ_p}^{\ad}(\dbQ_p)=G_{2,\dbQ_p}^{\ad}(\dbQ_p)$ coincide.

By replacing $E$ with a finite field extension of it, we can assume that the image of the $p$-adic Galois representation $\rho_2:\Gal(E)\to \text{\bf GL}_{L_1\otimes_{\dbZ} \dbQ_p}(\dbQ_p)$ attached to $B$ factors through $G_4(\dbQ_p)$ (i.e., all $\dbQ$--endomorphisms of $B_{\dbC}$ defined by elements of $\scrB$, are defined over $E$). Each $\dbQ$--endomorphism of $B$ defined by an element of $\scrB$ defines naturally a $\dbQ$--endomorphism of $B_v$. Let $\lambda_B$ be the principal polarization of $B$ that is the natural pull back of $\lambda_{\scrA_2}$. Let $\lambda_{B_v}$ be the principal polarization of $B_v$ that is defined naturally by $\lambda_B$.

\medskip\noindent
{\bf 5.5.1. Lemma.} {\it There exist a finite field extension $k(v^\prime)$ of $k(v)$, a complete discrete valuation ring $R$ of mixed characteristic $(0,l)$ and of residue field $k(v^\prime)$, and a principally polarized abelian scheme $(B^\prime,\lambda_{B^\prime})$ over $R$, such that the following three properties hold:

\medskip
{\bf (i)} the special fibre $(B^\prime_{v^\prime},\lambda_{B^\prime_{v^\prime}})$ of $(B^\prime,\lambda_{B^\prime})$ is $\dbZ[{1\over l}]$-isomorphic to $(B_{v^\prime},\lambda_{B_{v^\prime}}):=(B_v,\lambda_{B_v})_{k(v^\prime)}$ (i.e., it is obtained from $(B_{v^\prime},\lambda_{B_{v^\prime}})$ via isogenies of order a power of $l$);

\smallskip
{\bf (ii)} the abelian scheme $B^{\prime}$ has complex multiplication;

\smallskip
{\bf (iii)} the $\dbQ$--endomorphisms of $B_{v^\prime}^{\prime}$ that correspond to elements of $\scrB$, lift to $\dbQ$--endomorphisms of $B^{\prime}$ in such a way that the elements of $\scrB_{(l)}$ correspond to $\dbZ_{(l)}$-endomorphisms of $B^{\prime}$ and the trace forms on $\scrB_{(l)}\otimes_{\dbZ_{(l)}} k(v^\prime)$ defined by its representations on the tangent spaces of $B_{v^\prime}$ and $B^\prime_{v^\prime}$ at identity elements, are equal.}

\medskip
\proof
By using a product decomposition of $B$ (and thus implicitly of $(L_1,\psi_1,\scrB_{(l)})$) as in Subsubsection 5.4.3, to prove the Lemma we can assume that $\scrB$ is a simple $\dbQ$--algebra. Let $\scrI_1$ be the involution of $\End(W_1)$ defined by the identity $\psi_1(e(x),y)=\psi_1(x,\scrI_1(e)(y))$, where $e\in \End(W_1)$ and $x,y\in W_1$. As $G_4$ is a subgroup of $\text{\bf GSp}(W_1,\psi_1)$, the involution $\scrI_1$ normalizes $\scrB$ and thus also $\scrB_{(l)}$. As the injective map $f_4:(G_4,X_4)\hookrightarrow (\text{\bf GSp}(W_1,\psi_1),S_1)$ of Shimura pairs is a PEL type embedding, the group $G_4$ is the centralizer of $\scrB$ in $\text{\bf GSp}(W_1,\psi_1)$. From this and the fact that the Shimura pair $(G_4,X_4)$ is unitary, we get that the restriction of the involution $\scrI_1$ to $\scrB$ is of type $A$ in the sense of [Zi, Subsect. 2.3]. Thus the Lemma is a particular case of [Zi, Thm. 4.4] applied to the data provided by the injective map $f_4:(G_4,X_4)\hookrightarrow (\text{\bf GSp}(W_1,\psi_1),S_1)$ of Shimura pairs, the perfect alternating form $\psi_1:L_1\otimes_{\dbZ} L_1\to \dbZ$, the involution $\scrI_1$, the prime $l$, the semisimple $\dbZ_{(l)}$-algebra $\scrB_{(l)}$, and the principally polarized abelian variety $(B_v,\lambda_{B_v})$ over $k(v)$ which is the good reduction of the pull back $(B,\lambda_B)$ of $(\scrA_2,\lambda_{\scrA_2})$. Strictly speaking loc. cit. states only that $(B^\prime_{v^\prime},\lambda_{B^\prime_{v^\prime}})$ is obtained from $(B_{v^\prime},\lambda_{B_{v^\prime}})$ via isogenies. But one can always assume that these isogenies are of order a power of $l$.\endproof

\bigskip\noindent
{\bf 5.6. Applications of 5.5.1.} We can assume that the field $R[{1\over l}]$ contains $E$. We fix an embedding $R[{1\over p}]\hookrightarrow\dbC$ that extends the embedding $i_E$ of Subsection 1.1. Let $W_1^\prime:=H_1(B^\prime(\dbC),\dbQ)$, let $\psi_1^\prime$ be the alternating form on $W_1^\prime$ induced by $\lambda_{B^{\prime}}$, and let $\scrB\hookrightarrow\End(W_1^\prime)$ be the $\dbQ$--monomorphism defined by the identification (see property 5.5.1 (iii)) of $\scrB$ with a $\dbQ$--subalgebra of $\End(B_{\dbC}^\prime)\otimes_{\dbZ} \dbQ$. Let $G_4^\prime$ be the subgroup of $\text{\bf GSp}(W_1^\prime,\psi_1^\prime)$ that fixes all elements of $\scrB$. From the property 5.5.1 (i) we get that for each prime $\tilde l\in\dbN$ different from $l$, the two triples $(W_1^\prime,\psi_1^\prime,\scrB)$ and $(W_1,\psi_1,\scrB)$ become isomorphic over $\dbQ_{\tilde l}$. Thus the ``difference" between the isomorphism classes of the two triples $(W_1^\prime,\psi_1^\prime,\scrB)$ and $(W_1,\psi_1,\scrB)$ is ``measured" by a class $\delta_4\in H^1(\dbQ,G_4)$ whose image in $H^1(\dbQ_{\tilde l},G_{4,\dbQ_{\tilde l}})$ is trivial for all primes $\tilde l$ different from $l$. In particular, we get that $G_4^\prime$ is an inner form of $G_4$.

\medskip\noindent
{\bf 5.6.1. Fact.} {\it There exists a reductive subgroup $G_4^{(1)}$ of $\text{\bf GL}_{W_1}$ that is the extension of $\dbG_{m,\dbQ}^{s_{\scrB}-1}$ by $G_4$ and such that the images of the class $\delta_4$ in $H^1(\dbR,G^{(1)}_{4,\dbR})$ and $H^1(\dbQ_q,G^{(1)}_{4,\dbQ_q})$ are trivial, for all primes $q\in\dbN$ (see Subsubsection 5.4.3 for the definition of $s_{\scrB}$).}

\medskip
\proof
By using a product decomposition of $B$ (and thus implicitly of $(L_1,\psi_1,\scrB_{(l)})$) as in Subsubsection 5.4.3 and due to the first part of Fact 2.4.1 (a), to prove the Lemma we can assume that $\scrB$ is a simple $\dbQ$--algebra. Thus $s_{\scrB}=1$ and $G_4^{(1)}=G_4$. We already know that the image of $\delta_4$ in $H^1(\dbQ_q,G_{4,\dbQ_q})$ is trivial provided $q\neq l$. The fact that the images of $\delta_4$ in both pointed sets $H^1(\dbR,G^{(1)}_{4,\dbR})$ and $H^1(\dbQ_l,G^{(1)}_{4,\dbQ_l})$ are trivial is checked in the third and fourth (respectively) paragraphs of [Ko2, Sect. 8]. To loc. cit. we only have to add that the inequality $l>\dim(B)$ implies that it is irrelevant if in the context of the data provided in the proof of Lemma 5.5.1, we work with the determinant condition of [Ko2, Sect. 5] or with the trace condition of the property 5.5.1 (iii) as in [Zi].\endproof

\medskip
From Fact 5.6.1 we get that we have an identification $G_{4,\dbQ_p}=G_{4,\dbQ_p}^{\prime}$ (unique up to inner conjugation) as well as a canonical identification $Z(G_4)=Z(G_4^\prime)$.

Let $(G_4^\prime,X_4^\prime)$ be the Shimura pair, where $X_4^\prime$ is the $G_4^\prime(\dbR)$-conjugacy class of the homomorphism $\dbS\to G_{4,\dbR}^\prime$ that defines the Hodge $\dbQ$--structure on $W_1^\prime$. If we have $s_{\scrB}=1$, then there exists an isomorphism $W_1\otimes_{\dbQ} \dbR\arrowsim W_1^\prime\otimes_{\dbQ} \dbR$ which (cf. [Ko2, Lem. 4.2 and the third paragraph of Sect. 8]) has the following two properties:

\medskip
{\bf (a)} it takes $\psi_1$ to $\psi_1^\prime$ and it takes $b$ to $b$ for all $b\in\scrB$ (thus it takes the subgroup $G_{4,\dbR}$ of $\text{\bf GL}_{W_1\otimes_{\dbQ} \dbR}$ onto the subgroup $G_{4,\dbR}^\prime$ of $\text{\bf GL}_{W_1^\prime\otimes_{\dbQ} \dbR}$);

\smallskip
{\bf (b)} it takes $X_4$ onto $X_4^\prime$.

\medskip
Thus for $s_{\scrB}=1$, the canonical identification $Z(G_4)=Z(G_4^\prime)$ gives birth naturally (cf. (b)) to a canonical identification $(G_4^{\ab},X_4^{\ab})=(G_4^{\prime,\ab},X_4^{\prime,\ab})$ of $0$ dimensional Shimura pairs. Based on this and the product decomposition of Subsubsection 5.4.3, regardless of what the number $s_{\scrB}\in\dbN$ is, we get:

\medskip\noindent
{\bf 5.6.2. Fact.} {\it Let $Z_{\text{min}}$ be the smallest subtorus of $Z^0(G_4)$ with the property that all monomorphisms $\dbS\hookrightarrow G_{4,\dbR}$ defined by elements of $X_4$, factor through the extension to $\dbR$ of the subgroup of $G_4$ generated by $G_4^{\der}$ and $Z_{\text{min}}$. Let $Z_{\text{min}}^\prime$ be the subtorus of $Z^0(G_4^\prime)$ that is defined similarly to $Z_{\text{min}}$. Then the canonical identification $Z(G_4)=Z(G_4^\prime)$ restricts to a canonical identification $Z_{\text{min}}=Z_{\text{min}}^\prime$.}

\medskip\noindent
{\bf 5.6.3. Betti realization of $\pi_{v,2}$.} We use the above identification $G_{4,\dbQ_p}=G_{4,\dbQ_p}^\prime$. Coming back to the Frobenius elements of Subsection 5.5, from properties 5.5.1 (ii) and (iii) we get (cf. the third paragraph of Subsection 3.2 and the relation $\im(\rho_2)\leqslant G_4(\dbQ_p)$) the existence of a torus $T_4^\prime$ of $G_4^\prime$ and of an element $\pi_{v^\prime,2}^\prime\in T_4^\prime(\dbQ)$ such that the following three properties hold:

\medskip
{\bf (i)} the element $\pi_{v^\prime,2}^\prime$ is the Betti realization of an endomorphism of $B^\prime$ which lifts a suitable integral power of the Frobenius endomorphism of $B^\prime_{k(v^\prime)}$;

\smallskip
{\bf (ii)} the torus $T_4^{\prime}$ is the smallest subgroup of $G_4^\prime$ which has $\pi_{v^\prime,2}^\prime$ as a $\dbQ$--valued point;

\smallskip
{\bf (iii)} there exists an element $g^\prime\in G_4^\prime(\dbQ_p)=G_4(\dbQ_p)$ that takes under conjugation $T_{4,\dbQ_p}^\prime$ to $T_{v,2,\dbQ_p}$ and $\pi_{v^\prime,2}^\prime$ to $\pi_{v,2}$. [Here we have to introduce the element $g^\prime$ as the torus $T_4^\prime$ of $G_4^\prime$ is uniquely determined only up to $G_4^\prime(\dbQ)$-conjugation and as $T_{v,2,\dbQ_p}$ and $\pi_{v,2}$ are uniquely determined (cf. Subsection 3.2) only up to $\im(\rho_2)$-conjugation].

\medskip\indent
The images of $T_{v,\overline{\dbQ_p}}$ and $T_{4,\overline{\dbQ_p}}^\prime$ in $G^{\ad}_{4,\overline{\dbQ_p}}=G^{\prime,\ad}_{4,\overline{\dbQ_p}}$ are $G_4^{\ad}(\overline{\dbQ_p})$-conjugate, cf. (iii) and the end of Subsection 5.5.

\medskip\noindent
{\bf 5.6.4. A crystalline property.} As $R$ is a complete discrete valuation ring of mixed characteristic $(0,l)$ and of residue field $k(v^\prime)$, it has a canonical structure of a $W(k(v^\prime))$-algebra. Let $(M,\phi)$ be the $F$-crystal over $k(v)$  of $B_v$. Thus the pair $(M\otimes_{W(k(v))} B(k(v^\prime)),\phi)$ is the $F$-isocrystal over $k(v^\prime)$ of $B^\prime_{v^\prime}$. The canonical and functorial identification $M\otimes_{W(k(v))} R[{1\over p}]=H^1_{\text{dR}}(B^\prime/R)[{1\over p}]$ (see [BO, Thm. 1.3]) allows us to identify the de Rham realization of the endomorphism of $B^\prime$ mentioned in the property 5.6.3 (i) with a suitable integral power of $\phi^{[k(v^\prime):\dbF_p]}$. From this and the property 3.4.1 (c) applied in the context of the pair $(B,v)$ (instead of $(A,v)$), we get that a suitable $G_4^\prime(\overline{\dbQ_l})$-conjugate of $T_{4,\overline{\dbQ_l}}^\prime$ is a torus of $\tilde G_{\overline{\dbQ_l}}$.

\bigskip\noindent
{\bf 5.7. Simple properties.}
{\bf (a)} Each simple factor of $\Lie(G_{4,\dbR}^{\der})$ is absolutely simple and only one simple factor of $\Lie(G_{2,\dbR}^{\der})$ maps injectively into it. Each simple factor of $\Lie(G_{2,\dbR}^{\der})$ is embedded (``diagonally") into a product of simple factors of $\Lie(G_{3,\dbR}^{\der})$ isomorphic to it. The same holds over $\overline{\dbQ_p}$ instead of over $\dbR$. Moreover, each simple factor of $\Lie(G_2^{\ad})$ projects injectively into a unique simple factor of $\Lie(G_4^{\ad})$. All these follow from the property (iii) of Variant 4.9 (a) and from Variants 4.9 (b) and (c).

\smallskip
{\bf (b)} We use the identification $G_{4,\overline{\dbQ_p}}=G_{4,\overline{\dbQ_p}}^{\prime}$ (see paragraph after Fact 5.6.1). Two simple factors of $G_{4,\overline{\dbQ_p}}^{\ad}$ are connected over $\dbQ$ through $G_4^{\ad}$ (i.e., are factors of the extension to $\overline{\dbQ_p}$ of the same simple factor of $G_4^{\ad}$) if and only if they are connected over $\dbQ$ through $G_4^{\prime{,\ad}}$. This is so as $G_4^\prime $ is an inner form of $G_4$.

\smallskip
{\bf (c)} If $s_{\scrB}=1$ and $(G_4^{\ad},X_4^{\ad})$ is a simple Shimura pair of some $A_n$ type with $n$ odd, then the Hasse principle holds for $G_4$ (cf. [Ko2, Sect. 7]). Thus the class $\delta_4\in H^1(\dbQ,G_4)$ is trivial, cf. Fact 5.6.1.

\smallskip
{\bf (d)} If all simple factors of $(G_4^{\ad},X_4^{\ad})$ are of some $A_n$ type with $n$ odd and if in the property 5.5.1 (i) we require only that $B^\prime_{v^\prime}$ is isogenous (instead of being $\dbZ[{1\over l}]$-isomorphic) to $B_{v^\prime}$ and that properties 5.5.1 (ii) and (iii) hold, then one can easily check based on (c) that we can choose $k(v^\prime)$, $B^\prime_{v^{\prime}}$, $R$, and $B^\prime$ in such a way that $\delta_4$ is the trivial class (regardless of what $s_{\scrB}\in\dbN$ is). The last sentence applies if there exists no element $j\in J$ such that $(H_j,X_j)$ is of $A_n$ type with $n$ even, cf. the definition of the group $G(A)$ in Subsection 4.3 and the property 4.8 (iii).

\bigskip\noindent
{\bf 5.8. Proofs of 1.3.2 and 1.3.3.} We are ready to prove Theorem 1.3.2 and Corollary 1.3.3. We can assume that the adjoint group $H_A^{\ad}$ is non-trivial. Let $N(A)\in\dbN$ be a number that depends only on $A$ and that has the property that for each prime $l\Ge N(A)$, both conditions 5.5 (i) and (ii) hold and moreover $A$ has good reduction with respect to all primes of $E$ that divide $l$. By replacing $E$ with a finite field extension of it, to prove Theorem 1.3.2 and Corollary 1.3.3 we can assume that $A$ has a principal polarization $\lambda_A$ (cf. [Mu2, Ch. IV, Sect. 23, Cor. 1]) and that (as $l\Ge N(A)$) we are in the context of Subsections 5.2 to 5.7. As $T_4^\prime$ is a torus of $G_4^\prime$, the images of $T_{4,\overline{\dbQ_p}}^\prime$ in simple factors of the extension to $\overline{\dbQ_p}$ of a fixed simple factor of $G_4^{\prime,\ad}$ (equivalently of $G_4^{\ad}$, cf. property 5.7 (b)), have equal ranks that do not depend on $p\neq l$. Thus also the images of $T_{v,\overline{\dbQ_p}}$ in simple factors of the extension to $\overline{\dbQ_p}$ of a fixed simple factor of $H_A^{\ad}=G_2^{\ad}$, have equal ranks (cf. end of Subsubsection 5.6.3 and the property 5.7 (a)) that do not depend on $p\neq l$. This proves Theorem 1.3.2.

In order to get Corollary 1.3.3 from Theorem 1.3.2, we choose the prime $v$ of $E$ such that $l\Ge N(A)$ and $T_v$ has the same rank as $G_{\dbQ_p}$. Thus the fact that Corollary 1.3.3 holds for $p\neq l$ follows from Theorem 1.3.2. As in the end of Lemma 3.3, using a prime $v_1$ of $E$ that divides a prime $l_1\in\dbN$ that is greater than $N(A)$ and different form $l$, we get that Corollary 1.3.3 holds even if $p=l$. This proves Corollary 1.3.3.\endproof

\bigskip\noindent
{\bf 5.9. Definition.} An irreducible representation of a split, simple Lie algebra of classical Lie type $\grL$ over a field $K$ of characteristic 0 is called an {\it SD-standard representation}, if the following two conditions hold: (i) its highest weight is a minuscule weight, and (ii) if $\grL=A_n$, then moreover its dimension is $n+1$ (i.e., its highest weight is either $\varpi_1$ or $\varpi_n$). Here SD stands for Satake--Deligne, cf. [Sa2] and [De2, Prop. 2.3.10], [Va1, Subsects. 6.5 and 6.6], and the proofs of Proposition 4.1 and Theorem 4.8. Similarly, we speak about SD-standard representations of split, semisimple groups of classical Lie type over $K$ whose adjoints are simple.

\medskip
Until the end, we continue to assume that the adjoint group $H_A^{\ad}$ is non-trivial. Next we list some extra advantages we get (besides proving Theorem 1.3.2) by replacing the abelian variety $A$ and the injective map $f_A:(H_A,X_A)\hookrightarrow (\text{\bf GSp}(W_A,\psi_A),S_A)$ of Shimura pairs with $B$ and $f_2:(G_2,X_2)\hookrightarrow (\text{\bf GSp}(W_1,\psi_1),S_1)$ (respectively).

\bigskip\noindent
{\bf 5.10. Proposition.} {\bf (a)} {\it The irreducible subrepresentations of the representation of $\Lie(G_{2,\overline{\dbQ}}^{\der})$ on $W_1\otimes_{\dbQ} \overline{\dbQ}$ factor through simple factors of $\Lie(G_{2,\overline{\dbQ}}^{\der})$ and are SD-standard representations of these simple factors.

\smallskip
{\bf (b)} By using Variant 4.9 (a) or (c) (resp. 4.9 (c)) to construct the injective maps of Shimura pairs of Subsection 5.2, we can assume that for each simple factor $(H_j,X_j)$ of $(H_A^{\ad},X_A^{\ad})=(G^{\ad}_2,X^{\ad}_2)$ of $D_n^{\dbR}$ type with $n\Ge 4$ (resp. of $A_n$ type without involution and with $n\Ge 2$) and every simple factor $\scrF_j$ of $H_{j,\overline{\dbQ}}$, there exists a $G_{2,\overline{\dbQ}}$-submodule $\scrW$ of $W_1\otimes_{\dbQ} \overline{\dbQ}$ that has the following two properties:

\medskip
{\bf (i)} the natural representation $Z^0(G_{2,\overline{\dbQ}})\to \text{\bf GL}_{\scrW}$ factors through $Z(\text{\bf GL}_{\scrW})$ (i.e., $Z^0(G_{2,\overline{\dbQ}})$ acts on $\scrW$ via scalar automorphisms), and

\smallskip
{\bf (ii)} both half spin representations (resp. both irreducible representations of the highest weights $\varpi_1$ and $\varpi_{n}$) of $\Lie(\scrF_j)$ are among the simple $\Lie(G_{2,\overline{\dbQ}}^{\ad})$-submodules of $\scrW$.

\medskip
{\bf (c)} If $H^{\ad}_{A,\dbR}$ has no simple, compact factor and if $(H_A^{\ad},X_A^{\ad})$ has no simple factor of some $A_n$ type with involution, then we can also choose $f_2:(G_2,X_2)\hookrightarrow (\text{\bf GSp}(W_1,\psi_1),S_1)$ such that we have $Z^0(G_2)=Z(\text{\bf GL}_{W_1})\arrowsim\dbG_{m,\dbQ}$.}

\medskip
\proof
Part (a) is a consequence of the proof of Proposition 4.1 and of Subsection 4.3, cf. property (iii) of Variant 4.9 (a) and Variants 4.9 (b) and (c). Part (b) for the $D_n^{\dbR}$ type with $n\Ge 4$ (resp. the $A_n$ type without involution and with $n\Ge 2$) is a consequence of the fact that in Subsection 4.3 we defined $h_{\dbR}$ using both half spin representations (resp. of the first paragraph of Subsubsection 4.1.2). Based on Fact 2.4.1 (a), to check (c) we can assume that the adjoint Shimura pair $(H_A^{\ad},X^{\ad}_A)$ is simple. Thus, if $(H_A^{\ad},X^{\ad}_A)$ is not (resp. is) of $A_n$ Lie type, then (c) is a consequence of Remark 4.1.3 (resp. of the first paragraph of Subsubsection 4.8.3).\endproof

\medskip\noindent
{\bf 5.10.1. Remark.}
By replacing $(A,f_A)$ with $(B,f_2)$, the role of the Mumford--Tate group $H_A$ gets replaced by the one of a Mumford--Tate group which is a reductive, normal subgroup of $G_2$ that contains $G_2^{\der}$ (more precisely by the one of the smallest subgroup of $G_2$ with the property that all homomorphisms $\dbS\to G_{2,\dbR}$ defined by elements of $X_2$ factor through its extension to $\dbR$). Thus under such a replacement we lose to a great extent the ``control'' on the torus $Z(H_A)$; moreover, we often have $\dim(B)>\dim(A)$. If this looks ``unpleasant" or if one would like to ``stick" to the injective map $f_A$ of Shimura pairs, then one can consider the PEL-envelope $(\tilde H_A,\tilde X_A)$ of $f_A$ as defined in [Va1, Rm. 4.3.12]. In other words, one can consider natural injective maps
$$f_A:(H_A,X_A)\hookrightarrow (\tilde H_A,\tilde X_A)\operatornamewithlimits{\hookrightarrow}\limits^{\tilde f_A} (\text{\bf GSp}(W_A,\psi_A),S_A)$$
 of Shimura pairs, where $\tilde H_A$ is the identity component of the intersection of $\text{\bf GSp}(W_A,\psi_A)$ with the double centralizer of $H_A$ in $\text{\bf GL}_{W_A}$. Thus $\tilde f_A$ is a PEL type embedding and we get a variant of Subsubsection 5.6.3 for its context: if the analogues of the conditions 5.5 (i) and (ii) hold in the context of $A$ and $\tilde f_A$ (instead of $B$ and $f_2$), then Subsubsections 5.5.1 to 5.6.3 can be adapted to give us that for $l\gg 0$ we can view $\pi_v$ as a $\dbQ$--valued point of a torus $\tilde T_A^\prime$ of a suitable form $\tilde H_A^\prime$ of $\tilde H_A$. This form is inner provided $\tilde H_{A,\dbR}^{\der}$ has no normal subgroup that is an $\text{\bf SO}^*(2n)_{\dbR}$ group with $n\Ge 2$, cf. the connectedness aspects of [Ko2, Sect. 7]. We do not know how one could get Theorem 1.3.2 via $\tilde T_A^\prime$ (i.e., via using $\tilde f_A$ instead of $f_2$ of Subsection 5.2).

Let $\tilde C_A$ be the centralizer of $Z^0(H_A)$ in $\text{\bf GSp}(W_A,\psi_A)$. Let $\tilde H_{A,1}^{\ad}$ be the maximal factor of $\tilde C_A^{\ad}$ with the property that no simple factor of it becomes compact over $\dbR$. Let $\tilde H_{A,1}$ be the maximal normal, reductive subgroup of $\tilde C_A$ whose adjoint is $\tilde H_{A,1}^{\ad}$. We have a variant of the previous paragraph in which we work with $\tilde H_{A,1}$ instead of $\tilde H_A$. The monomorphism $\tilde H_{A,1}\hookrightarrow \text{\bf GSp}(W_A,\psi_A)$ extends uniquely to an injective map $\tilde f_{A,1}:(\tilde H_{A,1},\tilde X_{A,1})\hookrightarrow (\text{\bf GSp}(W_A,\psi_A),S_A)$ of Shimura pairs through which $\tilde f_A$ factors. The injective map $\tilde f_{A,1}$ of Shimura pairs is a PEL type embedding and thus for $l\gg 0$ we can view $\pi_v$ as a $\dbQ$--valued point of a torus $\tilde T_{A,1}^\prime$ of a suitable form of $\tilde H_{A,1}$. As the group $\tilde C_{A,\overline{\dbQ}}^{\der}$ (and thus also $\tilde H_{A,1,\overline{\dbQ}}^{\ad}$) has no normal subgroup which is an  $\text{\bf SO}$ group, this form is inner. Also, one can use the torus $\tilde T_{A,1}^\prime$ and [Zi] in a way similar to Subsections 5.5 and 5.6, to re-obtain Remark 3.3.1 for $l\gg 0$.

\medskip\noindent
{\bf 5.10.2. Remark.}
Different variants of Proposition 5.10 can be obtained by working with suitable other abelian varieties over $E$ that are adjoint-isogenous to $A$ and that are obtained similarly to $B$ but via another injective map $\tilde f_2\colon (\tilde G_2,\tilde X_2)\hookrightarrow (\text{\bf GSp}(\tilde W_1,\tilde\psi_1),\tilde S_1)$ of Shimura pairs for which we have an identity $(\tilde G_2^{\ad},\tilde X_2^{\ad})=(H_A^{\ad},X_A^{\ad})$. Thus for getting Proposition 5.10, we could have only quoted [De2, Subsubsects. 2.3.10 to 2.3.13]. Also, by quoting loc. cit. and by considering PEL-envelopes as in Remark 5.10.1 one gets directly: (i) less explicit forms of Theorem 4.8 provided the simple, adjoint Shimura pair $(G,X)$ of Section 4 is of $B_n$, $C_n$, or $D_n^{\dbH}$ type, and (ii) less explicit and weaker variants of Theorem 4.8 and Remark 4.8.1 if $(G,X)$ is of $D_n^{\dbR}$ type. For instance, if $(G,X)$ is of $D_n^{\dbR}$ type with $n\in 2+2\dbN$, then the half spin representations of the $D_n$ Lie type are self dual (cf. [Bou2, p. 210]) and thus [De2, Subsubsects. 2.3.10 to 2.3.13] does not apply to get directly that in Subsections 4.7 and 4.8 we can choose $(G_1,X_1)=(G_4^{\ad},X_4^{\ad})$ to be a simple, adjoint Shimura pair.

In addition, the explicit character of Subsections  4.1 to 4.8 and the new features of Proposition 4.1 and Remarks 4.8.2 (b) and (c) offer simplifications in many practical situations (like in Definition 6.1 (b) and in the proofs of Theorems 7.1 and 7.2 below).

\bigskip\smallskip
\noindent
{\boldsectionfont 6. Non-special $A_n$ types}
\bigskip

In this Section we define and study simple, adjoint Shimura pairs of non-special $A_n$ type. See Subsection 6.2 for different numerical aspects and examples that pertain to them. The notations will be independent from the ones of Subsections 1.1 and 1.3. If $(G_0,X_0)$ is a Shimura pair and $x\in X_0$, let the cocharacter $\mu_x:\dbG_{m,\dbC}\to G_{0,\dbC}$ be defined as in Subsection 2.2. We use $+$ to denote the Lie type of products of split, absolutely simple, adjoint groups.

\bigskip\noindent
{\bf 6.1. Definitions.} {\bf (a)} Let $(H_1,X_1)$ be a Shimura pair. A torus $T_1$ of $H_1$ is called {\it $Sh$-good} for $(H_1,X_1)$, if there exists a set $\grC_1$ of cocharacters $\dbG_{m,\dbC}\to H_{1,\dbC}$ that are $H_1(\dbC)$-conjugate to the cocharacters $\mu_x:\dbG_{m,\dbC}\to H_{1,\dbC}$ ($x\in X_1$) and such that $T_1$ is the smallest subgroup of $H_1$ with the property that the elements of $\grC_1$ factor through $T_{1,\dbC}$.

\smallskip
{\bf (b)} Let $(H_0,X_0)$ be a simple, adjoint Shimura pair of either $A_1$ type or of $A_n$ type with $n\Ge 2$ and with involution (resp. of $A_n$ type with $n\Ge 2$ and without involution). Let $f_1\colon (H_1,X_1)\hookrightarrow (\text{\bf GSp}(W_1,\psi_1),S_1)$ be an injective map of Shimura pairs such that the following four conditions hold:

\medskip
$\bullet$ we have a natural identity $(H_1^{\ad},X_1^{\ad})=(H_0,X_0)$;

$\bullet$ property 5.10 (a) holds (resp. properties 5.10 (a) and (b) hold) for $f_2=f_1$;

$\bullet$ the intersection $\tilde H_1$ of $\text{\bf GSp}(W_1,\psi_1)$ with the double centralizer of $H_1$ in $\text{\bf GL}_W$ is connected and the monomorphism $H_{1,\overline{\dbQ}}^{\der}\hookrightarrow \tilde H_{1,\overline{\dbQ}}^{\der}$ is an isomorphism (resp. is isomorphic to the diagonal embedding of $H_{1,\overline{\dbQ}}^{\der}$ in $H_{1,\overline{\dbQ}}^{\der}\times_{\overline{\dbQ}} H_{1,\overline{\dbQ}}^{\der}$);

$\bullet$ $H_1$ is a Mumford--Tate group (i.e., $H_1$ is the smallest subgroup of $\text{\bf GSp}(W_1,\psi_1)$ such that all monomorphisms $\dbS\hookrightarrow H_{1,\dbR}$ that are elements of $X_1$, factor through $H_{1,\dbR}$).

\medskip
Let $q\in\dbN$ be a prime. The image of each element $x\in X_1$ in $\tilde H_{1,\dbR}^{\ad}$ is non-trivial. Thus the group $\tilde H_1^{\ad}$ has no simple factor which over $\dbR$ is compact. From this we get that the $\tilde H_1(\dbR)$-conjugacy class $\tilde X_1$ of monomorphisms $\dbS\hookrightarrow \tilde H_{1,\dbR}$ defined by elements of $X_1$, has the property that we get injective maps $(H_1,X_1)\hookrightarrow (\tilde H_1,\tilde X_1)\operatornamewithlimits{\hookrightarrow}\limits^{\tilde f_1} (\text{\bf GSp}(W_1,\psi_1),S_1)$ of Shimura pairs. Let $\tilde f_1^\prime:(\tilde H_1^\prime,\tilde X_1^\prime)\hookrightarrow (\text{\bf GSp}(W_1,\psi_1),S_1)$ be an injective map of Shimura pairs defined by the twist of $\tilde f_1$ via a class in $H^1(\dbQ,\tilde H_1)$ that is locally trivial. We identify $\tilde X_1^\prime=\tilde X_1$ and $\tilde H_{1,\dbQ_q}=\tilde H_{1,\dbQ_q}^\prime$.

By an {\it $MT$ pair} for the triple $(f_1,\tilde f_1^\prime,q)$ we mean a pair $(T_1,G_1)$, where $T_1$ is a torus of $\tilde H_1^\prime$ that is $Sh$-good for $(\tilde H_1^\prime,\tilde X_1^\prime)$ and where $G_1$ is a reductive subgroup of $H_{1,\dbQ_q}$ that has the following two properties:

\medskip
{\bf (i)} the group $G_1$ has a maximal torus that is $\tilde H_1^\prime(\dbQ_q)$-conjugate to $T_{1,\dbQ_q}$;

\smallskip
{\bf (ii)} we have an identity $\End(W_1\otimes_{\dbQ} \dbQ_q)^{G_1(\dbQ_q)}=\End(W_1\otimes_{\dbQ} \dbQ_q)^{H_{1,\dbQ_q}(\dbQ_q)}$.

\medskip
{\bf (c)} Let $(H_0,X_0)$ be as in (b). We say that $(H_0,X_0)$ is of {\it non-special $A_n$ type}, if there exists an injective map $f_1\colon (H_1,X_1)\hookrightarrow (\text{\bf GSp}(W_1,\psi_1),S_1)$ of Shimura pairs as in (b) and a prime $q\in\dbN$ such that for each twist $\tilde f_1^\prime$ of $\tilde f_1$ obtained as in (b) and for every $MT$ pair $(T_1,G_1)$ for $(f_1,\tilde f_1^\prime,q)$, the group $G_{1,\overline{\dbQ_q}}^{\ad}$ has at least one simple factor of $A_n$ Lie type.

\bigskip\noindent
{\bf 6.2. Numerical tests.}
Let $(H_0,X_0)$ be a simple, adjoint Shimura pair of $A_n$ type. In this Subsection we list some numerical properties which are implied by the assumption that $(H_0,X_0)$ is or is not of non-special $A_n$ type and which rely on the classification (see Subsection 2.1) of minuscule weights.

We write
$$H_{0,\dbR}=\prod_{j\in J_1} \text{\bf SU}(a_j,b_j)_{\dbR}^{\ad}\times\prod_{j\in J_0} \text{\bf SU}(0,n+1)_{\dbR}^{\ad},$$
where for $j\in J_1$ we have $a_j,b_j\in\dbN$, $a_j\Le b_j$, and $a_j+b_j=n+1$. Here $J_0$ and $J_1$ are finite sets that keep track of the simple, non-compact and the simple, compact (respectively) factors of $H_{0,\dbR}$. As $H_{0,\dbR}$ has simple, non-compact factors (see Subsection 2.2), the set $J_1$ is non-empty. For $j\in J_1$ let
$$c_j:={a_j\over b_j}.$$
Let
$$\grC:=\{c_j|j\in J_1\}\subseteq\{{1\over n},{2\over {n-1}},\ldots,{[{{n+1}\over 2}]\over {n+1-[{{n+1}\over 2}]}}\}\;\;\text{and}\;\;\grM:=\{(a_j,b_j),(b_j,a_j)|j\in J_1\}.$$
The set $\grC$ has at least one element and at most $[{{n+1}\over 2}]$ elements.
Let
$$c:=\min(\grC)\,\,\text{and}\,\,d:=\max(\grC).$$
If $1\in \grC$, then $n$ is odd. Let $\grS:=\{(r,s)\in\dbN\times\dbN|2s-1\Le r\}$. For $(r,s)\in \grS$ let
$$c_{(r,s)}:={{\binom {r}{s-1}}\over {\binom {r}{s}}}\,\, \text{and}\,\, e_{(r,s)}:=\binom {r+1}{s}.$$
\noindent
{\bf 6.2.1. General properties of $MT$ pairs.} Let $q$, $f_1$, $\tilde f_1$, and $\tilde f_1^\prime$ be as in Definition 6.1 (b). Let $(T_1,G_1)$ be an $MT$ pair for $(f_1,\tilde f_1^\prime,q)$. Let $\scrW_0$ be an irreducible $\tilde H_{1,\overline{\dbQ_q}}$-submodule of $W_1\otimes\overline{\dbQ_q}$; it is also an irreducible $H_{1,\overline{\dbQ_q}}$-module. Let $\scrG_0$ and $\scrH_0$ be the images of $G_{1,\overline{\dbQ_q}}$ and $H_{1,\overline{\dbQ_q}}$ (respectively) in $\text{\bf GL}_{\scrW_0}$. We get a faithful representation $\rho_0:\scrG_0\hookrightarrow \text{\bf GL}_{\scrW_0}$. We list few basic properties of $\rho_0$ and $\scrH_0$.

\medskip
{\bf (a)} As the property 5.10 (a) holds for $f_2=f_1$, we have $\dim_{\overline{\dbQ_q}}(\scrW_0)=n+1$.

\smallskip
{\bf (b)} The identity $\End(W_1\otimes_{\dbQ} \dbQ_q)^{G_1(\dbQ_q)}=\End(W_1\otimes_{\dbQ} \dbQ_q)^{H_{1,\dbQ_q}(\dbQ_q)}$ implies that the representation $\rho_0$ is irreducible.

\smallskip
{\bf (c)} As the torus $T_1$ is $Sh$-good for $(\tilde H_1^\prime,\tilde X_1^\prime)$, the pair $(\scrG_0,\rho_0)$ is an irreducible weak Mumford--Tate pair of weight $\{0,1\}$. More precisely, the group $\scrG_0$ is generated by images in $\text{\bf GL}_{\scrW_0}$ of cocharacters $\mu_0$ of $\tilde H_{1,\overline{\dbQ_q}}=\tilde H^\prime_{1,\overline{\dbQ_q}}$ whose extensions to $\dbC$ via an (any) embedding $\overline{\dbQ_q}\hookrightarrow\dbC$, belong to the set $\grC(\tilde H_1^\prime,\tilde X_1^\prime)$ introduced in Subsection 2.2. Thus all simple factors of $G_{1,\overline{\dbQ_q}}^{\ad}$ are of classical Lie type (cf. [Pi, Table 4.2]) and the pair $(\scrG_0,\rho_0)$ is non-trivial (cf. (a)).

\smallskip
{\bf (d)} Each cocharacter $\mu_0$ as in (c) acts on $\scrW_0$ via the identical and the trivial character of $\dbG_{m,\overline{\dbQ_q}}$ and the pair of multiplicities $(m_{\text{id}},m_{\text{triv}})$ of such characters either belongs to $\grM$ or is such that $\{m_{\text{id}},m_{\text{triv}}\}=\{0,n+1\}$, cf. (a) and (c). Therefore, if $d<1$, then we have $m_{\text{id}}\neq m_{\text{triv}}$ and thus from the list of multiplicities in [Pi, Table 4.2] we get that all simple factors of $G_{1,\overline{\dbQ_q}}^{\ad}$ are of some $A_m$ Lie type, where $m\in\{1,\ldots,n\}$.

\smallskip
{\bf (e)} As $T_1$ is defined over $\dbQ$, we can choose the cocharacter $\mu_0$ of $\tilde H_{1,\overline{\dbQ_q}}=\tilde H^\prime_{1,\overline{\dbQ_q}}$ such that either $(m_{\text{id}},m_{\text{triv}})$ or $(m_{\text{triv}},m_{\text{id}})$ is an arbitrary a priori given element of $\grM$.

\smallskip
{\bf (f)} We assume that the group $\scrG^{\ad}_0$ is simple and not of $A_n$ Lie type. As $\scrG^{\ad}_0$ is simple, the pair $(\scrG_0,\rho_0)$ is a non-trivial irreducible strong Mumford--Tate pair of weight $\{0,1\}$. From this and [Pi, Table 4.2] we get that the set $\{(m_{\text{id}},m_{\text{triv}}),(m_{\text{triv}},m_{\text{id}})\}$ is uniquely determined by $\rho_0$ (i.e., it does not depend on $\mu_0$). This implies that $\grC$ has only one element, cf. (e).

\smallskip
{\bf (g)} The pair $(\scrG_0,\rho_0)$ is in the natural way a product of non-trivial irreducible strong Mumford--Tate pairs of weight $\{0,1\}$ (cf. [Pi, Sect. 4, pp. 210--211]) which automatically involve groups whose adjoints are simple (cf. [Pi, Prop. 4.4]). If there exists $j\in J_1$ such that $g.c.d.(a_j,b_j)=1$, then from (e) applied with $(m_{\text{id}},m_{\text{triv}})\in\{(a_j,b_j),(b_j,a_j)\}$ we get that this natural product does not have two or more irreducible factors and therefore the group $\scrG^{\ad}_0$ is simple.

\medskip\noindent
{\bf 6.2.2. Lemma.} {\it  We assume that $d<1$ and that $(H_0,X_0)$ is not of non-special $A_n$ type. Then for all elements $j\in J_1$ there exists a number $t_j\in\dbN$ and there exist numbers $r_{j,m}\in\dbN$, with $m\in\{1,\ldots,t_j\}$, such that the following three properties hold:

\medskip
{\bf (a)} the sum $\sum_{m=1}^{t_j} r_{j,m}$ does not depend on $j\in J_1$;

\smallskip
{\bf (b)} we have $\sum_{m=1}^{t_j} r_{j,m}<n$;

\smallskip
{\bf (c)} there exists a map $f(j)\colon\{1,\ldots,t_j\}\to\{0,1\}$ that has the following four properties

\medskip\noindent
{\bf (c.i)} we have $n+1=\prod_{m\in\{1,\ldots,t_j\}, f(j)(m)=0} (r_{j,m}+1)\times\prod_{m\in\{1,\ldots,t_j\}, f(j)(m)=1} e_{(r_{j,m},s_{j,m})}$;

\smallskip\noindent
{\bf (c.ii)} if $f(j)(m)=0$, then $r_{j,m}+1=a_{j,m}+b_{j,m}$, with $a_{j,m},b_{j,m}\in\dbN$ such that $c_{j,m}:={a_{j,m}\over b_{j,m}}\in \grC$;

\smallskip\noindent
{\bf (c.iii)} if $f(j)(m)=1$, then there exists a pair $(r_{j,m},s_{j,m})\in \grS$ such that $2\Le s_{j,m}\Le r_{j,m}-2$ and $c_{j,m}:=c_{(r_{j,m},s_{j,m})}\in \grC$;

\smallskip\noindent
{\bf (c.iv)}  the map $q_j\colon \{1,\ldots,t_j\}\to \grC$ that takes $m$ to $c_{j,m}$ has the property that we have an inclusion
$\grC\subseteq\text{Im}(q_j)\bigcup_{m\in\{1,\ldots,t_j\},\; f(j)(m)=0}\{{1\over {r_{j,m}}},{2\over {r_{j,m}-1}},\ldots,{[{{r_{j,m}+1}\over 2}]\over {r_{j,m}+1-[{{r_{j,m}+1}\over 2}]}}\}$.}\

\medskip
\proof
See Subsection 2.2 for the totally real number field $F(H_0,X_0)$. We fix an identification of $\overline{\dbQ}$ with the algebraic closure of $\dbQ$ in $\overline{\dbQ_q}$. It allows us to identify naturally the set $\Hom(F(H_0,X_0),\dbR)=\Hom(F(H_0,X_0),\overline{\dbQ})$ with the set $\Hom(F(H_0,X_0),\overline{\dbQ_q})$. This implies that we can identify the set $J$ of simple factors of $H_{0,\dbR}$ with the set of simple factors of $H_{0,\overline{\dbQ_q}}$. To prove the Lemma we can use the notations of Definition 6.1 (b) and we can assume that the group $G^{\ad}_{1,\overline{\dbQ_q}}$ does not have simple factors of $A_n$ Lie type.

For $j\in J_1$ let $t_j$ be the number of simple factors of the image $G_0(j)$ of $G_{1,\overline{\dbQ_q}}$ in the simple factor $H_0(j)$ of $H_{1,\overline{\dbQ_q}}^{\ad}\tilde\to H_{0,\overline{\dbQ_q}}$ that corresponds to $j$. We will refer to these simple factors as the first, second, $\ldots$, and the $t_j$-th simple factors of $G_0(j)$. For $m\in\{1,\ldots,t_j\}$ let $r_{j,m}$ be the rank of the $m$-th simple factor of $G_0(j)$. Property (a) is implied by the fact that the images of $T_{1,\overline{\dbQ_q}}$ (and thus also of $G_{1,\overline{\dbQ_q}}$) in simple factors of $\tilde H_{1,\overline{\dbQ_q}}^{\ad}$ have equal ranks. We check that (b) holds. As the group $G^{\ad}_{1,\overline{\dbQ_q}}$ does not have simple factors of $A_n$ Lie type, we have $G_0(j)\neq H_0(j)$. As $\scrW_0$ is a simple $\Lie(G_0(j))$-module and a simple $\Lie(H_0(j))$-module, from Lemma 1.1.7 we get that the rank of $\Lie(G_0(j))$ is less than the rank of $\Lie(H_0(j))$ i.e., we have $\sum_{m=1}^{t_j} r_{j,m}<n$. Thus (b) holds

We define the map $f(j)$ in such a way that we have $f(j)(m)=0$ if and only if the $m$-th simple factor  of $G_0(j)$ is embedded in $H_0(j)$ via an SD-representation.

As $d<1$, for each $m\in\{1,\ldots,r_j\}$ there exists a number $r_{j,m}\in\dbN$ such that the $m$-th simple factor of $G_0(j)$ is of $A_{r_{j,m}}$ Lie type (cf. property 6.2.1 (d)). For $m\in\{1,\ldots,t_j\}$ such that $f(j)(m)=1$, let $s_{j,m}\in \dbN$ be such that $(r_{j,m},s_{j,m})\in\grS$ and the $m$-th simple factor  of $G_0(j)$ is embedded in $H_0(j)$ via an irreducible representation associated to the minuscule weight $\varpi_{s_{j,m}}$ of the $A_{r_{j,m}}$ Lie type (cf. Corollary 1.1.9); as this representation is not an SD representation we have $2\Le s_{j,m}\Le r_{j,m}-2$. For $m\in\{1,\ldots,t_j\}$ such that $f(j)(m)=0$, let $a_{j,m},b_{j,m}\in\dbN$ be numbers for which we have $a_{j,m}+b_{j,m}=r_{j,m}+1$ and $a_{j,m}\Le {{r_{j,m}+1}\over 2}$ and for which there exists a cocharacter $\mu_0$ as in the property 6.2.1 (c) whose image in the $m$-th simple factor of $G_0(j)$ has a centralizer in the $m$-th simple factor of $G_0(j)$ whose derived group is of $A_{a_{j,m}-1}+A_{b_{j,m}-1}$ Lie type (i.e., this image is defined by one of the cocharacters $\varpi_{a_{j,m}}^{\vee}$ or $\varpi_{b_{j,m}}^{\vee}$ of the $A_{r_{j,m}}$ Lie type); we emphasize that here we do not have a unique choice for the pair $(a_{j,m},b_{j,m})$ (and this explains the $\cup$ sign in (c.iv)).

We check that (c) holds. Condition (c.i) is implied by the property 6.2.1 (a), cf. the fact that for $(r,s)\in\grS$ the irreducible representation associated to the minuscule weight $\varpi_s$ of the $A_r$ Lie type has dimension $e_{r,s}$. The fact that we have $c_{j,m}\in\grC$ for all $m\in\{1,\ldots,t_j\}$ is implied by the first part of the property 6.2.1 (d), cf. the list of multiplicities in [Pi, Table 4.2]. Thus (c.ii) and (c.iii) hold. Condition (c.iv) is implied by the property 6.2.1 (e). Thus (c) also holds.\endproof

\medskip\noindent
{\bf 6.2.3. Remarks.}
{\bf (a)} We assume that $c<d=1$. If $(H_0,X_0)$ is not of non-special $A_n$ type, then one can list numerical conditions that are similar to the ones of Lemma 6.2.2. We do not stop to list them here. We only mention that (as in the property 6.2.2 (c)) we get maps $q_j$'s into $\grC$ which are ``almost surjective'' (cf. property 6.2.1 (e)) and that (as $d=1$) the last sentence of the property 6.2.1 (d) does not a priori hold.

\smallskip
{\bf (b)} If $n$ is odd, then we can always identify $\tilde f_1^\prime=\tilde f_1$ (cf. property 4.1.1 (a), Subsubsection 4.1.2, and property 5.7 (c)).

\medskip\noindent
{\bf 6.2.4. Examples.} We present some situations to which (the ideas of) Subsections 6.2.1 and 6.2.2 apply.

\smallskip
{\bf Example 1.} The Shimura pair $(H_0,X_0)$ is of non-special $A_n$ type if there exists an element $j\in J_1$ such that $g.c.d.(a_j,b_j)=1$ and $(a_j,b_j)\notin\{(\binom {r}{s-1},\binom {r}{s})|(r,s)\in \grS,\,\,s\Ge 2\}$. In particular, $(H_0,X_0)$ is of non-special $A_n$ type if there exists an element $j\in J_1$ such that $a_j=1$ (i.e., if ${1\over n}\in \grC$).

\smallskip
{\bf Example 2.} We assume that $c<d$ and that there exists $j\in J_1$ such that $g.c.d.(a_j,b_j)=1$. Then $(H_0,X_0)$ is of non-special $A_n$ type.

\smallskip
We argue Example 1 (resp. 2). As $g.c.d.(a_j,b_j)=1$, the group $\scrG_0^{\ad}$ is simple (cf. property 6.2.1 (g)). As $g.c.d.(a_j,b_j)=1$ and as the group $\scrG_0^{\ad}$ is simple, from the list of multiplicities in [Pi, Table 4.2] we get that $\scrG_0^{\ad}$ is of $A_m$ Lie type for some $m\in\{1,\ldots,n\}$. But as $(a_j,b_j)\notin\{(\binom {r}{s-1},\binom {r}{s})|(r,s)\in \grS,\,\,s\Ge 2\}$ (resp. as $\grC$ has at least two elements), from the mentioned list (resp. from the property 6.2.1 (f)) we get that $m=n$. Thus $\scrG_0^{\ad}$ has $A_n$ Lie type and therefore it is $\scrH_0^{\ad}$. We get that $G_{1,\overline{\dbQ_q}}^{\ad}$ has simple factors of $A_n$ Lie type. Therefore $(H_0,X_0)$ is of non-special $A_n$ type.

\smallskip
{\bf Example 3.} We assume that $c<d<1$ and that there exists an element $j\in J_1$ such that $g.c.d.(a_j,b_j)$ does not have factors that are of the form $\tilde a_j+\tilde b_j$, with $\tilde a_j,\tilde b_j\in\dbN$ satisfying ${\tilde a_j\over\tilde b_j}\in \grC$, or of the form $e_{(\tilde r_j,\tilde s_j)}$, with $(\tilde r_j,\tilde s_j)\in \grS$ satisfying $c_{(\tilde r_j,\tilde s_j)}\in \grC$. Then $(H_0,X_0)$ is of non-special $A_n$ type.

To check this, we show that the assumption that $\scrG_0^{\ad}$ is not of $A_n$ Lie type leads to a contradiction. The group $\scrG_0^{\ad}$ is not simple (cf. property 6.2.1 (f)) and thus $(\scrG_0,\rho_0)$ is a tensor product of at least two non-trivial irreducible strong Mumford--Tate pairs of weight $\{0,1\}$ (cf. property 6.2.1 (g)). From this and the property 6.2.1 (d), we get that $\rho_0$ is a tensor product of at least two irreducible representations associated to minuscule weights of some $A_m$ Lie types. The dimensions of these irreducible representations are equal to one of the numbers $\tilde a_j+\tilde b_j$ or $e_{(\tilde r_j,\tilde s_j)}$ of the previous paragraph, cf. the list of multiplicities in [Pi, Table 4.2]. As at least one of these dimensions must divide $g.c.d.(a_j,b_j)$ (cf. property 6.2.1 (e)), we reached a contradiction. Thus $(H_0,X_0)$ is of non-special $A_n$ type.

\smallskip
{\bf Example 4.} We assume that $c=d=1$. If $n=3$ or if $4$ does not divide $n+1$, then $(H_0,X_0)$ is of non-special $A_n$ type. To check this, we first remark that $c=d=1$ implies that $(H_0,X_0)$ is without involution (cf. Subsubsection 2.2.2). Thus we can appeal to Proposition 5.10 (b). With the notations of the beginning paragraph of Subsubsection 6.2.1, let $\scrW_0^d$ be a simple $H_{1,\overline{\dbQ_q}}$-submodule of $W_1\otimes_{\dbQ} \overline{\dbQ_q}$ which as an $H^{\der}_{1,\overline{\dbQ_q}}$-module is dual to $\scrW_0$ and on which $Z^0(H_{1,\overline{\dbQ_q}})$ acts in the same way as on $\scrW_0$. As ${{n+1}\over 4}\not\in\dbN\setminus\{1\}$ and as $c=d=1$, from the property 6.2.1 (c) and [Pi, Table 4.2] we get that the group $\scrG_0^{\ad}$ is of $C_{{n+1}\over 2}$, $D_{{n+1}\over 2}$, $A_1+A_1$, or some $A_r$ Lie type. It is not of $C_{{n+1}\over 2}$, $D_{{n+1}\over 2}$, or $A_1+A_1$ Lie type, as otherwise the $G_{1,\overline{\dbQ_q}}$-modules $\scrW_0$ and $\scrW_0^d$ would be isomorphic and thus we would have $\End(W_1\otimes_{\dbQ} \dbQ_q)^{G_1(\dbQ_q)}\neq\End(W_1\otimes_{\dbQ} \dbQ_q)^{H_{1,\dbQ_q}(\dbQ_q)}$. Thus $\scrG^{\ad}$ is of some $A_r$ Lie type. Let $\varpi_s$ be the highest weight of the representation of $\Lie(\scrG_0^{\ad})$ on $\scrW_0$, cf. Corollary 1.1.9. If $r<n$, then as $c=d=1$ we must have $c_{(r,s)}=1$ and $n+1=e_{(r,s)}$. But we have $c_{(r,s)}=1$ if and only if $r=2s-1$. If $r=2s-1<n$, then the irreducible representation of $\Lie(\scrG_0^{\ad})$ of the highest weight $\varpi_s$ is self dual (cf. [Bou2, pp. 188--189]) and this implies that the $G_{1,\overline{\dbQ_q}}$-modules $\scrW_0$ and $\scrW_0^d$ are isomorphic. Thus, in order not to get as above a contradiction with the identity $\End(W_1\otimes_{\dbQ} \dbQ_q)^{G_1(\dbQ_q)}=\End(W_1\otimes_{\dbQ} \dbQ_q)^{H_{1,\dbQ_q}(\dbQ_q)}$, we must have $r=n$. Therefore $\scrG^{\ad}$ is of $A_n$ Lie type. Thus $(H_0,X_0)$ is of non-special $A_n$ type.

\smallskip
{\bf Example 5.} We assume that $c<d=1$ and $(H_0,X_0)$ is not of non-special $A_n$ type. As $d=1$, $n+1$ is even. The pair $(\scrG_0,\rho_0)$ is a product of at least two non-trivial irreducible strong Mumford--Tate pairs of weight $\{0,1\}$ (cf. properties 6.2.1 (f) and (g)) such that at least one of them is even dimensional (as $n+1$ is even). If $(\scrG_0,\rho_0)$ is a product of even dimensional irreducible strong Mumford--Tate pairs of weight $\{0,1\}$ (resp. if $(\scrG_0,\rho_0)$ is a product of non-trivial irreducible strong Mumford--Tate pairs of weight $\{0,1\}$ such that at least one is odd dimensional and at least one is even dimensional), then  from the property 6.2.1 (e) (resp. from the generation part of the property 6.2.1 (c) and from the property 6.2.1 (d)), we get that for all $j\in J_1$ such that $c_j<1$, $2$ divides both $a_j$ and $b_j$ (resp. we get that there exists $j\in J_1$ such that $c_j<1$ and $2$ divides both $a_j$ and $b_j$). Thus if $n+1$ is a power of $2$, then for all $j\in J_1$ such that $c_j<1$, $2$ divides both $a_j$ and $b_j$.

\smallskip
{\bf Example 6.} If $n+1$ is either $4$ or a prime, then $(H_0,X_0)$ is of non-special $A_n$ type (regardless of what the set $\grC$ is). This is an elementary exercise which follows from Examples 1 and 4.

\smallskip
We do not stop here to restate Lemma 6.2.2 in terms of $n$, similar to the statements of Examples 4 and 6. Example 1.3.5 (i) (resp. 1.3.5 (iii)) is a particular case of Example 2 (resp. 5); moreover, Example 1.3.5 (ii) (resp. 1.3.5 (iv)) is only Example 1 (resp. 6).

\bigskip\smallskip
\noindent
{\boldsectionfont 7. The proof of the Main Theorem}
\bigskip

Theorems 7.1 and 7.2 are the very essence of the proofs of Theorems 1.3.4 and 1.3.6 which are completed in Subsection 7.3. Corollary 7.4 is an application to the $D_4^{\dbH}$ type of the proofs of Theorems 7.1 and 7.2. Subsection 7.5 lists the ``simplest" cases not settled by Theorem 1.3.4. Theorem 7.6 proves a stronger variant of Theorem 1.3.7 for direct sums of abelian motives in the sense of [DM]. We end with remarks in Subsection 7.7.

We use $+$ to denote the Lie type of products of split, absolutely simple, adjoint groups. Until the end we use the notations of Subsections 1.1 and of the paragraphs before Theorems 1.1.5 and 1.3.2. Let $W_A$, $\lambda_A$, $\psi_A$, and $f_A$ be as in Subsection 1.5.

\bigskip\noindent
{\bf 7.1. Theorem.} {\it We assume that for each element $j\in J$ the adjoint Shimura pair $(H_j,X_j)$ is of either non-inner $D_4^{\dbH}$ type or a type listed in Theorem 1.3.4. Then for each element $j\in J$, the group $G_{\overline{\dbQ_p}}$ surjects onto a simple factor of $H_{j,\overline{\dbQ_p}}$.}

\medskip
\proof
To prove the Theorem we can assume that $\lambda_A$ is a principal polarization of $A$ (cf. property 3.1 (a) and [Mu2, Ch. IV, Sect. 23, Cor. 1]) and that the adjoint group $H_A^{\ad}$ is non-trivial (cf. property 3.4.1 (a)). Thus the set $J$ that indexes the simple factors of $(H_A^{\ad},X_A^{\ad})$ is non-empty. To prove the Theorem we can also assume that the set $J$ has precisely one element, cf. Proposition 5.4.1 and Subsubsection 5.4.3. Thus $(H_j,X_j)=(H_A^{\ad},X_A^{\ad})$. We take $v$ such that it is unramified over $l$.

Let the triple $(\scrW_0,\scrG_0,\scrH_0)$ be as before Subsubsection 1.4.1. The groups $\scrG_0^{\ad}$ and $\scrH_0^{\ad}$ are isomorphic with simple factors of $G_{\overline{\dbQ_p}}^{\ad}$ and $H_{A,\overline{\dbQ_p}}^{\ad}$ (respectively). Thus if $\scrG_0^{\ad}=\scrH_0^{\ad}$, then $G_{\overline{\dbQ_p}}^{\ad}$ surjects onto a simple factor of $H_{A,\overline{\dbQ_p}}^{\ad}$. We are left to check that:

\medskip
{\bf (*)} we can choose the simple $H_{A,\overline{\dbQ_p}}$-module $\scrW_0$ such that we have $\scrG_0^{\ad}=\scrH_0^{\ad}$.

\medskip
To check (*) we will consider one by one the cases (a) to (e) of Theorem 1.3.4, with the non-inner $D_4^{\dbH}$ being added to the case (d'). Only for the case (d') we will have to make a good choice of $\scrW_0$.

\smallskip
{\bf (a)} We assume that $(H_A^{\ad},X_A^{\ad})$ is of non-special $A_n$ type. To check the property (*) we can assume that $E^{\text{conn}}=E$. To check the property (*) we can also assume that the injective map $f_A:(H_A,X_A)\hookrightarrow (\text{\bf GSp}(W_A,\psi_A),S_A)$ of Shimura pairs factors through an injective map $\tilde f_A:(\tilde H_A,\tilde X_A)\hookrightarrow (\text{\bf GSp}(W_A,\psi_A),S_A)$ of Shimura pairs which is a PEL type embedding such that the monomorphism $H_{A,\overline{\dbQ}}^{\der}\hookrightarrow \tilde H_{A,\overline{\dbQ}}^{\der}$ is (cf. Proposition 5.4.1, Subsubsection 5.4.2, Proposition 4.1, and Subsubsection 4.1.2): (i) an isomorphism if either $n=1$ or $(H_A^{\ad},X_A^{\ad})$ is with involution, and (ii) is a diagonal embedding if $n\Ge 2$ and $(H_A^{\ad},X_A^{\ad})$ is without involution. Based on Definition 6.1 (c), we can also assume that there exists a prime $q\in\dbN$ such that for each twist $\tilde f_{A}^\prime$ of $\tilde f_{A}$ via a locally trivial class in $H^1(\dbQ,\tilde H_A)$ and for every $MT$ pair $(T_1,G_1)$ for $(f_{A},\tilde f_{A}^\prime,q)$, the group $G_{1,\overline{\dbQ_q}}^{\ad}$ has simple factors of $A_n$ Lie type.

We choose $v$ such that we also have $k(v)=\dbF_l$, $l\notin \{p,q\}$, $l$ is greater than the number $N(A)$ of Theorem 1.3.2, and (cf. Theorem 1.1.5) $T_v$ has the same rank as $G_{\dbQ_p}$. The torus $T_{v,\dbQ_p}$ of $G_{\dbQ_p}$ is $\tilde H_A(\dbQ_p)$-conjugate to the extension to $\dbQ_p$ of a torus $T_4$ of $\tilde H_A$ that is naturally identifiable with $T_v$ (cf. Subsubsection 5.6.3 applied with $G_2=G_4=\tilde H_A$). Strictly speaking, we ought to substitute in the last sentence $\tilde H_A$ by a suitable form of it that is locally isomorphic. But from the point of view of what follows, such a substitution is irrelevant. Thus we will not complicate notations by introducing explicitly a form of $\tilde H_A$. A suitable $\tilde H_A(\overline{\dbQ_l})$-conjugate of $T_{4,\overline{\dbQ_l}}$ is a subtorus of $G_{\overline{\dbQ_l}}$, cf. Subsubsection 5.6.4 (applied in the context of $A$ instead of $B$). Thus from Corollary 3.4.2 we get that $T_{4,\overline{\dbQ_l}}$ is generated by images of cocharacters in the set $\grC_l(\tilde H_A,\tilde X_A)$ of cocharacters of $\tilde H_{A,\overline{\dbQ_l}}$ whose extensions to $\dbC$ via an embedding $i_{\overline{\dbQ_l}}:\overline{\dbQ_l}\hookrightarrow\dbC$ that extends the embedding $i_E$ of Subsection 1.1, belong to the set $\grC(\tilde H_A,\tilde X_A)$ defined in Subsection 2.2. This implies that the torus $T_1:=T_4$ is $Sh$-good for $(\tilde H_A,\tilde X_A)$. As the groups $T_1$ and $G_{\dbQ_p}$ have the same ranks, from Theorem 1.1.4 we get that $(T_1,G_{\dbQ_p})$ is an $MT$ pair for $(f_A,\tilde f_A,p)$.

If $p$ is $q$, then $G_{\overline{\dbQ_p}}^{\ad}$ has simple factors of $A_n$ Lie type and thus the image of $T_{1,\overline{\dbQ}}$ in each simply factor of $\tilde H_{A,\overline{\dbQ}}^{\ad}$ has rank $n$. This last property implies that, regardless of the fact that $p$ is not or is $q$, the images of $G_{\overline{\dbQ_p}}$ in the simply factors of $H^{\ad}_{A,\overline{\dbQ_p}}$ (equivalently the images of $T_{1,\overline{\dbQ_p}}$ in the simply factors of $\tilde H^{\ad}_{A,\overline{\dbQ_p}}$) have rank $n$. From this and Lemma 1.1.7 applied in the context of the Lie monomorphisms $\Lie(\scrG_0^{\ad})\hookrightarrow\Lie(\scrH_0^{\ad})\hookrightarrow\End(\scrW_0)$, we get that $\Lie(\scrG_0^{\ad})=\Lie(\scrH_0^{\ad})$. Thus $\scrG_0^{\ad}=\scrH_0^{\ad}$ (i.e., the property (*) holds).

\smallskip
Based on Proposition 5.4.1, for the below cases (b) to (e) we can assume that $f_A$ factors through an injective map $f_2:(G_2,X_2)\to (\pmb{\text{GSp}}(W_1,\psi_1),S_1)=(\pmb{\text{GSp}}(W_A,\psi_A),S_A)$ of Shimura pairs for which the properties 5.10 (a) and (b) hold and for which the monomorphism $H_A\hookrightarrow G_2$ induces an identity $H_A^{\der}=G_2^{\der}$ at the level of derived groups.

As the representation $\rho_0$ is irreducible, from Theorem 1.1.8 we get that we have a short exact sequence $1\to Z(\text{\bf GL}_{\scrW_0})\to\scrG_0\to\scrG_0^{\ad}\to 1$.

\smallskip
{\bf (b)} If $(H_A^{\ad},X_A^{\ad})$ is of $B_n$ type, then $\scrH_0$ is a $\text{\bf GSpin}_{2n+1,\overline{\dbQ_p}}$ group and $\scrW_0$ is its spin representation. As $\rho_0$ is irreducible, from [Pi, Prop. 4.3] we get that $\scrG_0=\scrH_0$. Thus the property (*) holds.

\smallskip
{\bf (c)} If $(H_A^{\ad},X_A^{\ad})$ is of $C_n$ type, then $\scrH_0$ is a $\text{\bf GSp}_{2n,\overline{\dbQ_p}}$ group and $\scrW_0$ is of dimension $2n$. We assume that $n$ is odd. Thus $2n$ can not be written as a product of two or more even numbers. As the pair $(\scrG_0,\rho_0)$ is an irreducible weak Mumford--Tate pair of weight $\{0,1\}$ and as the representation $\rho_0$ is symplectic, the pair $(\scrG_0,\rho_0)$ is a product of irreducible strong Mumford--Tate pairs of weight $\{0,1\}$ that involve (orthogonal or symplectic and thus) even dimensional irreducible representations. From the last two sentences we get that $(\scrG_0,\rho_0)$ is a strong Mumford--Tate pair of weight $\{0,1\}$. From this and [Pi, Prop. 4.7 or Table 4.2] we get $\scrG_0=\scrH_0$. Thus the property (*) holds.

\smallskip
{\bf (d)} We assume that $(H_A^{\ad},X_A^{\ad})$ is of $D_n^{\dbH}$ type with $n$ odd and such that $2n\notin\{\binom {2^{m+1}}{2^m}|m\in\dbN\}$. This case is the orthogonal variant of (c). The group $\scrH_0$ is a $\text{\bf GSO}_{2n,\overline{\dbQ_p}}$ group and $\scrW_0$ is of dimension $2n$. As the representation $\rho_0$ is irreducible and orthogonal, it is a tensor product of (orthogonal or symplectic and thus of) even dimensional irreducible representations. Thus, as $n$ odd, the group $\scrG_0^{\ad}$ is simple. As the number $2n\in 2\dbN+2$ is not of the form $\binom {2^{m+1}}{2^m}$, with $m\in\dbN$, it is also not of the form $\binom {4m}{2m}$, with $m\in\dbN$. From this and the fact that $n$ is not of the form $2^m$ with $m\Ge 2$, we get that $(\scrG_0,\rho_0)=(\scrH_0,\rho_0)$ (cf. [Pi, Table 4.2]). Thus the property (*) holds.

\smallskip
{\bf (d')} We assume that $(H_A^{\ad},X_A^{\ad})$ is of non-inner $D_{2n}^{\dbH}$ type, with $n$ a prime. The group $\scrH_0$ is a $\text{\bf GSO}_{4n,\overline{\dbQ_p}}$ group and $\scrW_0$ is of dimension $4n$.

We have $[I(H_A^{\ad},X_A^{\ad}):F(H_A^{\ad},X_A^{\ad})]=2$ (as in the proof of Lemma 4.4 one argues that for $n=2$ the degree $[I(H_A^{\ad},X_A^{\ad}):F(H_A^{\ad},X_A^{\ad})]$ is not divisible by $3$). Based on the independence part of Corollary 1.3.3 and on Lemma 1.1.7, we can assume that $p$ is such that the number fields $I(H_A^{\ad},X_A^{\ad})$ and $E$ are unramified over $p$, that $H^{\ad}_{A,\dbQ_p}$ is unramified, and that (cf. the definition of the non-inner type and [La, p. 168]) there exists a prime $w$ of $F(H_A^{\ad},X_A^{\ad})$ which divides $p$ and which does not split in $I(H_A^{\ad},X_A^{\ad})$. We also take the prime $v$ such that $l=p$; thus $v$ is unramified over $p$. Let $F_{p,1}:=F(H_A^{\ad},X_A^{\ad})_w$. As $w$ does not split in $I(H_A^{\ad},X_A^{\ad})$ and as $H_{A,\dbQ_p}^{\ad}$ is unramified, there exists a non-split, absolutely simple factor $G_2$ of $H_{A,F_{p,1}}^{\ad}$. The group $G_2$ is unramified and splits over the unramified quadratic extension $F_{p,2}:=I(H_A^{\ad},X_A^{\ad})\otimes_{F(H_A^{\ad},X_A^{\ad})} F_{p,1}$ of $F_{p,1}$. We identify $\overline{F_{p,2}}=\overline{\dbQ_p}$. Let $G_3$ be the image of $G_{F_{p,1}}$ in $G_2$.

We choose $\scrW_0$ such that $\Lie(G_{2,\overline{\dbQ_p}})$ acts non-trivially on it. Thus we have natural identifications $G_{2,\overline{\dbQ_p}}=\scrH_0^{\ad}$ and $G_{3,\overline{\dbQ_p}}=\scrG_0^{\ad}$. Let $\mu_2$ be a non-trivial cocharacter of $G_{2,\overline{\dbQ_p}}=\scrH_0^{\ad}$ that is the image of a cocharacter of $\scrH_0$ which acts on $\scrW_0$ via the trivial and the identity characters of $\dbG_{m,\overline{\dbQ_p}}$. As $\scrW_0$ has dimension $4n$ and as $n\Ge 2$, the centralizer of $\mu_2$ in $G_{2,\overline{\dbQ_p}}$ has a derived subgroup of $A_{2n-1}$ Lie type. Based on Theorem 1.1.8, we can choose $\mu_2$ such that it factors through $G_{3,\overline{\dbQ_p}}$. As the centralizer of $\mu_2$ in $G_{2,\overline{\dbQ_p}}$ has a derived subgroup of $A_{2n-1}$ Lie type, the group $\Gal(F_{p,1})$ does not fix $[\mu_2]$ (this holds even if $n=2$, as the non-split group $G_2$ splits over $F_{p,2}=I(H_A^{\ad},X_A^{\ad})\otimes_{F(H_A^{\ad},X_A^{\ad})} F_{p,1}$).

The representation of $\Lie(\scrG_0^{\der})$ on $\scrW_0$ is a tensor product of $N$ irreducible representations that are associated to minuscule weights (cf. Corollary 1.9) and that are either symplectic or orthogonal, the number $N_S$ of symplectic ones being even. As $n$ is either odd or $2$, $4n$ is not divisible by $16$. From the last two sentences we get that $N\in\{1,2,3\}$. If $N=3$, then $n=2$ and therefore $N_S=N=3$; this contradicts the fact that $N_S$ is even. Thus $N\in\{1,2\}$.

We show that the assumption $G_3\neq G_2$ leads to a contradiction. As $n$ is a prime, we have $4n\notin\{\binom {4m}{2m}|m\in\dbN\}$. Thus from [Pi, Table 4.2] we get that either $N\neq 1$ or $n=2N=2$ and $\scrG_0^{\ad}$ is of $B_3$ Lie type. We first assume that $N=2$. The only way to factor $4n$ as a product of two even numbers is $2\times 2n$. Thus $2\Ge N_S\Ge 1$ and therefore $N_S=2$. Thus $\scrG_0^{\ad}$ is of $C_1+C_n$ Lie type (to be compared with (c)) and therefore each simple factor of $G_3$ is absolutely simple. Each group pair over $F_{p,1}$ that involves a product of absolutely simple, adjoint groups that over $\overline{\dbQ_p}$ are of some $C_m$ Lie types, has the property that its conjugacy class of cocharacters is fixed by $\Gal(F_{p,1})$. Thus $\Gal(F_{p,1})$ fixes $[\mu_2]$. Contradiction. If $n=2N=2$, then the same arguments lead to a contradiction: we only have to replace the $C_1+C_m$ Lie type with the $B_3$ Lie type. Therefore $G_3=G_2$ (i.e., $\scrG_0^{\ad}=\scrH_0^{\ad}$). Thus the property (*) holds.

\smallskip
{\bf (e)} We assume that $(H_A^{\ad},X_A^{\ad})$ is of $D_n^{\dbR}$ type. We show that the assumption $\scrG_0\neq \scrH_0$ leads to a contradiction. The faithful representation $\scrH_0^{\der}\hookrightarrow \text{\bf GL}_{\scrW_0}$ is a half spin representation. From this and [Pi, Prop. 4.3] we get that $\scrG_0^{\ad}$ is of $B_{n_1}+B_{n-n_1-1}$ Lie type, where $n_1\in\{0,\ldots,n-1\}$ (the $B_0$ Lie type corresponds to the trivial group).

Let $\scrF_0$ be the simple factor of $H^{\ad}_{A,\overline{\dbQ_p}}$ that is naturally isomorphic to $\scrH_0^{\ad}$. Let $\scrF_1$ be the image of $G_{\overline{\dbQ_p}}$ in $\scrF_0$. We identify naturally $\scrF_1=\scrG_0^{\ad}$. There exists an $H_{A,\overline{\dbQ_p}}$-submodule $\scrW_1$ of $W_A\otimes_{\dbQ} \overline{\dbQ_p}$ such that the $Z^0(H_{A,\overline{\dbQ_p}})$-modules $\scrW_0$ and $\scrW_1$ are isomorphic and the  representation of $\Lie(\scrF_0)$ on $\scrW_0\oplus\scrW_1$ is the spin representation, cf. Proposition 5.10 (b). The $H_{A,\overline{\dbQ_p}}$-modules $\scrW_0$ and $\scrW_1$ are not isomorphic. The representations of $\Lie(\scrF_1)$ on $\scrW_0$ and $\scrW_1$ are isomorphic to the tensor product of the spin representations of the two simple factors of $\Lie(\scrF_1)$, cf. [Pi, Prop. 4.3]. As $Z^0(G_{\overline{\dbQ_p}})=Z^0(H_{A,\overline{\dbQ_p}})$ (cf. Theorem 1.3.1) and as the representation of $\Lie(H^{\ad}_{A,\overline{\dbQ_p}})$ on $\scrW_0\oplus\scrW_1$ factors through $\Lie(\scrF_0)$ (cf. Proposition 5.10 (a)), we get that the $G_{\overline{\dbQ_p}}$-modules $\scrW_0$ and $\scrW_1$ are isomorphic. Thus the $H_{A,\overline{\dbQ_p}}$-modules $\scrW_0$ and $\scrW_1$ are isomorphic, cf. Theorem 1.1.4. Contradiction. Therefore $\scrG_0=\scrH_0$. Thus the property (*) holds.\endproof

\bigskip\noindent
{\bf 7.2. Theorem.} {\it We assume that the following two conditions hold:

\medskip
{\bf (i)} there exists no element $j\in J$ such that $(H_j,X_j)$ is of $D_4^{\dbH}$ type;

\smallskip
{\bf (ii)} for each element $j\in J$ the group $G_{\overline{\dbQ_p}}$ surjects onto a simple factor of $H_{j,\overline{\dbQ_p}}$.

\medskip
Then we have $G^{\ad}_{\dbQ_p}=H^{\ad}_{A,\dbQ_p}$.}

\medskip
\proof
We choose $v$ such that: it is unramified over a prime $l\in\dbN$ different from $p$ and greater than the number $N(A)$ of Theorem 1.3.2 and moreover the ranks of $T_v$ and $G_{\dbQ_p}$ are equal. From (ii) and Theorem 1.3.2 we get that the image of $G_{\overline{\dbQ_p}}$ in an arbitrary simple factor $\scrF_0$ of $H_{A,\overline{\dbQ_p}}^{\ad}$ has the same rank as $\scrF_0$. Let $\tilde H_{A,\overline{\dbQ_p}}$ be the subgroup of $H_{A,\overline{\dbQ_p}}$ that contains $Z^0(H_{A,\overline{\dbQ_p}})$ and whose image in $H_{A,\overline{\dbQ_p}}^{\ad}$ is the product of the images of $G_{\overline{\dbQ_p}}$ in the simply factors of $H_{A,\overline{\dbQ_p}}^{\ad}$. We have $G_{\overline{\dbQ_p}}\leqslant\tilde H_{A,\overline{\dbQ_p}}\leqslant H_{A,\overline{\dbQ_p}}$ and $\tilde H_{A,\overline{\dbQ_p}}$ has the same rank as $H_{A,\overline{\dbQ_p}}$. The centralizers of $\tilde H_{A,\overline{\dbQ_p}}$ and $H_{A,\overline{\dbQ_p}}$ in $\text{\bf GL}_{W_A\otimes_{\dbQ} \overline{\dbQ_p}}$ coincide, cf. Theorem 1.1.4. From the last two sentences we get that $\Lie(\tilde H_{A,\overline{\dbQ_p}})=\Lie(H_{A,\overline{\dbQ_p}})$, cf. Lemma 1.1.7. Thus $\tilde H_{A,\overline{\dbQ_p}}=H_{A,\overline{\dbQ_p}}$ and therefore $G_{\overline{\dbQ_p}}$ surjects onto $\scrF_0$. This implies that:

\medskip
{\bf (iii)} For each simple factor $\scrF_0$ of $H_{A,\overline{\dbQ_p}}^{\ad}$, there exists a unique simple factor of $G^{\ad}_{\overline{\dbQ_p}}$ that surjects (isomorphically) onto $\scrF_0$.

\medskip
We now start to check that we have an identity $G^{\ad}_{\dbQ_p}=H^{\ad}_{A,\dbQ_p}$. Based on the property 3.1 (c) and on (iii), the identity $G^{\ad}_{\dbQ_p}=H^{\ad}_{A,\dbQ_p}$ is equivalent to the fact that no simple factor of $G^{\ad}_{\overline{\dbQ_p}}$ surjects onto two or more simple factors of $H^{\ad}_{A,\overline{\dbQ_p}}$. To check this last statement, we can assume that (cf. Subsubsections 5.4.1 to 5.4.3 and Subsection 5.2):

\medskip
{\bf (iv)} all simple factors of either $H_{A,\overline{\dbQ}}^{\ad}$ or $G^{\ad}_{\overline{\dbQ_p}}$ have the same Lie type $\grL$,

\smallskip
{\bf (v)} and $f_A$ factors through an injective map
$$f_A^{(2)}:(H_A^{(2)},X_A^{(2)})\to (\pmb{\text{GSp}}(W_A,\psi_A),S_A)$$
of Shimura pairs for which the properties 5.10 (a) and (b) hold, with $(H_A^{(2),\ad},X_A^{(2),\ad})=(H_A^{\ad},X_A^{\ad})$ and with $f_A^{(2)}$ as a Hodge quasi product $\times^{\scrH}_{j\in J} f_j^{(2)}$ of injective maps
$$f_j^{(2)}:(H_j^{(2)},X_j^{(2)})\hookrightarrow (\text{\bf GSp}(W_j,\psi_j),S_j)$$ of Shimura pairs;
here $j\in J$ and $(H_j^{(2)},X_j^{(2)})$
is such that we can identify $(H_j^{(2),\ad},X_j^{(2),\ad})=(H_j,X_j)$.

\medskip
 Let $\scrF_1$ be a simple factor of $G^{\ad}_{\overline{\dbQ_p}}$.

\medskip\noindent
{\bf 7.2.1. Example.} We first treat the following particular case. We assume that all simple factors of $H_{A,\dbR}^{\ad}$ are non-compact and that all simple factors of $(H_A^{\ad},X_A^{\ad})$ are of $D_n^{\dbR}$ type. Thus we can also assume that $Z^0(H_A)=Z(\text{\bf GL}_W)\arrowsim\dbG_{m,\dbQ}$, cf. Propositions 5.10 (c) and 5.4.1. We assume that $\scrF_1$ surjects onto two distinct simple factors $\scrF_0^\prime$ and $\scrF_0^{\prime\prime}$ of $H^{\ad}_{A,\overline{\dbQ_p}}$. From Proposition 5.10 (a) and (b) we get that there exist four non-isomorphic, simple $\Lie(H^{\ad}_{A,\overline{\dbQ_p}})$-submodules $\scrW_0^\prime$, $\scrW_0^{\prime\prime}$, $\scrW_1^\prime$, and $\scrW_1^{\prime\prime}$ of $W_A\otimes_{\dbQ} \overline{\dbQ_p}$ with the property that for $*\in\{\prime,\prime\prime\}$, the $\Lie(\scrF_0^*)$-modules $\scrW_0^*$ and $\scrW_1^*$ are simple and associated to the two half spin representations of $\Lie(\scrF_0^*)$. We can assume that these four $\Lie(H^{\ad}_{A,\overline{\dbQ_p}})$-submodules are indexed in such a way that the $\Lie(\scrF_1)$-module $\scrW_0^\prime$ and $\scrW_1^\prime$ are isomorphic. As $Z^0(G_{\overline{\dbQ_p}})=Z^0(H_{A,\overline{\dbQ_p}})=Z(\text{\bf GL}_{W_A\otimes_{\dbQ} \overline{\dbQ_p}})$, $\scrW_0^\prime$ and $\scrW_1^\prime$ are also isomorphic as $Z^0(G_{\overline{\dbQ_p}})$-modules. Thus the simple $G_{\overline{\dbQ_p}}$-modules $\scrW_0^\prime$ and $\scrW_1^\prime$ are isomorphic. Therefore the centralizers of $G_{\overline{\dbQ_p}}$ and $H_{A,\overline{\dbQ_p}}$ in $\text{\bf GL}_{W_A\otimes_{\dbQ}\overline{\dbQ_p}}$ are distinct and this contradicts Theorem 1.1.4. Thus $G^{\ad}_{\dbQ_p}=H^{\ad}_{A,\dbQ_p}$. We come back to the general case.

\medskip\noindent
{\bf 7.2.2. General construction.} In order to show that no simple factor of $G^{\ad}_{\overline{\dbQ_p}}$ surjects onto two or more simple factors of $H^{\ad}_{A,\overline{\dbQ_p}}$, we need to have a very good understanding of the injective map $f_A^{(2)}$ and of the action of the torus $Z^0(H_A^{(2)})$ on $W_A$. Thus next we will list some basic properties of $f_A^{(2)}$. Based on Variant 4.9 (b) we can assume that $f_A^{(2)}$ is such that there exists a totally real number subfield $F$ of $\overline{\dbQ}$ with the property that each $f_j^{(2)}$ factors through the composite of two injective maps
$$(H_j^{(3)},X_j^{(3)})\hookrightarrow (H_j^{(4)},X_j^{(4)})\operatornamewithlimits{\hookrightarrow}\limits^{f_j^{(4)}} (\text{\bf GSp}(W_j,\psi_j),S_j)$$
of Shimura pairs for which the following three properties hold.

\medskip
{\bf (i)} The group $H_j^{(3),\ad}$ is the $\Res_{F/\dbQ}$ of an absolutely simple, adjoint group over $F$ whose extension to $\overline{\dbQ}$ is of $\grL$ Lie type.  Moreover the monomorphism $H_j^{(2)}\hookrightarrow H_j^{(3)}$ gives birth to an injective map $(H_j^{(2),\ad},X_j^{(2),\ad})\hookrightarrow (H_j^{(3),\ad},X_j^{(3),\ad})$ of Shimura pairs which is obtained as in (1) of Subsubsection 2.2.3.

\smallskip
{\bf (ii)} The group $H_j^{(4),\ad}$ is the $\Res_{F/\dbQ}$ of an absolutely simple, adjoint group over $F$ whose extension to $\overline{\dbQ}$ is of $\grL_4$ Lie type. Here $\grL_4$ is $\grL$ if $\grL=A_n$, is $A_{2n-1}$ if $(H_j,X_j)$ is of $C_n$ or $D_n^{\dbH}$ type, and is $A_{2^n-1}$ if $(H_j,X_j)$ is of $B_n$ or $D_n^{\dbR}$ type. If $\grL\neq A_n$, then the monomorphism $H_{j,\overline{\dbQ}}^{(3),\der}\hookrightarrow H_{j,\overline{\dbQ}}^{(4),\der}$ is a product of $[F:\dbQ]$ monomorphisms associated naturally to homomorphisms of the form $h_{\overline{\dbQ}}$, where $h$ is as in the property 4.3.1 (c). If $\grL=A_n$ with $n\ge 2$ and if $(H_j,X_j)$ is without involution (thus $(H_j^{(3),\ad},X_j^{(3),\ad})$ and $(H_j^{(4),\ad},X_j^{(4),\ad})$ are also without involution), then the monomorphism $H_{j,\overline{\dbQ}}^{(3),\der}\hookrightarrow H_{j,\overline{\dbQ}}^{(4),\der}$ is a diagonal embedding (cf. Subsubsection 4.1.2). If either $\grL=A_n$ and $(H_j,X_j)$ is with involution (thus $(H_j^{(3),\ad},X_j^{(3),\ad})$ and $(H_j^{(4),\ad},X_j^{(4),\ad})$ are also with involution) or $\grL=A_1$, then the monomorphism $H_{j,\overline{\dbQ}}^{(3),\der}\hookrightarrow H_{j,\overline{\dbQ}}^{(4),\der}$ is an isomorphism.

\smallskip
{\bf (iii)} The map $f_j^{(4)}$ is a PEL type embedding constructed as in Proposition 4.1.

\medskip
The map $f_A^{(2)}$ factors through $f_A^{(4)}:=\times^{\scrH}_{j\in J} f_j^{(4)}:(H_A^{(4)},X_A^{(4)})\hookrightarrow (\text{\bf GSp}(W_A,\psi_A),S_A)$.
From properties (i), (ii), and 7.2 (iii) we easily get that:

\medskip
{\bf (iv)} Two simple factors of either $\Lie(H_{A,\overline{\dbQ_p}}^{\ad})=\Lie(H_{A,\overline{\dbQ_p}}^{(2),\ad})$ or $\oplus_{j\in J} \Lie(H_{j,\overline{\dbQ_p}}^{(3),\ad})$ do not map injectively into the same simple factor of $\Lie(H_{A,\overline{\dbQ_p}}^{(4),\ad})$.

\smallskip
{\bf (v)} Two simple factors of $\Lie(G_{\overline{\dbQ_p}}^{\ad})$ do not map injectively into the same simple factor of $\Lie(H_{A,\overline{\dbQ_p}}^{(4),\ad})$. Moreover all monomorphisms from a simple factor of $\Lie(G_{\overline{\dbQ_p}}^{\ad})$ into a simple factor $\Lie(H_{A,\overline{\dbQ_p}}^{(3),\ad})$, are isomorphisms.

\medskip\noindent
{\bf 7.2.3. Notations.} Let $s$ and $s_{\dbR}$ be the number of simple factors of $H_A^{(4),\ad}$ and $H_{A,\dbR}^{(4),\ad}$ (respectively). Thus $s$ is the number of elements of the set $J$ and $s_{\dbR}=[F:\dbQ]s$.

Let $\scrS$ be the set of extremal vertices of the Dynkin diagram of either $\Lie(H_{A,\overline{\dbQ}}^{(4),\ad})$ or $\Lie(H_{A,\overline{\dbQ_p}}^{(4),\ad})$. The group $\Gal(\dbQ)$ acts on $\scrS$ (see Subsection 2.1). Until Subsubsection 7.2.6 we will assume that $\grL\neq A_1$; thus $\grL_4\neq A_1$. As $\grL_4\neq A_1$, we can identify $\scrS$ with the $\Gal(\dbQ)$-set $\Hom_{\dbQ}(K_\scrS,\overline{\dbQ})$, for $K_{\scrS}$ a suitable product indexed by $j\in J$ of totally imaginary quadratic extensions of $F$. By replacing $F$ with a totally real finite field extension of it, we can assume that there exists a totally imaginary quadratic extension $K$ of $F$ such that we have $K_{\scrS}=K^s$ (cf. Variant 4.9 (b)). The $\dbQ$--algebra $K_{\scrS}$ acts faithfully on $W_A$ (see the proof of Proposition 4.1). Thus we can identify $\Res_{K_{\scrS}/\dbQ} \dbG_{m,K_{\scrS}}$ with a torus of $\text{\bf GL}_{W_A}$. All monomorphisms $\dbS\hookrightarrow H_{A,\dbR}^{(4)}$ that are elements of $X_A^{(4)}$, factor through the extension to $\dbR$ of the subgroup of $\text{\bf GL}_{W_A}$ generated by $H_A^{(4),\der}$, by $Z(\text{\bf GL}_{W_A})$, and by the maximal subtorus $T_c$ of $\Res_{K_{\scrS}/\dbQ} \dbG_{m,K_{\scrS}}$ which over $\dbR$ is compact (cf. the proofs of Proposition 4.1 and Theorem 4.8 and the definition of Hodge quasi products).

\medskip\noindent
{\bf 7.2.4. Proposition.} {\it  We recall that $\grL\neq A_1$ and that $\scrF_1$ is a fixed simple factor of $G_{\overline{\dbQ_p}}^{\ad}$. Let $W_A\otimes_{\dbQ} \overline{\dbQ_p}=\oplus_{u\in \chi_A} \scrV_u$ be the minimal direct sum decomposition such that the torus $Z^0(H_A)_{\overline{\dbQ_p}}$ acts on $\scrV_u$ via the character $u$ of it. Let $\scrN_0$ be the unique normal, semisimple subgroup of $H_{A,\overline{\dbQ_p}}^{(2),\der}$ whose adjoint is the product of those simple factors of $H_{A,\overline{\dbQ_p}}^{\ad}=H_{A,\overline{\dbQ_p}}^{(2),\ad}$ onto which $\scrF_1$ surjects. Then there exists $u_0\in \chi_A$ such that $\scrN_0$ acts faithfully on $\scrV_{u_0}$. If $\grL=A_n$, then all irreducible $\Lie(\scrF_1)$-submodules of $\scrV_{u_0}$ are isomorphic.}

\medskip
\proof
We recall $f_A^{(4)}$ is a PEL type embedding, cf. property 7.2.2 (iii) and Fact 2.4.1 (b). Thus the Frobenius torus $T_v$ is naturally identified with a torus of a particular type of inner form of $H_A^{(4)}$, cf. the choice of $v$ and Subsubsections 5.6.1 to 5.6.3 applied to $(A,v)$ (instead of $(B,v)$ of Subsection 5.5). We will use the torus $T_v$ to show that the torus $Z^0(H_A)$ is small enough so that the Proposition holds; thus based on Fact 5.6.2, to ease notations we will assume that this inner form of $H_A^{(4)}$ is trivial and therefore that $T_v$ is naturally a torus of $H_A^{(4)}$. As the arguments are long, we will group them in six parts.

\smallskip
{\bf Part I. The torus $T_v^{(4)}$.} We check that the centralizer $T^{(4)}_v$ of $T_v$ in $H_A^{(4)}$ is a maximal torus of $H_A^{(4)}$. It is well known that $T^{(4)}_v$ is a reductive subgroup of $H_A^{(4)}$ (cf. [Bor, Ch. IV, 13.17, Cor. 2 (a)]) of the same rank as $H_A^{(4)}$. Thus to check that $T^{(4)}_v$ is a torus, it suffices to check the following statement: the intersection of $T_{v,\overline{\dbQ_p}}^{(4)}$ with each factor of $H_{A,\overline{\dbQ_p}}^{(4),\der}$ that has a simple adjoint, is a torus. All factors of $G^{\ad}_{\overline{\dbQ_p}}$ are of $\grL$ type and an $H_A^{(4)}(\overline{\dbQ_p})$-conjugate of $T_{v,\overline{\dbQ_p}}$ is a maximal torus of $G_{\overline{\dbQ_p}}$. Thus, if $\grL=A_n$, then $\grL=\grL_4$ and the statement can be easily checked over $\overline{\dbQ_p}$. If $\grL\neq A_n$, then based on properties 7.2.2 (i) and (ii) and on the property 7.2 (iii), the statement follows from Fact 4.3.2. Thus $T^{(4)}_v$ is a maximal torus of $H_A^{(4)}$.

\smallskip
{\bf Part II. The equivalence relation $\scrR$.} Let $\scrR$  be the unique equivalence relation on $\scrS$ that has the following property. Two distinct elements $\grn_1$ and $\grn_2$ of $\scrS$ are in relation $\scrR$ if and only if the following two conditions hold:

\medskip
{\bf (a)} {\it they are vertices of Dynkin diagrams of the Lie algebras of two distinct simple factors $\scrF_0^\prime$ and $\scrF_0^{\prime\prime}$ of $H_{A,\overline{\dbQ_p}}^{(4),\ad}$ and there exists a simple factor $\scrF_1^\prime$ of $G^{\ad}_{\overline{\dbQ_p}}$ whose Lie algebra maps injectively into both $\Lie(\scrF_0^\prime)$ and $\Lie(\scrF_0^{\prime\prime})$;}

\smallskip
{\bf (b)} {\it with respect to $\scrF_1^\prime$ and to the image of $T_{v,\overline{\dbQ_p}}^{(4)}$ in $H_{A,\overline{\dbQ_p}}^{(4),\ad}$, they are both either original vertices or ending vertices of Dynkin diagrams.}

\medskip
We explain the meaning in (b) of  ``are both either original vertices or ending vertices". As $T_{v,\overline{\dbQ_p}}^{(4)}$ is a maximal torus of $H_A^{(4)}$, there exist isomorphisms $\scrF_0^\prime\arrowsim\scrF_0^{\prime\prime}$ under which the images in $\scrF_0^\prime$ and $\scrF_0^{\prime\prime}$ of an a priori fixed maximal torus of the semisimple subgroup of $G_{\overline{\dbQ_p}}$ whose adjoint is $\scrF_1^\prime$, are canonically identified. Each such isomorphism is uniquely determined up to an automorphism of $\scrF_0^\prime$ that fixes a maximal torus of it and thus up to an inner automorphism of $\scrF_0^\prime$. Thus $\scrF_0^\prime$ and $\scrF_0^{\prime\prime}$ are naturally identified up to inner automorphisms and therefore it makes sense to speak in (b) about ``are both either original vertices or ending vertices". Based on this and the first part of the property 7.2.2 (v), properties (a) and (b) define uniquely an equivalence relation on $\scrS$. Each equivalence class of $\scrR$ contains at most one vertex of each Dynkin diagram of a simple factor of $\Lie(H^{(4),\ad}_{\overline{\dbQ_p}})$, cf. the very definition of $\scrR$.

\smallskip
{\bf Part III. The Galois invariance of $\scrR$.} In this Part III we check that the relation $\scrR$ is $\Gal(\dbQ)$-invariant. In (a) and (b) (of Part II), the roles of $G_{\overline{\dbQ_p}}$ and $T_{v,\overline{\dbQ_p}}^{(4)}$ can be also played by an $H_A^{(4)}(\overline{\dbQ_p})$-conjugate of $G_{\overline{\dbQ_p}}$ and by an arbitrary maximal torus of this $H_A^{(4)}(\overline{\dbQ_p})$-conjugate. But there exists an $H_A^{(4)}(\overline{\dbQ_p})$-conjugate of $G_{\overline{\dbQ_p}}$ that contains the extension to $\overline{\dbQ_p}$ of the subtorus of $T_v$ of $H_A^{(4)}$. Thus not to introduce extra notations, in the next paragraph we can assume that the monomorphism $T_v\hookrightarrow H_A^{(4)}$ has the property that its extension to $\overline{\dbQ_p}$ factors as a composite monomorphism $T_{v,\overline{\dbQ_p}}\hookrightarrow G_{\overline{\dbQ_p}}\hookrightarrow H_{A,\overline{\dbQ_p}}^{(4)}$.

Let ${\gr n}_1^\prime$ and ${\gr n}_2^{\prime}$ (resp. ${\gr n}_1^{\prime\prime}$ and ${\gr n}_2^{\prime\prime}$) be the two vertices of $\scrS$ that correspond to $\scrF_0^\prime$ (resp. to $\scrF_0^{\prime\prime}$). We can assume that $\grn_1^\prime$ and $\grn_1^{\prime\prime}$ are original vertices and that $\grn_2^{\prime}$ and $\grn_2^{\prime\prime}$ are ending vertices in the sense of Part II. The equivalence relation $\scrR$ on $\scrS$ is the smallest equivalence relation on $\scrS$ that contains all pairs of the form $(\grn_i^\prime,\grn_i^{\prime\prime})$ with $i\in\{1,2\}$. Let $T_{v,1}$ be the unique subtorus of $T_{v,\overline{\dbQ_p}}$ with the property that $\Lie(T_{v,1})$ is a Cartan Lie subalgebra of $\Lie(\scrF_1^\prime)$. Let $T_{v,1}^\perp$ be the unique subtorus of $T_{v,\overline{\dbQ_p}}$ with the property that the natural homomorphism $T_{v,1}\times_{\overline{\dbQ_p}} T_{v,1}^\perp\to T_{v,\overline{\dbQ_p}}$ is an isogeny. The torus $T_{v,1}^\perp$ has a trivial image in $\scrF_1^\prime$ and thus the image of $\Lie(T_{v,1}^\perp)$ in $\Lie(\scrF_0^\prime)\oplus\Lie(\scrF_0^{\prime\prime})$ is trivial. The Galois group $\Gal(\dbQ)$ acts naturally on the set of subtori of $T_{v,\overline{\dbQ_p}}$ and on the set of simple factors of $H_{A,\overline{\dbQ_p}}^{(4)}$. Let $\gamma\in\Gal(\dbQ)$. The torus $\gamma(T_{v,1})$ (resp. $\gamma(T_{v,1}^\perp)$) is a subtorus of $T_{v,\overline{\dbQ_p}}$ which has non-trivial images (resp. has trivial images) in each one of the two simple factors $\Lie(\gamma(\scrF_0^\prime))$ and $\Lie(\gamma(\scrF_0^{\prime\prime}))$ of $\Lie(H_{A,\overline{\dbQ_p}}^{(4),\ad})$. Let $\scrF^\prime_{1,1}$ and $\scrF^\prime_{1,2}$ be the unique simple factors of $G^{\ad}_{\overline{\dbQ_p}}$ whose Lie algebras map injectively into $\Lie(\gamma(\scrF_0^\prime))$ and $\Lie(\gamma(\scrF_0^{\prime\prime}))$ (respectively), cf. property 7.2.2 (v). We show that the assumption that $\scrF^\prime_{1,1}\neq\scrF^\prime_{1,2}$ leads to a contradiction. The torus $\gamma(T_{v,1}^{\perp})$ has a trivial image in $\scrF^\prime_{1,1}\times_{\overline{\dbQ_p}}\scrF^\prime_{1,2}$. Thus the rank of $G_{\dbQ_p}$ (i.e., of either $T_v$ or $T_{v,1}\times_{\overline{\dbQ_p}} T_{v,1}^\perp$) is at least equal to the rank of $T_{v,1}^{\perp}\times_{\overline{\dbQ_p}}\scrF^\prime_{1,1}\times_{\overline{\dbQ_p}}\scrF^\prime_{1,2}$. Therefore the rank of $T_{v,1}$ is at least equal to the rank of $\scrF^\prime_{1,1}\times_{\overline{\dbQ_p}}\scrF^\prime_{1,2}$. But as the groups $T_{v,1}$, $\scrF^\prime_{1,1}$, and $\scrF^\prime_{1,2}$ have the same ranks as $\grL$, we reached a contradiction. Thus $\scrF^\prime_{1,1}=\scrF^\prime_{1,2}$. Therefore we can define $\gamma(\scrF_1^\prime):=\scrF^\prime_{1,1}=\scrF^\prime_{1,2}$; it is a simple factor of $G^{\ad}_{\overline{\dbQ_p}}$ whose Lie algebra maps injectively into both $\Lie(\gamma(\scrF_0^\prime))$ and $\Lie(\gamma(\scrF_0^{\prime\prime}))$. As the tori $T_v$ and $T_v^{(4)}$ are defined over $\dbQ$ and as $T_{v,1}$ is defined over $\overline{\dbQ}$, it is easy to see that the vertices $\gamma(\grn_1^\prime)$ and $\gamma(\grn_1^{\prime\prime})$ are both original vertices and that the vertices $\gamma(\grn_2^\prime)$ and $\gamma(\grn_2^{\prime\prime})$ are both ending vertices of $\scrS$. Thus we have $(\gamma(\grn_1^\prime),\gamma(\grn_1^{\prime\prime}))$, $(\gamma(\grn_2^\prime),\gamma(\grn_2^{\prime\prime}))\in\scrR$. This implies that the equivalence relation $\scrR$ is $\Gal(\dbQ)$-invariant.

The Galois group $\Gal(\dbQ)$ acts on the quotient set $\scrS_1:=\scrS/\scrR$. Let $K_{\scrS_1}$ be the \'etale $\dbQ$--algebra such that we can identify $\scrS_1$ with the $\Gal(\dbQ)$-set $\Hom_{\dbQ}(K_{\scrS_1},\overline{\dbQ})$. Let $T_{c,1}$ be the maximal subtorus of $\Res_{K_{\scrS_1}/\dbQ} \dbG_{m,K_{\scrS_1}}$ which over $\dbR$ is compact. We identify $K_{\scrS_1}$ with a $\dbQ$--subalgebra of $K_{\scrS}$ and thus $K_{\scrS_1}$ is a product of at most $s$ fields. We have $T_{c,1}\leqslant T_c$.
Let $H_A^{(5)}$ be the reductive subgroup of $H_A^{(4)}$ generated by $H_A^{(4),\der}$, $Z(\text{\bf GL}_{W_A})$, and $T_{c,1}$.

\smallskip
{\bf Part IV. Cocharacters.} We write $H^{(4),\der}_{A,\dbR}=\prod_{i=1}^{s_{\dbR}} F^{(4)}_i$ as a product of normal, semisimple subgroups that have absolutely simple adjoints. Let $W_A\otimes_{\dbQ} \dbR=\oplus_{i=1}^{s_{\dbR}} V^{(4)}_i$ be the direct sum decomposition into $H^{(4)}_{A,\dbR}$-modules such that for all $i\in\{1,\ldots,s_{\dbR}\}$ the group $F^{(4)}_i$ acts trivially on $V^{(4)}_{i_1}$ if $i_1\in\{1,\ldots,s_{\dbR}\}\setminus\{i\}$ and acts non-trivially on each irreducible $H_{A,\dbR}^{(4)}$-submodule of $V^{(4)}_i$. Let $T^{(4)}_i$ be the direct factor of $T_{c,\dbR}$ that is isomorphic to $\text{\bf SO}(2)_{\dbR}=\dbS/\dbG_{m,\dbR}$ and that acts trivially on $\oplus_{i_1\in\{1,\ldots,s_{\dbR}\}\setminus\{i\}}V^{(4)}_{i_1}$. We identify $F^{(4)}_i$ and $T^{(4)}_i$ with their images in $\text{\bf GL}_{V^{(4)}_i}$. Let $x_4\in X^{(2)}_A\subseteq X^{(4)}_A$. Let $J^{(4)}_i$ be the image of the restriction of $x_4:\dbS\to H^{(4)}_{A,\dbR}$ to the compact subtorus $\text{\bf SO}(2)_{\dbR}$ of $\dbS$ in $T^{(4)}_i/(F^{(4)}_i\cap T^{(4)}_i)$; it is equipped with a homomorphism $j_i:\text{\bf SO}(2)_{\dbR}\to J^{(4)}_i$. Let $\mu_A^{(3),\ad}:\dbG_{m,\dbC}\to H^{(3),\ad}_{A,\dbC}$ and $\mu_A^{(4),\ad}:\dbG_{m,\dbC}\to H^{(4),\ad}_{A,\dbC}$ be the cocharacters defined naturally by $\mu_A:\dbG_{m,\dbC}\to H_{A,\dbC}$.

For $q\in\{l,p\}$ we consider an embedding $i_{\overline{\dbQ_q}}:\overline{\dbQ_q}\hookrightarrow\dbC$ that extends the embedding $i_E$ of Subsection 1.1. We use $i_{\overline{\dbQ_p}}$ to identity the sets of simple factors of $H_{A,\dbC}^{(4),\ad}$ and $H_{A,\overline{\dbQ_p}}^{(4),\ad}$. By applying the property 3.4.1 (b), we get the existence of a cocharacter $\mu^{\acute et}_{v,\overline{\dbQ_l}}:\dbG_{m,\overline{\dbQ_l}}\to G_{\overline{\dbQ_l}}$ whose extension via $i_{\overline{\dbQ_l}}$ is $H^{(4)}_A(\dbC)$-conjugate to $\mu_A$. Based on Subsubsection 5.6.4, we can assume that $\mu^{\acute et}_{v,\overline{\dbQ_l}}$ factors through an $H_A^{(4)}(\overline{\dbQ_l})$-conjugate of $T_{v,\overline{\dbQ_l}}$. As an $H_A^{(4)}(\overline{\dbQ_p})$-conjugate of $T_{v,\overline{\dbQ_p}}$ is a maximal torus of $G_{\overline{\dbQ_p}}$, we conclude that there exists a cocharacter of $G_{\overline{\dbQ_p}}$ whose extension via $i_{\overline{\dbQ_p}}$ is $H^{(4)}_A(\dbC)$-conjugate to $\mu_A$. From this and the second part of the property 7.2.2 (v) we get that:

\medskip
{\bf (c.i)} {\it two simple factors of the group pair $(H_{A,\dbC}^{(3),\ad},[\mu_A^{(3),\ad}])$ that correspond to two simple factors of $\Lie(H_{A,\overline{\dbQ_p}}^{(3),\ad})$ into which the same simple factor of $\Lie(G^{\ad}_{\overline{\dbQ_p}})$ injects, are isomorphic.}

\medskip
From the property (c.i) and the description of the monomorphisms $H_{j,\overline{\dbQ}}^{(3),\der}\hookrightarrow H_{j,\overline{\dbQ}}^{(4),\der}$ in the property 7.2.2 (ii), we get that:

\medskip
\medskip
{\bf (c.ii)} {\it two simple factors of the group pair $(H_{A,\dbC}^{(4),\ad},[\mu_A^{(4),\ad}])$ that correspond to two simple factors of $\Lie(H_{A,\overline{\dbQ_p}}^{(4),\ad})$ into which the same simple factor of $\Lie(G^{\ad}_{\overline{\dbQ_p}})$ injects, are isomorphic.}

\medskip
In properties (c.i) and (c.ii) we can replace $(H_{A,\dbC}^{(\bigstar),\ad},[\mu_A^{(\bigstar),\ad}])$ by $(H_{A,\dbR}^{(\bigstar),\ad},[\mu_A^{(\bigstar),\ad}])$ (cf. [De2, Prop. 1.2.2]), where $\bigstar\in\{3,4\}$.

\smallskip
{\bf Part V. Downsizing the torus $Z^0(H_A)$.} We show that we can modify $x_4$ (and thus implicitly $X_A^{(2)}$ and $X_A^{(4)}$) such that its images in $X_A^{\ad}=X_A^{(2),\ad}$ and $X_A^{(4),\ad}$ are the same and the composite homomorphism
$$(j_1,\ldots,j_{s_{\dbR}}):\text{\bf SO}(2)_{\dbR}\to \prod_{i=1}^{s_{\dbR}} J^{(4)}_i\hookrightarrow\prod_{i=1}^{s_{\dbR}} T^{(4)}_i/(F^{(4)}_i\cap T^{(4)}_i)=T_{c,\dbR}/(T_{c,\dbR}[2])$$
factors through $T_{c,1,\dbR}/(T_{c,1,\dbR}[2])$. Such a modification of $x_4$ is allowed (cf. Proposition 5.4.1 and the constructions of Subsection 5.2) and it is defined naturally by a replacement of the $s_{\dbR}$-tuple $(j_1,\ldots,j_{s_{\dbR}})$ with another $s_{\dbR}$-tuple $(j_1^{\eps_1},\ldots,j_{s_{\dbR}}^{\eps_{s_{\dbR}}})$, where $\eps_1,\ldots,\eps_{s_{\dbR}}\in\{-1,1\}$ and where we can have $\eps_i=-1$ only if $F^{(4)}_i$ is compact (cf. the paragraph after properties (a) and (b) of the proof of Proposition 4.1). The homomorphism $(j_1,\ldots,j_{s_{\dbR}})$ factors through $T_{c,1,\dbR}/(T_{c,1,\dbR}[2])$ if and only if the following property holds:

\medskip
{\bf (d)} {\it for arbitrary distinct elements $i^\prime$, $i^{\prime\prime}\in\{1,\ldots,s_{\dbR}\}$ that correspond to factors $\scrF_0^\prime$ and $\scrF_0^{\prime\prime}$ as in (a) of Part II, the representations of $\text{\bf SO}(2)_{\dbR}$ on $V^{(4)}_{i^\prime}$ and $V^{(4)}_{i^{\prime\prime}}$  defined by $j_{i^\prime}$ and $j_{i^{\prime\prime}}$ (respectively)  are such that we have identifications
$$V^{(4)}_{i^\prime}\otimes_{\dbR} \dbC=V^{(4)+}_{i^\prime}\oplus V^{(4)-}_{i^\prime}=V^{(4)+}_{i^{\prime\prime}}\oplus V^{(4)-}_{i^{\prime\prime}}=V^{(4)}_{i^{\prime\prime}}\otimes_{\dbR} \dbC$$
of $\text{\bf SO}(2)_{\dbC}=\dbG_{m,\dbC}$-modules, that are defined naturally by (2) and that have the property that the representations of $F^{(4)}_{i^\prime,\dbC}$ and $F^{(4)}_{i^{\prime\prime},\dbC}$ on $V^{(4)+}_{i^\prime}$ and  $V^{(4)+}_{i^{\prime\prime}}$ (respectively) are both via original vertices of the Dynkin diagrams of the Lie algebras of $F^{(4)}_{i^\prime,\dbC}$ and $F^{(4)}_{i^{\prime\prime},\dbC}$ (here original vertices are in the sense of (b) of Part II, the role of $\overline{\dbQ_p}$ being replaced by the one of $\dbC$).}

\medskip
Based on the property (c.ii) the groups $F^{(4)}_{i^\prime}$ and $F^{(4)}_{i^{\prime\prime}}$ are either both compact or both non-compact. It is easy to see that we can choose the $\eps_i$'s such that (d) holds for all pairs $(i^\prime,i^{\prime\prime})$ such that both groups $F^{(4)}_{i^\prime}$ and $F^{(4)}_{i^{\prime\prime}}$ are compact. If $F^{(4)}_{i^\prime}$ and $F^{(4)}_{i^{\prime\prime}}$ are both non-compact, then from the property (c.ii) and the uniqueness of $\eps_{i^\prime}$ and $\eps_{i^{\prime\prime}}$ (they are both $1$) we get that (d) automatically holds for the pair $(i^\prime,i^{\prime\prime})$. Therefore we can assume that (d) holds and thus that the homomorphism $(j_1,\ldots,j_{s_{\dbR}})$ factors through $T_{c,1,\dbR}/(T_{c,1,\dbR}[2])$. We get that $\mu_A$ factors through the reductive group $H^{(5)}_{A,\dbC}$ defined at the end of Part III. As $H_A$ is a Mumford--Tate group, we get that $Z^0(H_A)\leqslant Z^0(H_A^{(5)})$.

\smallskip
{\bf Part VI. End of the proof of 7.2.4.} Let $W_A\otimes_{\dbQ} \overline{\dbQ_p}=\oplus_{t\in \chi_A^{(5)}} \scrV_t$ be the minimal direct sum decomposition such that the torus $Z^0(H_{A,\overline{\dbQ_p}}^{(5)})$ acts on each $\scrV_t$ via its character $t$. See proof of Proposition 4.1 for how the torus $\Res_{K_{\scrS}/\dbQ} \dbG_{m,K_{\scrS}}$ (and thus also how its subtorus $T_{c,1}$) acts on $W_A$. The direct sum decompositions $V^{(4)}_i\otimes_{\dbR} \dbC=V^{(4)+}_i\oplus V^{(4)-}_i$ (see (d)) are such that ${\Res_{K_{\scrS}/\dbQ} \dbG_{m,K_{\scrS}}}_{\dbC}$ acts on each $V^{(4),u}_i$ via a unique character, where $i\in\{1,\ldots,s_{\dbR}\}$ and $u\in\{-,+\}$ (see (2) of the proof of Proposition 4.1). The upper indices $+$ and $-$ are such that $T_{c,1,\dbC}$ acts on $V^{(4),u^\prime}_{i^\prime}$ and $V^{(4),u^{\prime\prime}}_{i^{\prime\prime}}$ via the same character if and only if $u^\prime=u^{\prime\prime}\in\{-,+\}$ and the elements $i^\prime, i^{\prime\prime}\in\{1,\ldots,s_{\dbR}\}$ are as in (d). We get the existence of two distinct characters $t_1$, $t_2\in \chi_A^{(5)}$ such that $\scrN_0$ acts faithfully on both $\scrV_{t_1}$ and $\scrV_{t_2}$ and acts trivially on $\scrV_t$ if $t\in \chi_A^{(5)}\setminus\{t_1,t_2\}$. As $Z^0(H_A)$ is a subtorus of $Z^0(H_A^{(5)})$, each $\overline{\dbQ_p}$-vector space $\scrV_u$ (of Proposition 7.2.4) is a direct sum of some of the $\scrV_t$'s. Thus we can take the character $u_0\in \chi_A$ such that we have either $\scrV_{t_1}\subseteq \scrV_{u_0}$ or $\scrV_{t_2}\subseteq \scrV_{u_0}$.

We assume that $\grL=A_n$. Due to the vertices parts of (b) and (d) and to the identity $\grL=\grL_4$, the irreducible $\Lie(\scrF_1)$-submodules of $\scrV_{u_0}$ are all isomorphic.\endproof

\medskip\noindent
{\bf 7.2.5. End of the proof for $\grL\neq A_1$.} From properties 7.2.2 (iii) and (iv) we get that the representation of $\scrN_0$ on $\scrV_{u_0}$ is a direct sum of irreducible representations of its factors that have absolutely simple adjoints (to be compared with the property 5.7 (a)). Let $\Lie(\scrF_{1,0})$ be a simple factor of $\Lie(\scrN_0)$. If $\grL$ is $B_n$ or $C_n$ with $n\Ge 2$, then from Proposition 5.10 (a) we get that the simple $\Lie(\scrF_{1,0})$-submodules of $\scrV_{u_0}$ are all isomorphic. If $\grL=D_4$, then each $(H_j,X_j)$ is of $D_4^{\dbR}$ type and thus from Proposition 5.10 (b) we get that both half spin representations of $\Lie(\scrF_{1,0})$ are among the simple $\Lie(\scrN_0)$-submodules of $\scrV_{u_0}$. Let now $\grL$ be $D_n$, with $n\Ge 5$. If $\scrN_0$ has a normal factor that is a $\text{\bf Spin}_{2n,\overline{\dbQ_p}}$ group, then from the property 7.2.4 (c.ii) we get that each simple, adjoint Shimura pair $(H_j,X_j)$ such that the intersection $\scrN_0^{\ad}\cap H_{j,\overline{\dbQ_p}}$ (taken inside $H^{\ad}_{A,\overline{\dbQ_p}}$) is non-trivial, is of $D_n^{\dbR}$ type. Thus $\scrN_0$ is a product of $\text{\bf Spin}_{2n,\overline{\dbQ_p}}$ groups (cf. property (iv) of Variant 4.9 (a)) and therefore again from Proposition 5.10 (b) we get that both half spin representations of $\Lie(\scrF_{1,0})$ are among the simple $\Lie(\scrN_0)$-submodules of $\scrV_{u_0}$. If $\scrN_0$ has a normal factor that is an $\text{\bf SO}_{2n,\overline{\dbQ_p}}$ group, then a similar argument shows that each simple, adjoint Shimura pair $(H_j,X_j)$ such that the intersection $\scrN_0^{\ad}\cap H_{j,\overline{\dbQ_p}}$ (taken inside $H^{\ad}_{A,\overline{\dbQ_p}}$) is non-trivial, is of $D_n^{\dbH}$ type, and that $\scrN_0$ is a product of $\text{\bf SO}_{2n,\overline{\dbQ_p}}$ groups. Thus the simple $\Lie(\scrF_{1,0})$-submodules of $\scrV_{u_0}$ are all isomorphic to the standard $2n$ dimensional simple $\Lie(\scrF_1)$-module.

Regardless of what $\grL\neq A_1$ is, based on the previous paragraph we conclude that the $\Lie(\scrF_1)$-module $\scrV_{u_0}$ is either isotypic or a direct sum of two non-isomorphic, isotypic $\Lie(\scrF_1)$-modules. The last situation happens if and only if $\scrN_0$ is a product of $\text{\bf Spin}_{2n,\overline{\dbQ_p}}$ groups.

Let $s_0$ be the number of simple factors of $\scrN^{\ad}_0$. The Lie monomorphism $\Lie(\scrF_1)\hookrightarrow\Lie(\scrN^{\ad}_0)$ is isomorphic to the diagonal embedding $\Lie(\scrF_{1,0})\hookrightarrow \Lie(\scrF_{1,0})^{s_0}$. If $\scrN_0$ is not (resp. is) a product of $\text{\bf Spin}_{2n,\overline{\dbQ_p}}$ groups, then the $\Lie(\scrN_0)$-module $\scrV_{u_0}$ is a direct sum of $2s_0$ (resp. $4s_0$) non-isomorphic, isotypic $\Lie(\scrN_0)$-modules; thus the $H_{A,\overline{\dbQ_p}}$-module $\scrV_{u_0}$ is also a direct sum of $2s_0$ (resp. $4s_0$) non-isomorphic, isotypic $H_{A,\overline{\dbQ_p}}$-modules. But the $G_{\overline{\dbQ_p}}$-module $\scrV_{u_0}$ is a direct sum of $2$ (resp. $4$) non-isomorphic, isotypic $G_{\overline{\dbQ_p}}$-modules, cf. the very definition of $\scrN_0$ in Proposition 7.2.4. Based on this and Theorem 1.1.4, we get that $s_0=1$. Thus the simple factor $\scrF_1$ of $G_{\overline{\dbQ_p}}^{\ad}$ does not surject onto two or more simple factors of $H^{\ad}_{A,\overline{\dbQ_p}}$. This proves the Theorem 7.2 for $\grL\neq A_1$.

\medskip\noindent
{\bf 7.2.6. End of the proof for $\grL=A_1$.} If $\grL=A_1$, then the only modifications required to be made in Subsubsections 7.2.3 to 7.2.5 are as follows. For $j\in J$ we have $H_j^{(3)}=H_j^{(4)}$ and we have to replace the factor of $K_{\scrS}$ that corresponds to $j$ by a totally imaginary quadratic extension $K_j$ of $F$ obtained as in the case $n=1$ of the proof of Proposition 4.1; by replacing $F$ with a totally real finite field extension of it (cf. Variant 4.9 (b)), we can assume that there exists a totally imaginary quadratic extension $K_0$ of $\dbQ$ such that we have $K_j=K:=F\otimes_{\dbQ} K_0$ for all $j\in J$. We can identify $\scrS$ with the $\Gal(\dbQ)$-set $\Hom_{\dbQ}(F^s,\overline{\dbQ})$. We define $K_{\scrS}:=K^s$. The equivalence class $\scrR$ on $\scrS$ is defined this time by the property 7.2.4 (a) alone. If $\scrS_1:=\scrS/\scrR$, let $F_{\scrS_1}$ be the $\dbQ$--subalgebra of $F^s$ such that we can identify $\scrS_1$ with the $\Gal(\dbQ)$-set $\Hom_{\dbQ}(F_{\scrS_1},\overline{\dbQ})$. Let $K_{\scrS_1}:=F_{\scrS_1}\otimes_{\dbQ} K_0$. Let $T_{c,1}$ be as in the proof of Proposition 7.2.4. The rest of the arguments do not have to be modified. This ends the proof of Theorem 7.2 for $\grL=A_1$ and thus also for all $\grL$.\endproof

\bigskip\noindent
{\bf 7.3. Proofs of 1.3.4 and 1.3.6.} If for all elements $j\in J$ the type of $(H_j,X_j)$ is among those listed in Theorem 1.3.4, then the property 7.2 (i) obviously holds and Theorem 7.1 implies that the property 7.2 (ii) also holds. Thus $G^{\ad}_{\dbQ_p}=H_{A,\dbQ_p}^{\ad}$, cf. Theorem 7.2. Therefore $G_{\dbQ_p}=H_{A,\dbQ_p}$, cf. Theorem 1.3.1. This ends the proof of Theorem 1.3.4.

To prove Theorem 1.3.6, we can assume that the condition 7.2 (ii) holds (cf. Proposition 5.4.1 and Subsubsection 5.4.3). The hypothesis of Theorem 1.3.6 implies that the conditions 7.2 (i) also holds. Thus Theorem 1.3.6 follows from Theorem 7.2.\endproof

\bigskip\noindent
{\bf 7.4. Corollary.} {\it Let $J_4:=\{j\in J|(H_j,X_j)\,\text{is of either}\,D_4^{\dbH}\,\text{or}\,D_4^{\dbR}\,\text{type}\}$. If $j\in J_4$ (resp. if $j\in J\setminus J_4$) we assume that $(H_j,X_j)$ is of non-inner $D_4^{\dbH}$ type (resp. is of one of the types listed in Theorem 1.3.4). We also assume that there exists a prime $q\in\dbN$ such that the following two properties hold:

\medskip
{\bf (i)} the simple factors of $\prod_{j\in J_4} H_{j,\dbQ_q}$ are absolutely simple and pairwise non-isomorphic;

\smallskip
{\bf (ii)} for all elements $j\in J_4$ the Shimura pair $(H_j,X_j)$ is with $\dbQ_q$-involution.

\medskip
Then Conjecture 1.1.1 holds for $A$.}

\medskip
\proof
The proof of the Corollary is very much the same as the (d') part of the proof of Theorem 7.1 and the proof of Theorem 7.2. We can assume that the set $J_4$ is non-empty, cf. Theorem 1.3.4. Based on Theorem 1.1.6 we can assume that $p=q$. As in the (d') part of the proof of Theorem 7.1, based on (ii) and the absolutely simple part of (i) we get that for each element $j\in J_4$ the group $G_{\dbQ_p}$ surjects onto a simple factor of $H_{j,\dbQ_p}$.

Similar to the arguments that checked the property 7.2 (iii), based on Theorem 1.1.4 and Lemma 1.1.7 we argue that $G_{\overline{\dbQ_p}}$ surjects onto each simple factor of $\prod_{j\in J_4} H_{j,\overline{\dbQ_p}}$. Thus $G_{\dbQ_p}$ surjects onto $\prod_{j\in J_4} H_{j,\dbQ_p}$, cf. the pairwise non-isomorphic part of (i). Therefore in the proof of Theorem 7.2 we can work only with a Lie type $\grL\neq D_4$ and thus Theorem 7.2 applies (i.e., we have $G^{\ad}_{\dbQ_p}=H_{A,\dbQ_p}^{\ad}$). Therefore $G_{\dbQ_p}=H_{A,\dbQ_p}$, cf. Theorem 1.3.1. \endproof

\bigskip\noindent
{\bf 7.5. The simplest cases not settled by 1.3.4.}
For the sake of future references we list below the ``simplest"  cases  not settled by Theorem 1.3.4 (they involve adjoint groups whose ranks are small and of fixed parity). We assume that $H_A^{\ad}$ is a simple group over $\dbQ$ but we do not impose any restriction on either $p$ or $G_{\dbQ_p}$. The image $\grI$ of $G_{\overline{\dbQ_p}}$ in a fixed simple factor of $H_{A,\overline{\dbQ_p}}^{\ad}$ is an adjoint group, cf. property 3.1 (c); let $\grL$ be its Lie type. In A5 to A9 below, the other possible values for $\grL$ allowed by [Pi, Table 4.2] are excluded as in Example 4 of Subsubsection 6.2.4.  If $(H_A^{\ad},X_A^{\ad})$ is of $A_n$ type, then let $(H_0,X_0):=(H_A^{\ad},X_A^{\ad})$ and let the set $\grC$ be as in Subsection 6.2.

\medskip\noindent
{\bf A5 and A7.} The Shimura pair $(H_A^{\ad},X_A^{\ad})$ is of $A_5$ (resp. $A_7$) type and $\grC=\{1,{1\over 2}\}$ (resp. $\grC$ is $\{1\}$ or $\{1,{1\over 3}\}$). We could theoretically have $\grL=A_1+A_2$ (resp. $\grL=A_1+A_3$).

\smallskip\noindent
{\bf A8.} The Shimura pair $(H_A^{\ad},X_A^{\ad})$ is of $A_8$ type and $\grC=\{{1\over 2}\}$. We could theoretically have $\grL=A_2+A_2$.

\smallskip\noindent
{\bf A9.} The Shimura pair $(H_A^{\ad},X_A^{\ad})$ is of  $A_9$ type. We could theoretically have either $\grL=A_4$ and $\grC=\{{2\over 3}\}$ or $\grL=A_1+A_4$ and $\grC\in\{\{1,{1\over 4}\},\{1,{2\over 3}\}\}$.

\smallskip\noindent
{\bf A11.} The Shimura pair $(H_A^{\ad},X_A^{\ad})$ is of  $A_{11}$ type. We could theoretically have $\grL\in\{A_1+A_5,A_1+A_1+A_2,A_2+A_3,A_2+C_2\}$. The Lie types $A_1+A_1+A_2$ and $A_2+C_2$ have the same rank; if $\grL$ is one of these two types, then $\grC=\{1,{1\over 2}\}$.

\smallskip\noindent
{\bf C4 and C6.} The Shimura pair $(H_A^{\ad},X_A^{\ad})$ is of $C_4$ (resp. $C_6$) type. We could theoretically have $\grL=A_1+A_1+A_1$ (resp. $\grL=A_1+D_3$).

\smallskip\noindent
{\bf D4 and D6.} The Shimura pair $(H_A^{\ad},X_A^{\ad})$ is of $D_4^{\dbH}$ type and Corollary 7.4 does not apply (resp. is of inner $D_6^{\dbH}$ type). We could theoretically have $\grL=B_3$ or $\grL=B_1+B_2$ (resp. $\grL=C_1+C_3$).

\smallskip\noindent
{\bf D35.} The Shimura pair $(H_A^{\ad},X_A^{\ad})$ is of $D_{35}^{\dbH}$ type. We could theoretically have $\grL=A_7$.

\smallskip
The $A_n$ types with $n\in\{1,2,3,4,6,10,12\}$ are settled by Example 6 of Subsubsection 6.2.4.

\medskip\noindent
{\bf 7.5.1. Remarks.}
{\bf (a)} The numerical tests of Subsection 6.2 ``work" regardless of what $q$ is. However, due to the fact that in Definition 6.1 (c) we are allowed to choose the prime $q$, there exist large classes of simple, adjoint Shimura pairs of non-special $A_n$ type for which the numerical tests of Subsection 6.2 fail. We include one example related to A8 of Subsection 7.5. We assume that the group $H_A^{\ad}$ is simple and  $H_{A,\dbQ_p}^{\ad}$ has a simple factor $\scrC$ that is the adjoint group of the group of invertible elements of a central division algebra over $\dbQ_p$ of dimension $81$. If $\grI$ is a subgroup of $\scrC_{\overline{\dbQ_p}}$, then we have $\grL\neq A_2+A_2$ as otherwise $\grI$ (and therefore also $\scrC$) splits over a finite field extension of $\dbQ_p$ of degree a divisor of $12$ and thus not divisible by $9$. Thus $\grL=A_8$ and therefore $(H_A^{\ad},X_A^{\ad})$ is of non-special $A_8$ type. Similar such examples can be constructed in connection to A7 of Subsection 7.5. This remark and Subsections 6.2 and 7.5 motivate the terminology ``non-special" used in Definition 6.1 (c).

\smallskip
{\bf (b)} Let $(H_0,X_0)$ be a simple, adjoint, unitary Shimura pair. If by chance $F(H_0,X_0)\otimes_{\dbQ} \dbQ_q$ is an unramified field extension of $\dbQ_q$, then one can re-obtain easily most of the numerical properties of Subsection 6.2 without being ``bothered" to mention $T_1$ in Definition 6.1 (b) or to appeal to Proposition 5.10. But this does not suffice to get them all (like Example 4 of Subsubsection 6.2.4) or to get Theorem 1.3.4 (a). Thus even if we have a good understanding of the (assumed to be very simple) arithmetics of the field $F(H_0,X_0)$, Theorem 7.2 (and thus implicitly Proposition 5.10) looks to us unavoidable i.e., the proof of Theorem 7.2 can not be simplified. However, in practice it is often easy to check that the numbers $t_j$ of Lemma 6.2.2 do not depend on $j\in J_1$, and accordingly, that the pairs $(r_{j,m},s_{j,m})$ (counted with multiplicities) and defined when $f(j)(m)=1$, do not depend on $j\in J_1$. As in A9 of Subsection 7.5, we might not have a unique possibility for some of $t_j$'s.

\medskip
In what follows the abelian motives over a number subfield of $\dbC$ are defined using absolute Hodge cycles, see [DM, Sect. 6]. Using realizations in the Betti and \'etale cohomologies, as for abelian varieties we speak about the Mumford--Tate groups and the attached (adjoint) Shimura pairs of their extensions to $\dbC$ and about the Mumford--Tate conjecture for them. Theorem 1.3.7 is a particular case of the following more general result.

\bigskip\noindent
{\bf 7.6. Theorem (the independence property for abelian motives).}
{\it We recall that in Subsection 1.1 we fixed an embedding $i_E:E\hookrightarrow\dbC$. Let $(\scrD_i)_{i\in I}$ be a finite set of abelian motives over the number field $E$. We assume that for each $i\in I$, the Mumford--Tate conjecture is true for $\scrD_i$ and no simple factor of the adjoint Shimura pair attached to the extension of $\scrD_i$ to $\dbC$ is of $D_4^{\dbH}$ type. Then the Mumford--Tate conjecture is true for $\scrD:=\oplus_{i\in I} \scrD_i$.}

\medskip
\proof
We take the abelian variety $A$ over $E$ such that $\scrD$ is an object of the Tannakian category of abelian motives over $E$ generated by the $H_1$-motive $W_A=L_A\otimes_{\dbZ} \dbQ$ of $A$. By passing to a finite field extension of $E$, we can choose $A$ such that it has a principal polarization $\lambda_A$ (cf. [Mu2, Ch. IV, Sect. 23, Cor. 1]). Let $h_A:\dbS\to H_{A,\dbR}$ be as in Subsection 1.3. Let $f_A:(H_A,X_A)\hookrightarrow (\text{\bf GSp}(W_A,\psi_A),S_A)$ be as in Subsection 1.5. Let $H_{\scrD}$ be the reductive group over $\dbQ$ that is the image of $H_A$ in $\text{\bf GL}_{W_{\scrD}}$, where $W_{\scrD}$ is the rational vector space that is the Betti realization of $\scrD$. We have to show that the identity component $G_{\scrD,\dbQ_p}$ of the algebraic envelope of the image of the $p$-adic Galois representation attached to $\scrD$ is $H_{\scrD,\dbQ_p}$. Always $G_{\scrD,\dbQ_p}$ is a subgroup of $H_{\scrD,\dbQ_p}$, cf. Theorem 1.1.2. Based on this and Theorem 1.3.1, we only have to show that $H_{\scrD,\dbQ_p}^{\ad}=G_{\scrD,\dbQ_p}^{\ad}$.

Let $H_B^{\der}$ and $H_C^{\der}$ be the commuting normal subgroups of $H_A^{\der}$ such that their adjoints are naturally identified with $H_{\scrD}^{\ad}$ and $H_A^{\ad}/H^{\ad}_{\scrD}$ (respectively). Let $H^\prime_B$ be a reductive subgroup of $H_A$ of the same rank as $H_A$ such that (cf. [Ha, Lem. 5.5.3]) the $\dbR$-ranks of $H_B^{\prime,\ab}$ and $H_A^{\ab}$ are both $1$ and we have $H^{\prime,\der}_B=H_B^{\der}$. Two maximal compact tori of $H^{\der}_C$ are $H^{\der}_C(\dbR)$-conjugate. This implies that there exists a homomorphism $h_B:\dbS\to H^\prime_{B,\dbR}$ such that the following two properties hold: (i) its composite with the inclusion $i_A^B:H^\prime_{B,\dbR}\hookrightarrow H_{A,\dbR}$ is $H_A(\dbR)$-conjugate to $h_A$, and (ii) the homomorphisms $\dbS\to H^{\ad}_{B,\dbR}=H^{\ad}_{\scrD,\dbR}$ defined naturally by $h_B$ and $h_A$ coincide. We choose $h_B$ such that the smallest subgroup $H_B$ of $H_B^\prime$ with the property that $h_B$ factors through $H_{B,\dbR}$, is a reductive group that has $H_B^{\der}$ as its derived group. By replacing $E$ with a finite field extension of it, we can chose $h_B$ such that moreover the polarized abelian variety over $\dbC$ defined (see Riemann's theorem of [De1, Thm. 4.7]) by the triple $(L_A,\eps\psi_A,h_B)$ is the pull back of a polarized abelian variety $(B,\lambda_B)$ over $E$; here $\eps\in\{-1,1\}$ depends on the connected component of $S_A$ to which $f_{A,\dbR}\circ i^B_{A,\dbR}\circ h_B$ belongs.

Let $C$ be the abelian variety over $E$ obtained similarly to $B$ but starting from $H_C^{\der}$ (here we might have to replace once more $E$ by a finite field extension of it). The abelian varieties $A$ and $B\times_E C$ are adjoint-isogenous. For each simple factor $H_j$ of $H_B^{\ad}=H_B^{\prime\ad}=H^{\ad}_{\scrD}$, there is $i\in I$ such that $H_j$ is isomorphic to a simple factor of the adjoint group $H_{\scrD_i}^{\ad}$ defined similarly to $H_{\scrD}^{\ad}$ but for $\scrD_i$ instead of $\scrD$. As the Mumford--Tate conjecture holds for $\scrD_i$, the identity component of the algebraic envelope of the image of the $p$-adic Galois representation attached to $B$ surjects naturally onto $H_{j,\dbQ_p}$. Thus this identity component surjects also onto $H_{B,\dbQ_p}^{\ad}$, cf. Theorem 7.2. From this and Theorem 1.3.1 we get that the Mumford--Tate conjecture holds for $B$. Thus we have $H_{B,\dbQ_p}^{\ad}=H_{\scrD,\dbQ_p}^{\ad}=G_{\scrD,\dbQ_p}^{\ad}$. \endproof

\medskip\noindent
{\bf 7.6.1. Example.}
This example is only meant to justify the name ``independence property" for Theorems 1.3.7 and 7.6. We assume that $\End(A_{\overline{\dbQ}})=\dbZ$; thus the adjoint Shimura pair $(H_A^{\ad},X_A^{\ad})$ is simple (see [Pi, Prop. 5.12]). We also assume that $(H_A^{\ad},X_A^{\ad})$ is of $D_{2n+2}^{\dbR}$ type and that the group $H_{A,\dbR}^{\ad}$ is not absolutely simple. Let $f_2:(G_2,X_2)\hookrightarrow (\text{\bf GSp}(W_1,\psi_1),S_1)$ be an injective map as in Subsections 5.2 and 5.3 (thus $(G_2^{\ad},X_2^{\ad})=(H_A^{\ad},X_A^{\ad})$). We assume that there exists an abelian variety $B$ over $E$ that is the pull back of an abelian scheme $\scrA_2$ as in Subsubsection 5.4.2. Let $C:=A\times_E B$. We also assume that $B$ is such that the adjoint group of the Mumford--Tate group $H_C$ of $C_{\dbC}$ is $H_A^{\ad}\times_{\dbQ} G_2^{\ad}=H_A^{\ad}\times_{\dbQ} H_A^{\ad}$. The identity component of the algebraic envelope of the image of the $p$-adic Galois representation attached to $C$ could theoretically be (without contradicting the results of Subsections 1.1.2 to 1.1.9 and of Section 3) a reductive group over $\dbQ_p$ whose adjoint is isomorphic to $H^{\ad}_{A,\dbQ_p}$ (and thus it is different from $H_{C,\dbQ_p}^{\ad}$). But Theorem 1.3.4 (e) implies that this identity component is $H_{C,\dbQ_p}$.

\bigskip\noindent
{\bf 7.7. Concluding remarks. (a)}
Theorems 1.3.4 and 7.6 extend (this is standard) to abelian varieties over finitely generated fields of characteristic 0.

\smallskip
{\bf (b)} We assume that $\End(A_{\overline{E}})=\dbZ$ and that the adjoint Shimura pair $(H_A^{\ad},X_A^{\ad})$ is simple of $D_{2n+2}^{\dbR}$ type. Property 4.8 (vi) and the shifting process of Subsection 5.4 allow another proof of [Pi, Sect. 7, pp. 230--236] i.e., allow us to get directly the lifting property of $l$-adic Galois representations mentioned in [Pi, Sect. 7, top of p. 230] and due to Wintenberger.

\smallskip
{\bf (c)} One can not get a general form of [Pi, Thm. 5.14] in the pattern of Theorem 1.3.4 (c) by using only the ideas of Sections 1 to 7. In z future paper we will present the general principle behind loc. cit. that extends significantly loc. cit. and this paper in the contexts when for each $j\in J$ the group $H_{j,\dbR}$ has only one simple, non-compact factor.

\smallskip
{\bf (d)} For the left cases of the Mumford--Tate conjecture, one can assume that all simple factors of $H_{j,\overline{\dbQ}}$ have the same Lie type $\grL$ and (if it helps) that properties 5.10 (a), (b), or (c) holds (cf. Proposition 5.4.1, Subsubsection 5.4.3, and Theorem 7.2). If $\grL\neq D_4$, then we can also assume that the group $H_A^{\ad}$ is $\dbQ$--simple (cf. Theorem 7.2).

\smallskip
{\bf (e)} Corollary 7.4 can be used to strengthen both Theorem 7.6 and (d) for the case when $\grL=D_4$.

\smallskip
{\bf (f)} Many ideas of this paper can be used for other classes of motives, provided (some of) the analogues of [De2, Subsubsects. 2.3.9 to 2.3.13], Theorem 1.1.4, [Sa1], and [Zi, Thm. 4.4] are known to be true in some form. In connection to Conjecture 1.1.1, the results [De2, Subsubsects. 2.3.9 to 2.3.13] and [Zi, Thm. 4.4] were neither used nor quoted prior to (earlier versions of) this work.

\smallskip
{\bf (g)} By combining [An, Thm. 1.5.1 and Cor. 1.5.2] with Theorem 1.3.4, one gets a different proof of the result [An, Thm. 1.6.1 3)] that pertains to the Mumford--Tate conjecture for hyperk\"ahler varieties over number fields.

\smallskip
{\bf (h)} We have a variant of the notion of non-special $A_n$ type (see Definition 6.1 (c)) in which the order of quantifiers is changed: there exists an $f_1$ such for all $\tilde f_1^\prime$ and $T_1$ as in Definition 6.1 (b), there exist two distinct primes $q_1$, $q_2\in\dbN$ with the property that for all $MT$ pairs $(T_1,G_1)$ and $(T_2,G_2)$ for $(f_1,\tilde f_1^\prime,q_1)$ and $(f_1,\tilde f_1^\prime,q_2)$ (respectively), both groups $G_{1,\overline{\dbQ_{q_1}}}^{\ad}$ and $G_{2,\overline{\dbQ_{q_2}}}^{\ad}$ have simple factors of $A_n$ Lie type. Presently we do not know if such a variant could lead to significant improvements to Theorem 1.3.4 (a).

\bigskip\smallskip
\references{37}
{\nspace{

\Ref[An]
Y. Andr\'e,
\sl On the Shafarevich and Tate conjectures for hyperk\"ahler varieties,
\rm Math. Ann. {\bf 305} (1996), no. 2, pp. 205--248.

\Ref[Bl]
D. Blasius,
\sl A p-adic property of Hodge cycles on abelian varieties,
\rm Motives (Seattle, WA, 1991),  pp. 293--308, Proc. Symp. Pure Math., Vol. {\bf 55}, Part 2, Amer. Math. Soc., Providence, RI, 1994.

\Ref[Bog]
F. A. Bogomolov,
\sl Sur l'alg\'ebricit\'e des repr\'esentations $l$-adiques,
\rm C. R. Acad. Sci. Paris S\'er. A-B Math. {\bf 290} (1980), no. 15, pp. 701--703.

\Ref[Bor]
A. Borel,
\sl Linear algebraic groups. Second enlarged edition,
\rm Grad. Texts in Math., Vol. {\bf 126}, Springer-Verlag, New York, 1991.

\Ref[Bou1]
N. Bourbaki,
\sl Groupes et alg\`ebres de Lie,
\rm Chapitre {\bf 4--6}, Actualit\'es Scientifiques et Industrielles, No. {\bf 1337}, Hermann, Paris, 1968.

\Ref[Bou2]
N. Bourbaki,
\sl Groupes et alg\`ebres de Lie,
\rm Chapitre {\bf 7--8}, Actualit\'es Scientifiques et Industrielles, No. {\bf 1364}, Hermann, Paris, 1975.

\Ref[BHC]
A. Borel and Harish-Chandra,
\sl Arithmetic subgroups of algebraic groups,
\rm Ann. of Math. (2) {\bf 75} (1962), no. 3, pp. 485--535.

\Ref[BLR]
S. Bosch, W. L\"utkebohmert, and M. Raynaud,
\sl N\'eron models,
\rm Ergebnisse der Mathematik und ihrer Grenzgebiete (3), Vol. {\bf 21}, Springer-Verlag, Berlin, 1990.

\Ref[BO]
P. Berthelot and A. Ogus,
\sl F-crystals and de Rham cohomology. I,
\rm Invent. Math. {\bf 72} (1983), no. 2, pp. 159--199.

\Ref[Ch]
W. C. Chi,
\sl $l$-adic and $\lambda$-adic representations associated to abelian varieties defined over number fields,
\rm Amer. J. Math. {\bf 114} (1992), no. 2, pp. 315--353.

\Ref[De1]
P. Deligne,
\sl Travaux de Shimura,
\rm S\'eminaire  Bourbaki, 23\`eme ann\'ee (1970/71), Exp. No. 389, pp. 123--163, Lecture Notes in Math., Vol. {\bf 244}, Springer-Verlag, Berlin, 1971.

\Ref[De2]
P. Deligne,
\sl Vari\'et\'es de Shimura: interpr\'etation modulaire, et
techniques de construction de mod\`eles canoniques,
\rm Automorphic forms, representations and $L$-functions (Oregon State Univ., Corvallis, OR, 1977), Part 2, pp. 247--289, Proc. Sympos. Pure Math., Vol. {\bf 33}, Amer. Math. Soc., Providence, RI, 1979.

\Ref[De3]
P. Deligne,
\sl Hodge cycles on abelian varieties,
\rm Hodge cycles, motives, and Shimura varieties, pp. 9--100, Lecture Notes in Math., Vol. {\bf 900}, Springer-Verlag, Berlin-New York, 1982.

\Ref[DM]
P. Deligne and J. Milne,
\sl Tannakian categories,
\rm Hodge cycles, motives, and Shimura varieties, pp. 101--128, Lecture Notes in Math., Vol. {\bf 900}, Springer-Verlag, Berlin-New York, 1982.

\Ref[Fa]
G. Faltings,
\sl Endlichkeitss\"atze f\"ur abelsche Variet\"aten \"uber Zahlk\"orpern,
\rm Invent. Math. {\bf 73} (1983), no. 3, pp. 349--366.

\Ref[Fo]
J.-M. Fontaine,
\sl Le corps des p\'eriodes p-adiques,
\rm P\'eriodes $p$-adiques (Bures-sur-Yvette, 1988), pp. 59--111, J. Ast\'erisque {\bf 223}, Soc. Math. de France, Paris, 1994.

\Ref[FH]
W. Fulton and J. Harris,
\sl Representation theory. A first course,
\rm Grad. Texts in Math., Vol. {\bf 129}, Springer-Verlag, New York, 1991.

\Ref[Ha]
G. Harder,
\sl \"Uber die Galoiskohomologie halbeinfacher Matrizengruppen. II,
\rm  Math. Z. {\bf 92} (1966), pp. 396--415.

\Ref[He]
S. Helgason,
\sl Differential geometry, Lie groups, and symmetric spaces,
\rm Pure and Applied Mathematics, Vol. {\bf 80}, Academic Press, Inc. [Harcourt Brace Jovanovich, Publishers], New York-London, 1978.

\Ref[Ko1]
R. E. Kottwitz,
\sl Stable trace formula: elliptic singular terms,
\rm Math. Ann.  {\bf 275}  (1986), no. 3, pp. 365--399.

\Ref[Ko2]
R. E. Kottwitz,
\sl Points on some Shimura Varieties over finite fields,
\rm J. of Amer. Math. Soc. {\bf 5} (1992), no. 2, pp. 373--444.

\Ref[La]
S. Lang,
\sl Algebraic number theory. Second edition,
\rm Grad. Texts in Math., Vol. {\bf 110}, Springer-Verlag, New York, 1994.

\Ref[LP]
M. Larsen and R. Pink,
\sl Abelian varieties, $l$-adic representations and $l$-independence,
\rm Math. Ann. {\bf 302} (1995), no. 3, pp. 561--579.

\Ref[Mi1]
J. S. Milne,
\sl The action of an automorphism of $\dbC$ on a Shimura
variety and its special points,
\rm Arithmetic and geometry, Vol. I, pp. 239--265, Progr. in Math., Vol. {\bf 35}, Birkh\"auser, Boston, MA, 1983.

\Ref[Mi2]
J. S. Milne,
\sl The points on a Shimura variety modulo a prime of good
reduction,
\rm The Zeta functions of Picard modular surfaces, pp. 151--253, Univ. Montr\'eal, Montreal, Quebec, 1992.

\Ref[Mi3]
J. S. Milne,
\sl Shimura varieties and motives,
\rm Motives (Seattle, WA, 1991),  pp. 447--523, Proc. Symp. Pure Math., Vol. {\bf 55}, Part 2, Amer. Math. Soc., Providence, RI, 1994.

\Ref[Mi4]
J. S. Milne,
\sl Descent for Shimura varieties,
\rm Mich. Math. J. {\bf 46} (1999), no. 1, pp. 203--208.

\Ref[Mu1]
D. Mumford,
\sl Families of abelian varieties,
\rm Algebraic Groups and Discontinuous Subgroups (Boulder, CO, 1965),  pp. 347--351, Proc. Sympos. Pure Math., Vol. {\bf 9}, Amer. Math. Soc., Providence, RI, 1966.

\Ref[Mu2]
D. Mumford,
\sl Abelian varieties,
\rm Tata Inst. of Fund. Research Studies in Math., No. {\bf 5}, Published for the Tata Institute of Fundamental Research, Bombay; Oxford Univ. Press, London, 1970 (reprinted 1988).

\Ref[No]
R. Noot,
\sl Classe de conjugaison du Frobenius des vari\'et\'es ab\'eliennes \`a r\'eduction ordinaire,
\rm Ann. Inst. Fourier (Grenoble)  {\bf 45}  (1995),  no. 5, pp. 1239--1248.

\Ref[Pi]
R. Pink,
\sl $l$-adic algebraic monodromy groups, cocharacters, and the Mumford--Tate conjecture,
\rm J. reine angew. Math. {\bf 495} (1998), pp. 187--237.

\Ref[Sa1]
I. Satake,
\sl Holomorphic imbeddings of symmetric domains into a Siegel space,
\rm  Amer. J. Math. {\bf 87} (1965), pp. 425--461.

\Ref[Sa2]
I. Satake,
\sl Symplectic representations of algebraic
groups satisfying a certain analyticity condition,
\rm Acta Math. {\bf 117} (1967), pp. 215--279.

\Ref[Se1]
J.-P. Serre,
\sl Groupes alg\'ebriques associ\'es aux modules de Hodge--Tate,
\rm Journ\'ees de G\'eom. Alg. de Rennes. (Rennes, 1978), Vol. III, pp. 155--188, J. Ast\'erisque {\bf 65}, Soc. Math. de France, Paris, 1979.

\Ref[Se2]
J.-P. Serre,
\sl Propri\'et\'es conjecturales des groupes de Galois motiviques et des repr\'esentations $l$-adiques,
\rm Motives (Seattle, WA, 1991), pp. 377--400, Proc. Sympos. Pure Math., Vol. {\bf 55}, Part 1, Amer. Math. Soc., Providence, RI, 1994.

\Ref[Se3]
J.-P. Serre,
\rm Collected papers, Vol. IV, 1985--1998, Springer-Verlag, Berlin, 2000.

\Ref[Sh]
G. Shimura,
\sl Moduli of abelian varieties and number theory,
\rm Algebraic Groups and Discontinuous Subgroups (Boulder, CO, 1965), pp. 312--332, Proc. Sympos. Pure Math., Vol. {\bf 9}, Amer. Math. Soc., Providence, RI, 1966.

\Ref[Ta1]
S. G. Tankeev,
\sl Cycles on abelian varieties and exceptional numbers,
\rm Izv. Ross. Akad. Nauk. Ser. Mat. {\bf 60} (1996), no. 2, pp. 159--194, translation in  Izv. Math. {\bf 60}  (1996),  no. 2, pp. 391--424.

\Ref[Ta2]
S. G. Tankeev,
\sl On the weights of an $l$-adic representation and the arithmetic of Frobenius eigenvalues,
\rm (Russian) Izv. Ross. Akad. Nauk
Ser. Mat. {\bf 63} (1999), no. 1, pp. 185--224, translation in Izv. Math. {\bf 63} (1999), no. 1, pp. 181--218.

\Ref[Tat]
J. Tate,
\sl Endomorphisms of abelian varieties over finite fields,
\rm Invent. Math. {\bf 2} (1966), pp. 134--144.

\Ref[Ts]
T. Tsuji,
\sl p-adic \'etale cohomology and crystalline cohomology
in the semi-stable reduction case,
\rm Invent. Math. {\bf 137} (1999), no. 2, pp. 233--411.

\Ref[Va1]
A. Vasiu,
\sl Integral canonical models of Shimura varieties of preabelian type,
\rm Asian J. Math. {\bf 3} (1999), no. 2, pp. 401--518.

\Ref[Va2]
A. Vasiu,
\sl Surjectivity criteria for $p$-adic representations, Part I,
\rm Manuscripta Math. {\bf 112} (2003), no. 3, pp. 325--355.

\Ref[Za1]
Y. G. Zarhin,
\sl Abelian varieties having a reduction of K3 type,
\rm  Duke Math. J. {\bf 65} (1992), no. 3, pp. 511--527.

\Ref[Za2]
Y. Zarhin,
\sl Very simple $2$-adic representations and hyperelliptic Jacobians,
\rm Moscow Math. J. {\bf 2} (2002), no. 2, pp. 403--431.

\Ref[Zi]
T. Zink,
\sl Isogenieklassen von Punkten von Shimuramannigfaltigkeiten mit Werten in einem endlichen K\"orper,
\rm Math. Nachr. {\bf 112} (1983), pp. 103--124.

\Ref[Wi]
J.-P. Wintenberger,
\sl Un scindage de la filtration de Hodge pour certaines vari\'et\'es algebriques sur les corps locaux,
\rm Ann. of Math. (2) {\bf 119} (1984), no. 3, pp. 511--548.

}}

\bigskip
\hbox{Adrian Vasiu,}
\hbox{Department of Mathematical Sciences, Binghamton University,}
\hbox{Binghamton, New York 13902-6000, U.S.A.}
\hbox{e-mail: adrian\@math.binghamton.edu}

\enddocument